\documentclass[a4paper,abstracton]{scrartcl}
\usepackage[natbibapa]{apacite}
\usepackage{amsfonts}
\usepackage{amssymb}
\usepackage{amsmath}
\usepackage{amsthm}
\usepackage[utf8]{inputenc}
\usepackage{booktabs}
\usepackage{xspace}
\usepackage{soul}	
\usepackage{xcolor}
\usepackage{graphicx}
\usepackage{color}
\usepackage{subfigure}
\usepackage{a4wide}
\usepackage[normalem]{ulem}
\usepackage{authblk}

\allowdisplaybreaks
\DeclareMathOperator*{\argmin}{arg\,min}

\DeclareMathOperator{\OWA}{OWA}
\DeclareMathOperator{\opt}{Opt}
\DeclareMathOperator{\infeas}{Inf}

\begin{document}
\newcommand{\X}{{\mathcal{X}}}
\newcommand{\cU}{{\mathcal{U}}}
\newcommand{\cW}{{\mathcal{W}}}
\newcommand{\cI}{{\mathcal{I}}}
\newcommand{\cC}{{\mathcal{C}}}
\newcommand{\U}{(UW)\xspace}
\newcommand{\UA}{(UA)\xspace}
\newcommand{\UO}{(UO)\xspace}
\newcommand{\Uniform}{Uniform $\pmb{w}$}
\newcommand{\UniformAlpha}{Uniform $\alpha$}
\newcommand{\Orness}{Uniform orness}
\newcommand{\RVR}{(RvR)\xspace}
\newcommand{\PVP}{(PvP)\xspace}
\newcommand{\OVP}{(OvP)\xspace}
\newcommand{\WDist}{preference vector distance}
\newcommand{\XInDist}{in-sample solution distance}
\newcommand{\XOutDist}{out-of-sample solution distance}
\newcommand{\New}[1]{\textcolor{blue}{#1}}
\newcommand{\Removed}[1]{\textcolor{red}{\renewcommand{\ULdepth}{0.0em}\st{#1}\renewcommand{\ULdepth}{.75ex}}}

\title{A Preference Elicitation Approach for the Ordered Weighted Averaging Criterion using Solution Choice Observations}

\author[1]{Werner Baak\thanks{Corresponding author. Email: werner.baak@uni-passau.de} }
\author[1]{Marc Goerigk}
\author[2]{Michael Hartisch}

\affil[1]{Business Decisions and Data Science, University of Passau,\authorcr Dr.-Hans-Kapfinger-Str. 30, 94032 Passau, Germany}

\affil[2]{Network and Data Science Management, University of Siegen,\authorcr Unteres Schlo{\ss} 3, 57072 Siegen, Germany}

\date{}

\maketitle

\begin{abstract}
Decisions under uncertainty or with multiple objectives usually require the decision maker to formulate a preference regarding risks or trade-offs. If this preference is known, the ordered weighted averaging (OWA) criterion can be applied to aggregate scenarios or objectives into a single function. Formulating this preference, however, can be challenging, as we need to make explicit what is usually only implicit knowledge. We explore an optimization-based method of preference elicitation to identify appropriate OWA weights. We follow a data-driven approach, assuming the existence of observations, where the decision maker has chosen the preferred solution, but otherwise remains passive during the elicitation process. We then use these observations to determine the underlying preference by finding the preference vector that is at minimum distance to the polyhedra of feasible vectors for each of the observations. Using our optimization-based model, weights are determined by solving an alternating sequence of linear programs and standard OWA problems. Numerical experiments on risk-averse preference vectors for selection, assignment and knapsack problems show that our passive elicitation method compares well against having to conduct pairwise comparisons and performs particularly well when there are inconsistencies in the decision maker's choices.\end{abstract}


\noindent\textbf{Acknowledgements:} Supported by the Deutsche Forschungsgemeinschaft (DFG) through grant GO 2069/2-1.

\section{Introduction}

Decision making is a ubiquitous challenge, where often multiple conflicting objectives or the consequences over multiple scenarios need to be taken into account \citep{ehrgott}. In such settings, a popular approach is to use an aggregation function to combine several values into a single one. 
The Ordered Weighted Average (OWA) operator is one such method \citep{yager1988ordered}. Since its inception, it has been widely studied and applied in settings as diverse as fuzzy modeling \citep{o1988aggregating,yager1998including}, location planning \citep{malczewski2006ordered}, financial decision-making problems \citep{merigo2009induced,merigo2011induced}, geographic information system based site planning \citep{zabihi2019gis} or risk assignment \citep{chang2011evaluating}, see also the survey by \cite{emrouznejad2014ordered}.

The idea of the OWA function is to take a vector of values as input, sort this vector from largest to smallest value, and to calculate the scalar product of this sorted vector with a weight vector. Hence, the weights are assigned to ordered values and can be used to stress importance on high, low or mid-ranged inputs. These weights should represent the decision maker's preferences.

Solving problems with the OWA operator is usually more challenging than solving their single-criterion counterparts (usually called nominal counterpart in robust optimization), where a single objective function is given. Different solution methods have been developed, including linear programming (LP) or mixed-integer linear programming (MILP) models, see, e.g., \cite{ogryczak2003solving} and \cite{ogryczak2012milp}. In \cite{chassein2015alternative}, a compact reformulation of the model based on linear programming duality was established. The complexity of discrete decision making problems with the OWA criterion has been studied as well, see \cite{kasperski2015combinatorial} and \cite{chassein2020approximating}, where approximation algorithms have been developed, or \cite{Galand20121540}, where exact algorithms for the spanning tree problem are considered.

For the purpose of such theoretical analysis, it is usually assumed that the preference weights are given by the decision maker. In practice however, these weights are not simply given, but must first be determined. For this purpose, several methods have been developed, see the surveys by \cite{xu2005overview} and \cite{liu2011review}.

The idea of sample learning methods is to use empirical data to fit OWA weights. Given a set of observations consisting of alternatives and their aggregated values, an optimization model is used to fit preference vectors that satisfy these observations as far as possible, see, e.g., \cite{yager1996quantifier} and \cite{garcia2011generating}. Using a similar idea, \cite{ahn2008preference} assumes that pairwise comparisons between solutions are given to indicate preferences and then uses an optimization model to find preference vectors that adhere to these comparisons. More methods are introduced by \cite{bourdache2017anytime,bourdache2019active} and \cite{bourdache2020bayesian}, which cover incremental elicitation or active learning strategies, where again knowledge of pairwise preferences is presupposed.

The knowledge of past choices in sets or queries on sets of solutions has been assumed in several previous works on elicitation, where the user or decision maker is (semi-)actively taking actions during the elicitation process. \cite{dragone2018constructive} propose a Choice Perceptron for learning user preferences, whereas \cite{zintgraf2018ordered} propose ordered preference elicitation strategies based on ranking and clustering. In \cite{viappiani2020equivalence} a recommender system is introduced, exploring the connection  between providing good recommendation and asking informative choice queries, in order to compute optimal recommendation as well as query sets.

Many more methods exist to elicit preference weights for OWA criteria. We briefly summarize some of these. \cite{benabbou2015minimax} use a search tree model for regret-based optimization, including OWA, where a minimization of max regret is performed. \cite{adam2021possibilistic} guarantee a robust approach identifying preferences and also an error detection method for wrong preferences in which OWA models are utilized to test against the true model. \cite{labreuche2015extension} use binary alternatives (ordinal information) as extension to the so-called MACBETH method to elicit an OWA operator and where the weights given are trapezoidal but can be weakened to a convex fuzzy set and thus taking account of inconsistencies. \cite{kim2018implicit} elicit the decision maker's preferences by comparing answers given with extreme or arbitrary options and based on the results adding constraints to the OWA weights. \cite{wang2007aggregating} make use of the orness degree to determine the weights by analyzing the decision maker's optimism level. Further preference elicitation methods include maximal entropy methods \citep{fuller2001analytic}, data-driven approaches \citep{filev1994learning}, introduction of a weight generating function \citep{filev1998issue,yager2016some} or using kernel density estimations \citep{lin2020determine}. 

Extensions to the OWA operator have been studied as well. These include
OWAWA (Ordered Weighted Average Weighted Average), WOWA (Weighted Ordered Weighted Average) and IOWA (Induced Ordered Weighted Average). The latter variant enables the possibility to reorder variables in a more complex way. Building on that, 
\cite{merigo2010fuzzy} proposes the induced generalized ordered weighted averaging (IGOWA) operator. The aggregation operator combines the characteristics of the generalized OWA and the induced OWA operator. The IEOWAD (Induced Euclidean Ordered Weighted Averaging Distance) approach by \cite{merigo2011induced}, parameterizes distances measures using an IOWA operator, resulting in a modality which allows considering more complex attitudinal characters of the decision maker and resulting in different conclusions.

A particularly important class of OWA operators makes use of or-like, non-increasing preference weights \citep{yager1993families}. These weights are chosen in such a way that they exhibit a decreasing pattern, assigning higher weights to lower-ranked objectives. This characteristic captures the risk-averse behavior of decision makers, who prioritize objectives that perform worse over those that perform well. This is in alignment with the Pigou-Dalton principle in social choice theory, which states that a more equitable distribution of utilities of any two persons, while keeping the sum of their utilities the same, is weakly preferred.
Similarly, in the context of decision-making under risk, using non-increasing weights reflects a preference for reducing risks or pessimism among decision makers.

The use of non-increasing weights in the OWA criterion offers several advantages, making it a valuable approach in decision-making. In particular, it contains the worst-case decision criterion as a special case, which is widely used in robust optimization \citep{aissi2009min,kasperski2016robust}. Here non-increasing weights align with the risk-averse attitude by assigning greater importance to the worst-performing objectives. This formulation ensures robustness against unfavorable outcomes and enhances the decision maker's ability to handle uncertainties. Another relevant example is the application of the conditional value at risk (CVaR) \citep{bertsimas2009constructing}. By considering the worst-case scenarios beyond a specified confidence level, decision makers can address potential extreme losses. The use of non-increasing weights in OWA aligns with the risk-averse nature of CVaR, as it assigns higher weights to objectives associated with the tail events, thereby capturing the decision maker's aversion to extreme losses. The adoption of non-increasing preference vectors in risk-averse decision-making is well-documented in the literature. Numerous studies have explored the properties and applications of such vectors, see, e.g., \citep{ogryczak2003minimizing, ogryczak2003solving,Galand20121540, kasperski2015combinatorial}. These investigations further validate the effectiveness and wide-ranging applicability of non-increasing weights in capturing risk-averse preferences. Furthermore, risk-averse weight vectors allow for improved problem reformulations, see \cite{chassein2015alternative,chassein2020approximating}.

In this paper we propose a model which enables  the determination of a risk-averse decision maker's preferences in the form of non-increasing OWA weights using a passive elicitation method. This method uses observed choices and does not require the decision maker to actively take (corrective) actions or interfere during the elicitation process. The used observations include information that was available to the decision maker, that is, for a given decision making problem, we only require the preferred solution and data basis, and in particular no aggregated values have to be assigned to alternatives.   
We propose a novel optimization-based model to determine OWA weights, which is solved by a sequence of MILPs. The model aims at finding a preference vector that is at minimum distance to the polyhedra of feasible OWA vectors for the observations. This approach utilizes historical observations, which are easily accessible in most applications today. Our proposed approach avoids the need for time-consuming and often costly interviews required to obtain pairwise comparisons or aggregated values of alternatives, making it an efficient and cost-effective decision support system. Not requiring interaction with the decision maker, it is less intrusive and more consumer-friendly, which can lead to increased user acceptance and adoption. This approach is particularly beneficial for decision-making situations where similar choices must be made repeatedly, and historical observations are available. Examples of applications for this approach include purchasing and selection decisions, routing problems and scheduling.
 We compare the properties of this approach to optimization models based on pairwise preference comparisons in numerical experiments. The evaluation of the approaches is carried out by measuring the average distance between the true OWA weights and the estimated weights.

While we focus on risk-averse weight vectors in the following, our models can also be used in the context of general weight vectors. The main difference is that there are less efficient models available to solve the resulting OWA problems.

The paper is organized in the following order: In Section~\ref{sec:def} we define the OWA problem and the notation. We present our optimization approach for preference elicitation and discuss modeling alternatives in Section~\ref{sec:opt}. In Section~\ref{sec:exp} we discuss the computational experiments, considering setup, results and insights. Finally, in Section~\ref{Conclusion} we reflect and conclude the paper taking into account possible further research approaches.

\section{OWA Problem Definition}\label{sec:def}

In this section we give a formal description of the OWA criterion and related optimization problems. Throughout this paper, we use the notation $[K]$ to denote a set $\{1,\ldots,K\}$, write vectors in bold, and drop the transpose symbol for vector multiplication if the context is clear.

We consider an optimization problem over some set of feasible solutions $\X\subseteq\mathbb{R}^n$ with a linear objective function $\pmb{c}\pmb{x}$ that we would like to minimize. We explore the case where there is more than one cost coefficient vector that is relevant for the decision making process. This can be because there are multiple relevant objectives, or due to uncertainty.

Here we assume that every criterion can be described via a linear function. This is particularly common in combinatorial optimization under uncertainty, see the survey by \cite{kasperski2016robust}, which describes several applications for this setting. This includes shortest path \citep{kasperski2015combinatorial}, spanning tree \citep{Galand20121540}, or knapsack \citep{chassein2020approximating} problems. The linearity assumption offers several benefits, including simplifying the analysis and mathematical formulation of the problem. Furthermore, linearity in a combinatorial problem does not necessarily require linear utility functions. For example, the problem with discrete alternatives, one of which should be chosen (see, e.g., the case of choosing a mobile phone discussed by \cite{reimann2017well}), corresponds to a combinatorial problem with a linear objective function and the constraint that one variable should be active, irrespective of how utilities are calculated. We also note that linearity is mainly assumed so that the arising OWA problems can be solved efficiently; for the definition of our model, this is not a strict requirement.

We denote by $\{\pmb{c}^1,\ldots,\pmb{c}^K\}$ the set of $K$ cost coefficient vectors that we would like to consider simultaneously and let $C\in \mathbb{R}^{K\times n}$ be the respective cost matrix. This means that for a given decision vector $\pmb{x}\in\X$, there are $K$ objective values $\pmb{c}^1\pmb{x},\ldots,\pmb{c}^K\pmb{x}$. Here it should be noted that these cost coefficients (or utilities) need to be proportionate to each other, i.e.~there should be a commensurate standardization. The purpose of the OWA criterion is to aggregate these objective values into a single value.

To this end, we require a preference vector $\pmb{w}\in[0,1]^K$ with $\sum_{k\in[K]} w_k = 1$. We write $\cW = \{ \pmb{w}\in[0,1]^K : \sum_{k\in[K]} w_k = 1\}$. The purpose of $w_k$ is to assign an importance to the $k$th-largest objective value for $k\in[K]$. Let $\pi$ be a permutation that sorts the $K$ objective values from largest to smallest, i.e., $\pi$ is such that $\pmb{c}^{\pi(1)}\pmb{x} \ge \pmb{c}^{\pi(2)}\pmb{x} \ge \ldots \ge \pmb{c}^{\pi(K)}\pmb{x}$. Then, the OWA operator is defined as
\[ \OWA_{\pmb{w}}(\pmb{x},C) = \sum_{k\in[K]} w_k (\pmb{c}^{\pi(k)} \pmb{x}) \]
Note that the permutation $\pi$ depends on the solution $\pmb{x}$.
The OWA operator contains several well-known decision making criteria as special cases. By setting $\pmb{w}=(1,0,\ldots,0)$, all weights are assigned to the largest objective value, which means that OWA becomes the worst-case criterion (in our minimization setting) as used in robust optimization. On the other hand, setting $\pmb{w}=(0,\ldots,0,1)$ gives the best-case criterion. Additionally, $w=(\alpha,0,\ldots,0,1-\alpha)$ for some $\alpha\in[0,1]$ corresponds to the Hurwicz criterion (see, e.g., \citep{chassein2020approximating}), while $w=(1/K,\ldots,1/K)$ gives the average value.

A special case of OWA operators uses risk-averse preference vectors. Risk-averse preference vectors play an important role in practice, reflecting typical preferences of decision makers who assign a proportionally larger importance to bad outcomes than to good outcomes. Such preference vectors also have advantages from a modeling perspective. Let us define $\cW' = \{\pmb{w}\in\cW : w_1 \ge w_2 \ge \ldots \ge w_K\}$. Using a preference from $\cW'$, the permutation $\pi$ that defines OWA is also a permutation that maximizes the objective function, that is, we have 
\[ \OWA_{\pmb{w}}(\pmb{x},C) = \max_{\pi\in\Pi_K} \sum_{k\in[K]} w_k (\pmb{c}^{\pi(k)} \pmb{x}) \]
where $\Pi_K$ denotes the set of permutations of vectors of size $K$. As discussed in \cite{chassein2015alternative}, using the dual of the maximization problem allows us to reformulate the problem of minimizing OWA as follows:
\begin{align}
\min_{\pmb{x}\in\X} \OWA_{\pmb{w}}(\pmb{x},C) = \min\ & \sum_{k\in[K]} \alpha_k + \beta_k \label{eq:owa1} \\
\text{s.t. } & \alpha_j + \beta_k \ge \sum_{i\in[n]} w_j c^k_i x_i & \forall j,k\in[K] \label{eq:owa2}\\
& \pmb{x}\in\X \label{eq:owa3}\\
& \pmb{\alpha},\pmb{\beta}\in\mathbb{R}^K \label{eq:owa4}
\end{align}
In particular, if $\X$ defines an (integer) linear set of feasible solutions, minimizing OWA becomes an (integer) linear optimization problem as well.

Note that the existence of this easy reformulation is why we focus on the case of non-increasing OWA weights in the following. If desired, our models can be extended by using a more general reformulation of the OWA problem instead, see, for example, the models proposed in \cite{Galand20121540}.

\section{An Optimization Model for Preference Elicitation}\label{sec:opt}

\subsection{Basic Model and Solution Method}
\label{sec:basic}

We assume that the preference vector $\pmb{w}$ is not known. Instead, we would like to identify a suitable vector $\pmb{w}\in\cW'$ based on observations of how a decision maker chooses an alternative. That is, we assume that we are given pairs $(C^1,\pmb{x}^1),\ldots,(C^S,\pmb{x}^S)$ of $S$ historic decisions. The task of preference elicitation is to identify a suitable vector $\pmb{w}$ that can explain the choice of solutions for each observation. Note that we only assume knowledge of past situations and corresponding solution choices and in particular do not require pairwise comparisons between solutions, which may only be available through an interview process.

The underlying idea is to define for each observation $s\in[S]$ the set of preference vectors that can explain this observation, and to find a vector $\pmb{w}$ that is as close to each such set as possible. That is, we define the sets
\[ \opt_s = \left\{ \pmb{w} \in \cW' : \pmb{x}^s \in \argmin_{\pmb{x}\in\X} \OWA_{\pmb{w}}(\pmb{x},C^s) \right\}\]
which contain those weight vectors for which $\pmb{x}^s$ is optimal for the corresponding OWA problem,
and propose to solve
\begin{equation}
\min \left\{ \sum_{s\in[S]} D(\pmb{w}, \opt_s) : \pmb{w}\in\cW' \right\} \tag{Pref}\label{pref}
\end{equation}
where $D: [0,1]^K \times 2^{[0,1]^K} \to \mathbb{R}_+$ is a suitable distance measure between vector and set, given as the distance to its closest element, i.e. $D(x,Y) = \min_{y\in Y} d(x,y)$ for a distance metric $d$.

\begin{figure}[htb]
\begin{center}
\includegraphics[width=0.5\textwidth]{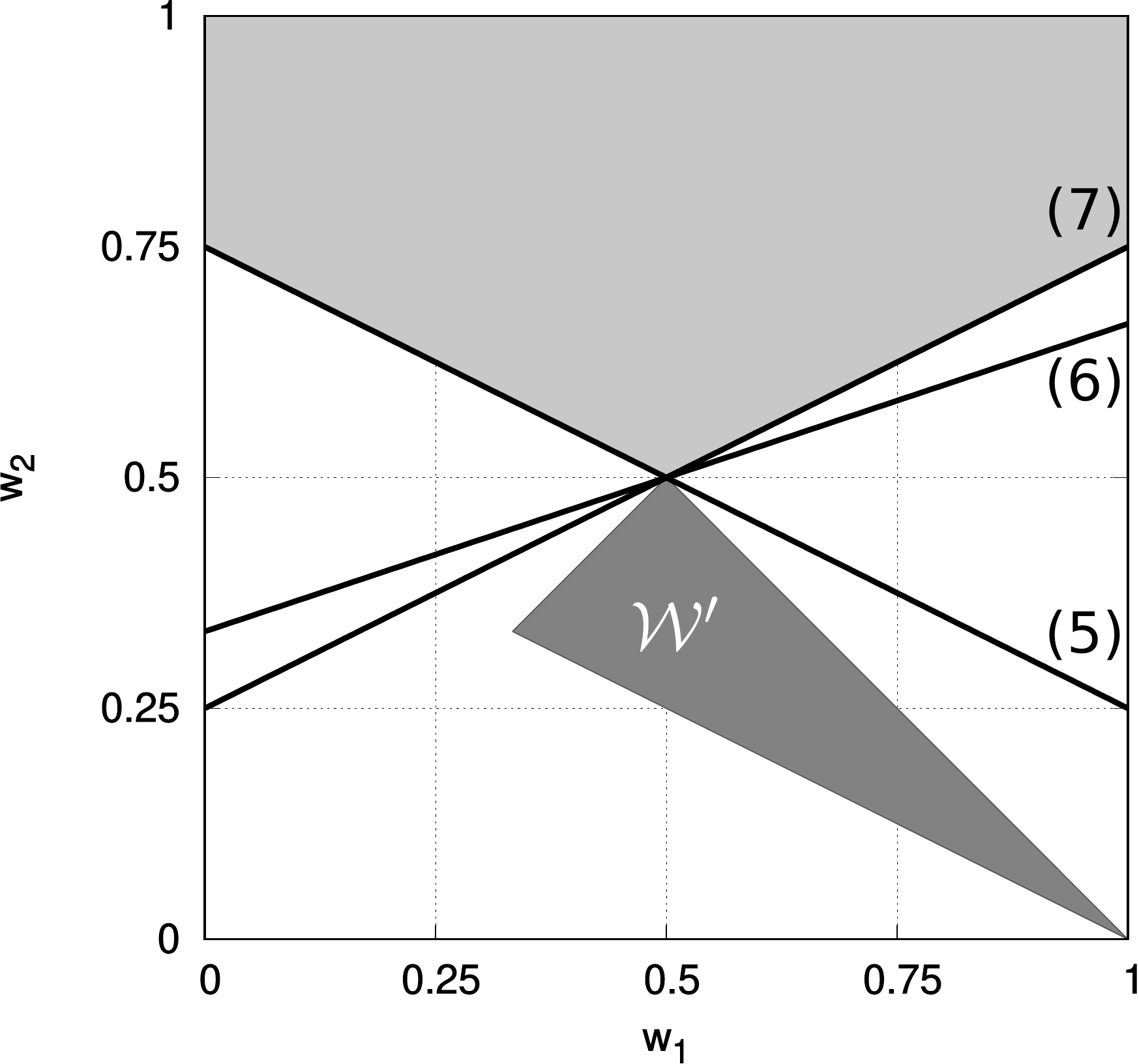}
\end{center}
\caption{Areas for preference vectors in the example problem.}\label{fig:example}
\end{figure}

To clarify this approach, we first present an example. Consider a decision making problem with $\X = \{ \pmb{x}\in\{0,1\}^4 : x_1 + x_2 + x_3 + x_4 = 3\}$, i.e., we need to select three out of $n=4$ given items. As there is only one such observation given, $S=1$ and we drop the index $s$ for simplicity. For this example problem, there are only four solutions: $\pmb{x}^1 = (1,1,1,0)$, $\pmb{x}^2 = (1,1,0,1)$, $\pmb{x}^3 = (1,0,1,1)$ and $\pmb{x}^4 = (0,1,1,1)$. We assume that there are $K=3$ scenarios. For a cost matrix
\[ C = \begin{pmatrix}
1 & 6 & 8 & 4 \\
6 & 7 & 8 & 3 \\
9 & 3 & 2 & 8 \end{pmatrix}\in\mathbb{R}^{K\times n} \]
we are told that the decision maker prefers the solution $\pmb{x}^1 = (1,1,1,0)$. What does this imply for the underlying preference vector? There are four possible solutions. Calculating the sorted vectors of objective values for each solution gives $(21,15,14)$, $(20,16,11)$, $(19,17,13)$ and $(18,18,13)$, respectively. From the choice of $\pmb{x}$ as the preferred solution, we can deduce that its OWA-value is not larger than the OWA value of any other solution. This gives three constraints on the $\pmb{w}$ vector of the form $\OWA_{\pmb{w}}(\pmb{x}^1,C) \le \OWA_{\pmb{w}}(\pmb{x}^i,C)$ for $i=2,3,4$. These are equivalent to the following system of linear equations.
\begin{align}
(21-20)w_1 + (15-16)w_2 + (14-11)w_3 = \phantom{2}w_1 - \phantom{2}w_2 + 3w_3 &\le 0 \label{ex:1} \\
(21-19)w_1 + (15-17)w_2 + (14-13)w_3 = 2w_1 - 2w_2 + \phantom{3}w_3  &\le 0 \label{ex:2}\\
(21-18)w_1 + (15-18)w_2 + (14-13)w_3 = 3w_1 - 3w_2 + \phantom{3}w_3  &\le 0 \label{ex:3}
\end{align}
Substituting $w_3 = 1 - w_1 - w_2$ allows us to plot these equations in two dimensions, see Figure~\ref{fig:example}, where the three black lines indicate points where each of the equations in (\ref{ex:1}-\ref{ex:3}) is fulfilled with equality.

The light gray area contains those vectors $(w_1,w_2)\in[0,1]^2$ that fulfill (\ref{ex:1}-\ref{ex:3}). As we only consider non-increasing weights to reflect a risk-averse decision maker, we furthermore require that $w_1 \ge w_2 \ge w_3$. The dark gray area hence indicates the set $\cW'$. In this example, there is only one element in $\cW'$ that also fulfills the system of equations (\ref{ex:1}-\ref{ex:3}), which is $\opt =  \{ (1/2,1/2,0) \}$. In other words, $\opt$ contains the only feasible preference vector that leads to the observed choice of $\pmb{x}$. It is noteworthy that this uniqueness is only a characteristic of this handcrafted example. In general, $\opt_s$ is a bounded polyhedron if the number of feasible solutions $|\X|$ is finite. Hence, when having additional knowledge or expectations regarding the arising preference vector (e.g. high/low orness) a lexicographic approach can be used to further restrict the set of preference vectors. For multiple observations $S>1$, it is possible that the intersection $\cap_{s\in[S]} \opt_s$ is empty. This may happen if there is no underlying ''true'' value for $\pmb{w}$ that the decision maker uses; decisions may be made intuitively, rather than systematically. For this reason, we would like to find a value of $\pmb{w}$ that is as close as possible to those values $\opt_s$ that can explain each observation. It may also happen that for some $s\in[S]$, $\opt_s = \emptyset$, i.e., there is no possible risk-averse preference vector that can explain that particular specific choice. In this case, Problem~\eqref{pref} needs to be slightly modified, which is explained later in this section. Also note that while \eqref{pref} is defined using the same set $\X$ for each observation $s$, this assumption might also be relaxed to allow different sets of feasible solutions for each decision making observation. We only require that there is the same number of scenarios $K$ for each problem, so that vectors $\pmb{w}$ are of the same dimension.

We now discuss how to solve Problem~\eqref{pref}. For each $\pmb{x}\in\X$, let us denote by $a^s_1(\pmb{x}),\ldots,a^s_K(\pmb{x})$ the objective values $C^s\pmb{x}$, sorted from largest to smallest value. Using this notation, we have
\[ \opt_s = \left\{ \pmb{w} \in \cW' : \sum_{k\in[K]} w_k a^s_k(\pmb{x}^s) \le \sum_{k\in[K]} w_k a^s_k(\pmb{x}) \quad \forall \pmb{x}\in\X  \right\} \]
This allows us to rewrite Problem~\eqref{pref} as follows:
\begin{align}
\min\ & \sum_{s\in[S]} d(\pmb{w},\pmb{w}^s) \label{mod:1} \\
\text{s.t. } 
&  \sum_{k\in[K]} w^{s}_k a^s_k(\pmb{x}^s) \le \sum_{k\in[K]} w^{s}_k a^s_k(\pmb{x}) & \forall s\in[S], \pmb{x}\in\X \label{mod:2} \\
& \sum_{k\in[K]} w^s_k = 1 & \forall s\in[S] \label{mod:3}\\
& \sum_{k\in[K]} w_k = 1  \label{mod:4}\\
& w^s_1 \ge w^s_2 \ge \ldots w^s_K \ge 0 & \forall s\in[S] \label{mod:5}\\
& w_1 \ge w_2 \ge \ldots w_K \ge 0  \label{mod:6}
\end{align}
Note that the constraints~(\ref{mod:2}-\ref{mod:6}) are linear, as values $a^s_k(\pmb{x})$ can be precomputed. However, depending on the number of elements in $\X$, there might be infinitely or exponentially many constraints of type~\eqref{mod:2}.

Considering the distance $d$, several choices result in tractable optimization problems. Using the 1-norm $d_1(\pmb{w},\pmb{w}^s)=\vert\vert \pmb{w}-\pmb{w}^s \vert\vert_1 = \sum_{k\in[K]}|w_k-w_k^s|$, objective function~\eqref{mod:1} can be linearized by introducing new variables $d^s_k\ge 0$ with constraints $-d^s_k \le w_k - w^s_k \le d^s_k$ for all $s\in[S]$ and $k\in[K]$. By minimizing $\sum_{s\in[S]} \sum_{k\in[K]} d^s_k$, we hence minimize the sum of absolute differences in each component. It is also possible to use weights on positions $k\in[K]$, e.g., one unit of difference in $w_1$ may be more important than one unit of difference in $w_K$. Alternatively, also the $\infty$-norm $d_\infty(\pmb{w},\pmb{w}^s)=\vert\vert \pmb{w}-\pmb{w}^s \vert\vert_\infty =\max_{k\in[K]}\vert w_k-w_k^s \vert $ can be linearized, while using the 2-norm $d_2(\pmb{w},\pmb{w}^s)=\vert\vert \pmb{w}-\pmb{w}^s \vert\vert_2 =\sqrt{\sum_{k\in[K]}(w_k-w_k^s )^2}$ results in a convex quadratic optimization problem.

\begin{figure}[htb]
\begin{center}
\includegraphics[width=0.8\textwidth]{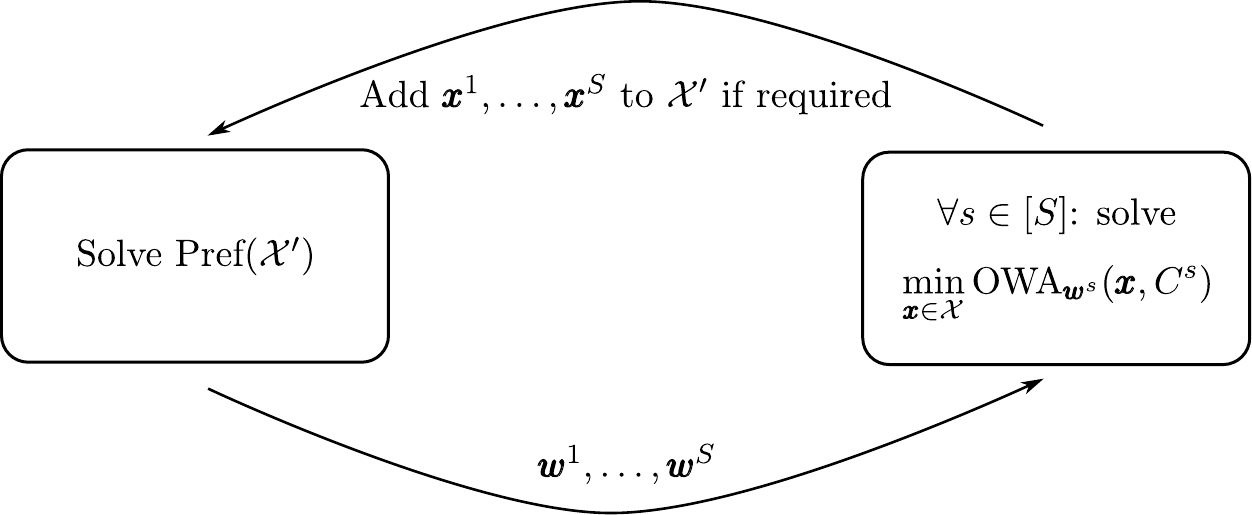}
\caption{Iterative solution algorithm for \eqref{pref}.}\label{fig:alg}
\end{center}
\end{figure}

We now discuss how to treat the large number of constraints of type~\eqref{mod:2}. We propose to use an iterative solution procedure that is summarized in Figure~\ref{fig:alg}. 

We begin with a finite subset of alternatives $\X' \subseteq \X$, for example, we may set $\X' = \{\pmb{x}^1,\ldots,\pmb{x}^S\}$. We solve the optimization problem \eqref{pref}, where we restrict constraints~\eqref{mod:2} to $\X'$. We refer to this problem as \ref{pref}($\X'$). This is a linear program, given that $d$ can be linearized. The result is a candidate preference vector $\pmb{w}\in\cW'$, along with preference vectors $\pmb{w}^s$ for each $s\in[S]$. We then consider each $s\in[S]$ separately to check if indeed it holds that $\pmb{w}^s \in\opt_s$. To this end, we need to compare the value $\OWA_{\pmb{w}^s}(\pmb{x}^s,C^s)$ with $\min_{\pmb{x}\in\X} \OWA_{\pmb{w}^s}(\pmb{x},C^s)$. If it turns out that there exists a solution that results in a better objective value than $\pmb{x}^s$ under $\pmb{w}^s$, we add this solution to $\X'$ and repeat the process.

Note that this solution method ends after a finite number of iterations, if $\X$ is finite. This holds, e.g., if $\X\subseteq\{0,1\}^n$ is the set of feasible solutions for a combinatorial optimization problem. In this case, we have found an optimal solution to \eqref{pref}. Alternatively, we may stop the method after a fixed number of iterations is reached or if the difference between consecutive solutions for $\pmb{w}$ becomes sufficiently small. Further note that in this iterative method, we alternate between solving a linear program to identify a candidate for $\pmb{w}$ and standard OWA problems. Hence, it can be extended to find weight vectors that do not have the non-increasing property in the same way.

To conclude the description of our approach, we consider the case that $\opt_s=\emptyset$ for some $s\in[S]$, i.e., the decision maker chooses a solution that is optimal with respect to no preference vector in $\cW'$. This means that some model \ref{pref}($\X'$) that we try to solve in the iterative procedure becomes infeasible. In that case, we can follow a lexicographic approach to consider those preferences $\pmb{w}^s$ that remain as close to optimality as possible: We first minimize the violation of constraint~\eqref{mod:2} (which so far was supposed to be zero) and then find the solution that minimizes the objective function with that best smallest possible violation. To this end, we solve the optimization problem
\begin{equation}
\infeas_s := \min_{\pmb{w}\in\cW'} \max_{\pmb{x}\in\X} \Big( \OWA_{\pmb{w}}(\pmb{x}^s,C^s) - \OWA_{\pmb{w}}(\pmb{x},C^s) \Big) \label{eq:inf}
\end{equation}
to calculate the smallest possible violation of the corresponding constraint~\eqref{mod:2}. We then replace the constraint with
\[  \sum_{k\in[K]} w_k^{s} a^s_k(\pmb{x}^s) \le \sum_{k\in[K]} w_k^{s} a^s_k(\pmb{x}) + \infeas_s \quad \forall \pmb{x}\in\X \]
to ensure feasibility of Problem~\eqref{pref}. To solve the problem in~\eqref{eq:inf}, we can use an iterative procedure analogously to the method described in Figure~\ref{fig:alg}. Alternatively, $\infeas_s$ is used as an additional non-negative variable in Problem~\eqref{pref}, where we modify the objective to additionally minimize $\sum_{s\in[S]} \infeas_s$ with a large weight to give it priority over minimizing the sum of distances in the preference vectors. Note, however, that the existence of an appropriate weight depends on the problem structure, making this a more heuristic approach, in general.

\subsection{Heuristic Compact Model}

We now reconsider model~\eqref{pref} and develop a heuristic formulation for the case that $\X$ is given by a polyhedron, i.e., we assume that $\X = \{ \pmb{x} \ge 0 : A\pmb{x} \ge \pmb{b}\}$ for $A\in\mathbb{R}^{m\times n}$ and $\pmb{b}\in\mathbb{R}^m$. Recall that some combinatorial optimization problems, such as the shortest path or the minimum spanning tree problems, also allow for such linear programming formulations where binary variables are not required. 

Constraints~\eqref{mod:2} ensure for each $s\in[S]$ that solutions $\pmb{x}^s$ are indeed optimal, that is, they are equivalent to
\[ \OWA_{\pmb{w}}(\pmb{x}^s) \le \min_{\pmb{x}\in\X} \OWA_{\pmb{w}}(\pmb{x}).\] 
Previously, we treated the right-hand side by constructing one constraint per solution $\pmb{x}\in\X$. As $\X$ is a polyhedron, we can use strong duality to reformulate the right-hand side as a maximization problem with the same optimal objective value.  The dual of the OWA problem with costs $C^s$, see \eqref{eq:owa1}-\eqref{eq:owa4}, when using $\pmb{\pi}$ and $\pmb{\sigma}$ as dual variables corresponding to constraints~\eqref{eq:owa2} and \eqref{eq:owa3}, respectively, is given as:
\begin{align*}
\max\ & \pmb{b}^t \pmb{\sigma} \\
\text{s.t. } & \sum_{k\in[K]} \pi_{jk} = 1 & \forall j\in[K] \\
& \sum_{j\in[K]} \pi_{jk} = 1 & \forall k\in[K] \\
& A^t\pmb{\sigma} \le \sum_{j\in[K]} \sum_{k\in[K]} w^s_j \pmb{c}^{s,k} \pi_{jk} \\
&\pmb{\sigma} \ge \pmb{0} \\
&\pmb{\pi} \ge \pmb{0}
\end{align*} Furthermore, by weak duality, the objective value of any feasible solution to the dual problem gives a lower bound to the optimal primal objective value, making the substitution valid not only for the optimal, but for any solution of the dual.
We can thus find the following compact reformulation of model~\eqref{pref}.
\begin{align*}
\min\ & \sum_{s\in[S]} d(\pmb{w},\pmb{w}^s) \\
\text{s.t. } & \sum_{k\in[K]} w^s_k a^s_k(\pmb{x}^s) \le \pmb{b}^t\pmb{\sigma}^s & \forall s\in[S] \\
& \sum_{k\in[K]} \pi^s_{jk} = 1 & \forall j\in[K], s \in[S] \\
& \sum_{j\in[K]} \pi^s_{jk} = 1 & \forall k\in[K], s \in[S] \\
& A^t\pmb{\sigma}^s \le \sum_{j\in[K]} \sum_{k\in[K]} w^s_j \pmb{c}^{s,k} \pi_{jk}  & \forall s\in[S] \\
& \pmb{w}\in\cW' \\
& \pmb{w}^s \in\cW' & \forall s\in[S] \\
&\pmb{\sigma}^s \ge \pmb{0} & \forall s\in[S]\\
&\pmb{\pi}^s \ge \pmb{0} & \forall s\in[S]
\end{align*}

Note that the product $w^s_j\pi_{jk}$ introduces a non-linearity due to the continuous nature of both variables. To address this non-linearity and facilitate analysis, we utilize McCormick envelopes as an approximation technique. By employing McCormick envelopes, we approximate $\tau^s_{jk} = w^s_j\pi_{jk}$ while imposing the following constraints: $\tau^s_{jk} \le w^s_k$ and $\tau^s_{jk} \le \pi_{jk}$. Note that this approximation does not yield an equivalent reformulation of the problem; instead, it provides a heuristic model that circumvents the need for an iterative solution method. The adoption of McCormick envelopes offers an effective means of approximating the non-linear term.

\subsection{Alternative Model Formulation}

In the preference elicitation model~\eqref{pref}, we minimize the distance in preference vectors $\pmb{w}$ and thus focused on finding preference vectors that reflect the OWA values of the observations. Alternatively, one might be interested in a preference vector that mimics the observations to the extend of reproducing similar solutions, rather than obtaining solutions with similar OWA values. In particular, a slight modification of $\pmb{w}$ can yield a significantly changed solution, while the change in the OWA value itself remains small. It is possible that a very different preference vector may result in only slightly
different optimal solutions.
Hence, if it is more important to accurately mimic observed solutions rather than obtaining solutions of similar OWA score we  may redefine our model as follows. The set of optimal solutions for a given preference vector $\pmb{w}$ is denoted as
\[ \opt'_s(\pmb{w}) = \Big\{ \pmb{y}\in\X: \OWA_{\pmb{w}}(\pmb{y},C^s) \le \OWA_{\pmb{w}}(\pmb{x},C^s) \ \forall \pmb{x}\in\X \Big\}, \]
that is, while $\opt_s$ is a subset of $\mathbb{R}^K$, we now consider a subset of the set of feasible solutions $\X$.  The modified preference elicitation model now asks for a vector $\pmb{w}\in\cW'$ such that the sets $\opt'_s(\pmb{w})$ are close to the given solutions $\pmb{x}^s$ in a distance metric $d':\X \times \X \to \mathbb{R}_+$. We thus define the following problem:
\begin{equation}
\min\ \left\{ \sum_{s\in[S]} d'(\pmb{y}^s,\pmb{x}^s) : \pmb{w}\in\cW',\ \pmb{y}^s\in \opt'_s(\pmb{w})\ \forall s\in[S] \right\} \tag{Pref'} \label{altpref}
\end{equation}
To solve this problem, we need a tractable formulation of constraints
\[ \pmb{y}^s\in \opt'_s(\pmb{w})\quad \forall s\in[S], \]
which are equivalent to
\begin{equation}
\sum_{k\in[K]} w_ka^s_k(\pmb{y}^s) \le \sum_{k\in[K]} w_k a^s_k(\pmb{x})\quad \forall \pmb{x}\in\X, s\in[S]. \label{auxcon}
\end{equation}
The values $a^s_1(\pmb{y}^s) \ge \ldots \ge a^s_K(\pmb{y}^s)$ give a worst-case sorting of the objective values of $\pmb{y}^s$ in observation $s$. Hence, we must ensure that the left-hand side value of equation~\eqref{auxcon} is not underestimated. Using the assumption that $\pmb{w}\in\cW'$, we have the equivalent formulation
\[
\max_{\pi\in\Pi_K} \sum_{k\in[K]} w_{\pi(k)} \pmb{c}^{s,\pi(k)} \pmb{y}^s \le \sum_{k\in[K]} w_k a^s_k(\pmb{x})\quad \forall \pmb{x}\in\X, s\in[S].
\]
Applying the same dualization technique as in the $\OWA$ reformulation in (\ref{eq:owa1}-\ref{eq:owa4}), we introduce new variables $\pmb{\alpha}^s$ and $\pmb{\beta}^s$ to reformulate Problem~\eqref{altpref} as follows.
\begin{align}
\min\ & \sum_{s\in[S]} d'(\pmb{y}^s,\pmb{x}^s) \label{eq:var1} \\
\text{s.t. } 
& \sum_{k\in[K]} \alpha^s_k + \beta^s_k \le \sum_{k\in[K]} w_k a^s_k(\pmb{x}) & \forall \pmb{x}\in\X, s\in[S] \label{eq:var2} \\
& \alpha^s_j + \beta^s_k \ge \sum_{i\in[n]} w_j c^{s,k}_i y^s_i & \forall s\in[S], j,k\in[K] \\
& \pmb{y}^s \in\X & \forall s\in[S] \label{eq:var3} \\
& \sum_{k\in[K]} w_k = 1 \label{eq:var4} \\
& w_1 \ge w_2 \ge \ldots \ge w_K \ge 0 \label{eq:var5}
\end{align}
As before, values $a^s_k(\pmb{x})$ can be precomputed for a given $\pmb{x}$. To treat the potentially exponential or infinite number of constraints in~\eqref{eq:var2}, an iterative solution method analogous to the one presented in Section~\ref{sec:basic} can be applied. If $\X$ is the set of feasible solutions of a combinatorial decision making problem, then we can linearize products $w_jy^s_i$ by introducing additional variables $\tau^s_{ji}\ge 0$ with $\tau^s_{ji} \ge w_j + y^s_i - 1$. In this case, a natural choice for a distance measure $d'$ may be the Hamming distance, where we simply count the number of differing entries. In our notation, this means we minimize
\[ d'(\pmb{y}^s,\pmb{x}^s) = \sum_{i\in[n]: x^s_i = 1} (1-y^s_i) + \sum_{i\in[n] : x^s_i = 0} y^s_i \]
which implies that Problem~\eqref{altpref} can be solved through a sequence of mixed-integer programming formulations.
While models \eqref{pref} and \eqref{altpref} may provide different solutions (in Appendix \ref{App::ExampleAlternative}, we provide such an example), we note that if there is a vector $\pmb{w}$ that can explain the choice of solutions for each observation (i.e., there is $\pmb{w}$ such that each $\pmb{x}^s$ is an optimizer with respect to $\OWA_{\pmb{w}}(\cdot,C^s)$), then this $\pmb{w}$ is optimal for both models \eqref{pref} and \eqref{altpref}. Hence, we can only expect to see different solutions if the decision maker did not use the same $\pmb{w}$ consistently throughout all observations. This result is formalized in Appendix~\ref{App::ExampleAlternative}.

\section{Computational Experiments}
\label{sec:exp}

\subsection{Setup}

The purpose of these experiments is to evaluate the degree to which decision maker preferences can be identified based on observations using the optimization approach compared to having pairwise comparisons at hand, which is the prevalent approach in preference elicitation.

We mainly focus on selection problems, where $p$ out of $n$ items must be chosen, i.e., $\X = \{ \pmb{x}\in\{0,1\}^n : \sum_{i\in[n]} x_i = p\}$.  Problems of this type are frequently considered, see, e.g., \cite{chassein2018recoverable}, and have the advantage that the decision maker is mostly unconstrained in her choice of items, which means that the underlying preference becomes more important in the decision making process.  In the appendix we also provide results for assignment as well as knapsack problems (see Appendix~\ref{Sec::AssignmentKnapsack}). We note that the experimental results exhibit a consistent behavior across all three problem settings.

To generate preference vectors $\pmb{w}$ we utilize their \textit{orness}, which is an important characterizing measures introduced by  \cite{yager1988ordered}, given by 
$$orness(\pmb{w})=\frac{1}{K-1}\sum_{k\in[K]} (K-k)w_k\, .$$
Note that as we only assume risk-averse weights in these experiments, the orness of such vectors is always greater than or equal to $0.5$, which is the orness of the average criterion $(\frac{1}{K},\ldots,\frac{1}{K})$, whereas $(1,0,\ldots,0)$, which is the worst-case criterion in our minimization setting, has an orness of $1$ \citep{yager1993families,liu2004properties}.
After sampling the orness using a uniform distribution over the interval $[0.5,1]$, this value is handed to the model proposed by \cite{wang2005minimax} which then yields a preference vector with the specified orness. The model is given in Appendix \ref{App::OWAOrness}. 

An instance is created by first sampling $S$ cost matrices $\bar{C}^s \in \{1,100\}^{K \times n}$, where each entry  is generated uniformly as a random integer in $\{1,100\}$. To ensure that all objectives are commensurate to each other, we normalize each of the $k\in[K]$ cost vectors using min-max normalization, to obtain the final  matrix $C^s \in [0,1]^{K \times n}$: 
$$C^{s}_{k,i}=\frac{\bar{C}^{s}_{k,i}-\min_{j\in[n]}\bar{C}^{s}_{k,j}}{\max_{j\in[n]}\bar{C}^{s}_{k,j}-\min_{j\in[n]}\bar{C}^{s}_{k,j}}$$
Additionally, we allow noise in the data, i.e.~that the decision maker may deviate from her underlying preference vector for each decision. To model this behavior, we include a noise parameter $\epsilon \in [0,1]$. Whenever an OWA problem is solved, we modify each value of $w_k$ by adding a uniformly random value in $[\max\{-w_k,-\epsilon\},+\epsilon]$. Afterwards, $\pmb{w}$ is normalized and sorted. 

Hence, the most important parameters for a selection instance are $n$, $p$, $K$, $S$, and $\epsilon$. Our \emph{basic setting} is $n=40$, $p=\frac{n}{2}$, $K=5$ objectives, $S=16$ observations and $\epsilon=0$ noise.

For the optimization approach, we implemented the iterative method presented in Section~\ref{sec:basic} to solve model~\eqref{pref} with the 1-norm as distance measure. The method was implemented in C++ using CPLEX version 22.1 to solve optimization models.

We compare our elicitation method based on observations with having to obtain pairwise preferences. In order to transform the pairwise comparisons into a preference vector we use the model proposed by \cite{ahn2008preference} that minimizes the violation of the decision maker’s judgment. The model can be found in Appendix~\ref{App::OwaAhn}. We compare our approach with having to conduct $1$, $5$, $10$, and $20$ pairwise comparisons \emph{for each} of the $S$ observed data sets.
 We generate solution pairs for each data set $s\in [S]$ by obtaining two random supported Pareto-optimal solutions, by randomly creating a vector $\pmb{w} \in \cW$ and computing the optimal solution for such preference.
In preliminary experiments we also looked at pairs consisting of two random solutions as well as one Pareto-optimal solution vs. the optimal solution. Here comparing two Pareto-optimal solutions yielded the best results, performing similarly good as comparing two random solutions. 

Recall that having knowledge of an optimal solution in our model implies an indirect comparison of this one optimal solution with all other possible solutions. Usually, this implies a much larger number of implicit comparisons. We thus compare how many direct comparisons result in roughly the same quality of results as the set of indirect comparisons implied by the set of observations for model~\eqref{pref}.

To evaluate the quality of these methods, we compare the difference between the predicted preference vector and the actual preference vector in the 2-norm. Using a different norm than the 1-norm used in the objective function of our model is supposed to ensure a fairer comparison. Additionally, we test the performance of the different approaches by measuring the in-sample and out-of-sample performance by calculating the average Hamming distance between the optimal solution and the solution obtained when using the proposed preference vector over all $S$ observations and over $100$ randomly created new data points, respectively.
As the two proposed models \eqref{pref} and \eqref{altpref} result in the same set of optimal solutions if there is no noise (see Appendix \ref{App::ExampleAlternative}), we only use \eqref{altpref} in experiments where the underlying preference vector is disturbed, resulting in observations that cannot be described perfectly with a single preference vector.

\subsection{Results}
\subsubsection{Orness Evaluation\label{sec::ExpOrness}}

\begin{figure}[!htb]
\centering
\subfigure[Average preference vector distance.\label{fig::ManyForOrness}]
{\includegraphics[width=.32\textwidth]{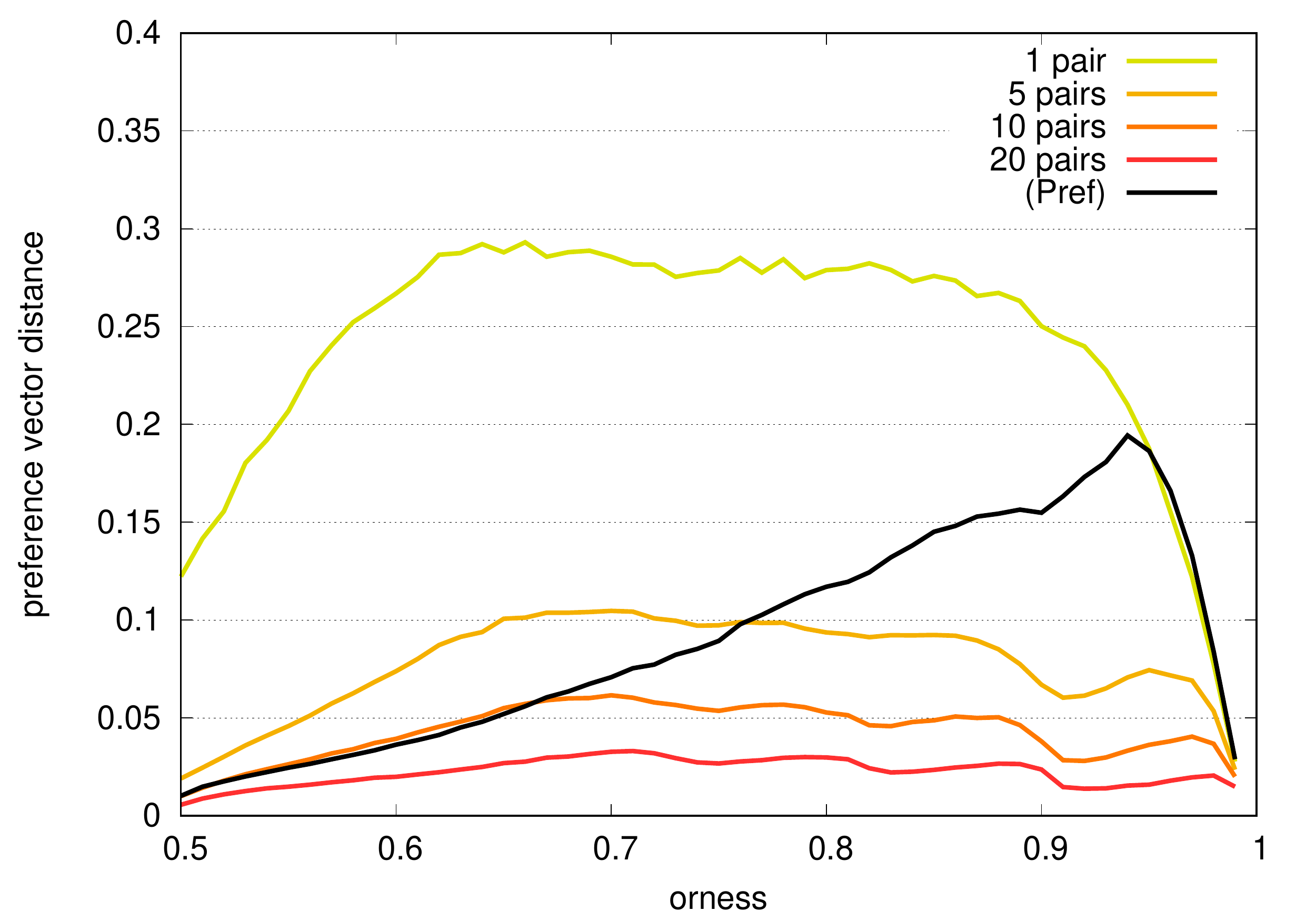}}%
\subfigure[Ratio of proposed preference vectors with orness 1. \label{fig::10_Orness}]
{
\includegraphics[width=.32\textwidth]{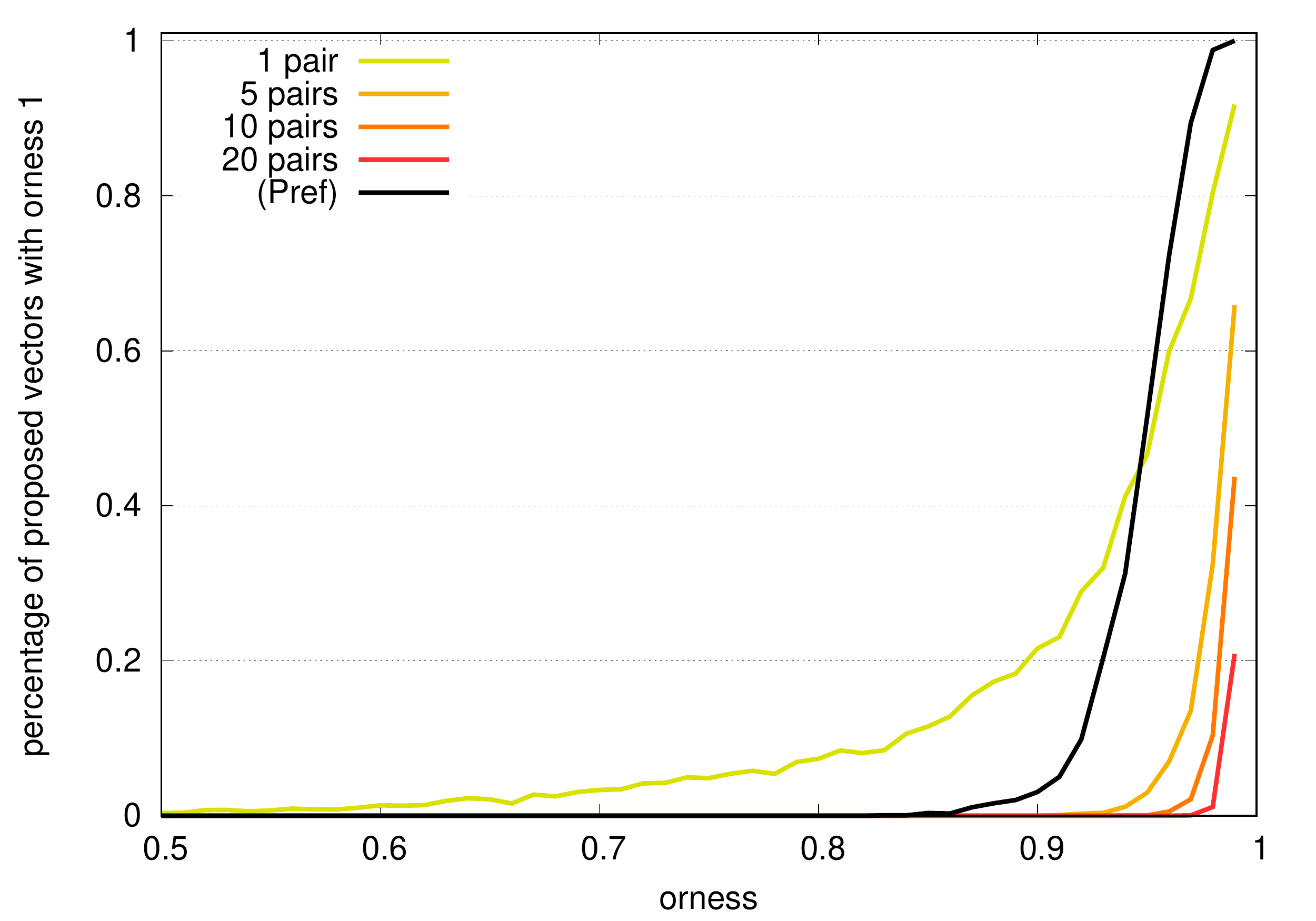}}%
\subfigure[Ratio of random optimal vectors. \label{fig::Optimal_Orness}]
{
\includegraphics[width=.32\textwidth]{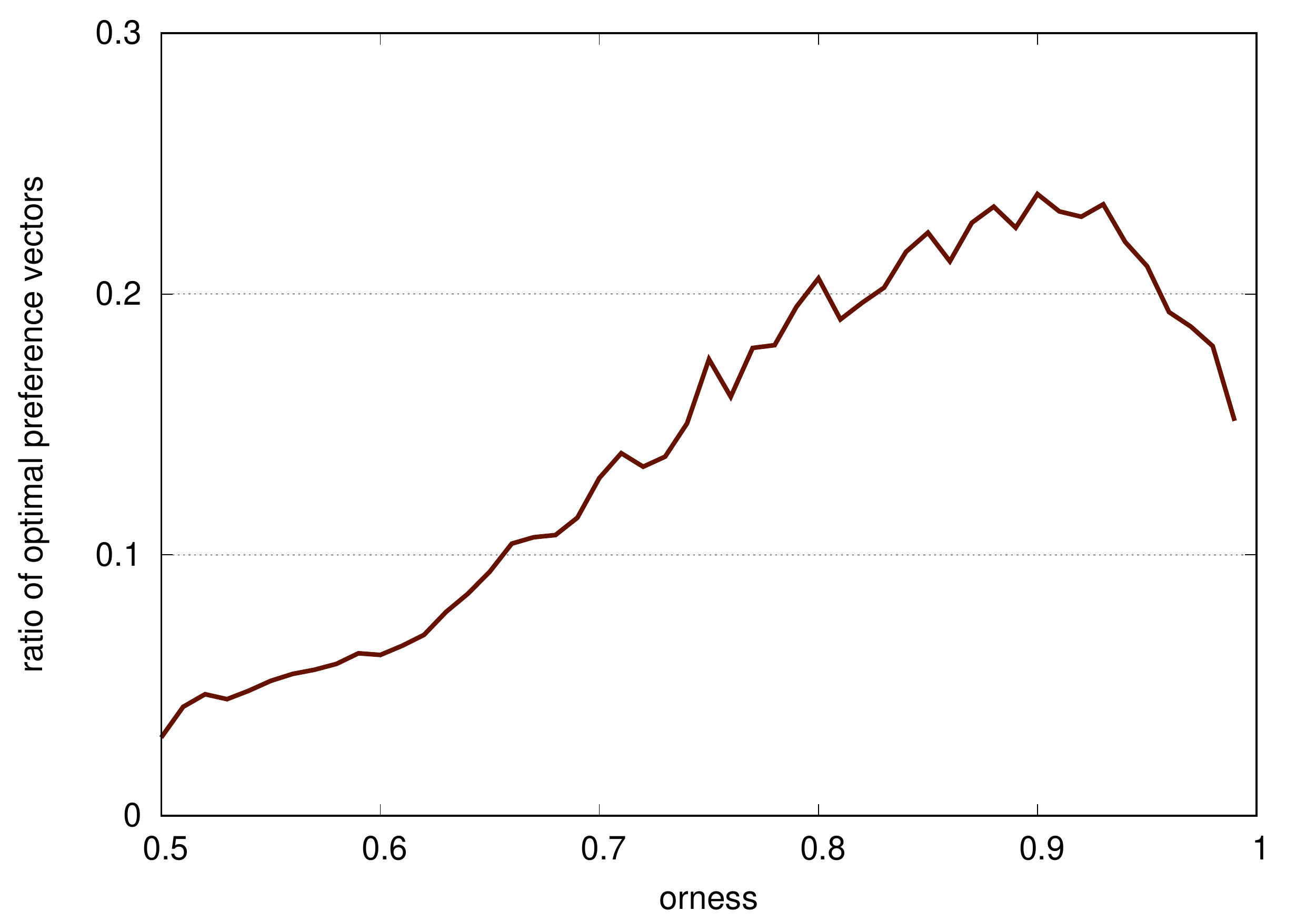}}
\caption{The average preference vector difference (left), the ratio of preference vectors with an orness of 1 (middle), and the ratio of random preference vectors that are optimal in \eqref{pref} (right), each dependent on the orness of the underlying true preference vector.}
\end{figure}

We first want to gain insights into the relevance of the orness of the resulting preference vector. For Figures~\ref{fig::ManyForOrness} and \ref{fig::10_Orness} we created $100\,000$ preference vectors for the basic setting ($n=40$, $S=16$, $K=5$, $\epsilon=0$). We first plotted the average preference vector distance of the true underlying preference vector compared to the proposed preference vector of our optimization model \eqref{pref} as well as of the approach using pairwise comparisons.

For low orness our model performs very well even compared to 20 pairwise comparisons per observation, which are overall $320$ comparisons by the decision maker, in contrast to using only $S=16$ observations. However, for increasing orness the average distance to the true preference vector increases even to the extend that our approach performs similar to having just a single comparison available. One reason for this is that solutions obtained by preference vectors with high orness can quite often be explained via the preference vector $(1,0,\ldots,0)$, i.e.~the vector with orness 1:
In Figure~\ref{fig::10_Orness} the percentage of proposed preference vectors that have an orness of 1 are shown dependent on the orness of the true OWA weights and it can be observed that when using pairwise comparisons the less comparisons are conducted the more often OWA weights with orness 1 are proposed. For our approach, even though we used $S=16$ observations, OWA weights with true orness larger than $0.95$ are explained via the vector $(1,0,\ldots,0)$ in more than fifty percent of the cases. This is largely an artifact of the used optimization solver, which returns $(1,0,\ldots,0)$, if it is feasible.
To also understand how many preference vectors are capable of explaining the observed scenarios we look at the ratio of optimal preference vectors depending on the orness of the true preference vector. For Figure \ref{fig::Optimal_Orness} we generated 10000 preference vectors with random orness for $S=4$ observations, $K=5$ and $n=40$. For each instance, we sample 1000 additional random preference vectors and check if this preference vector can explain the given observations.
We can observe that for preference vectors with orness of about $0.9$, about $20\%$ of randomly drawn vectors perfectly explain the observation. This also serves as an explanation of why our approach performs less well for such vectors.

\subsubsection{Instance Parameters}

\begin{figure}[!htb]
\centering
\subfigure[\WDist\label{fig::ExpNW}]{\includegraphics[width=.32\textwidth]{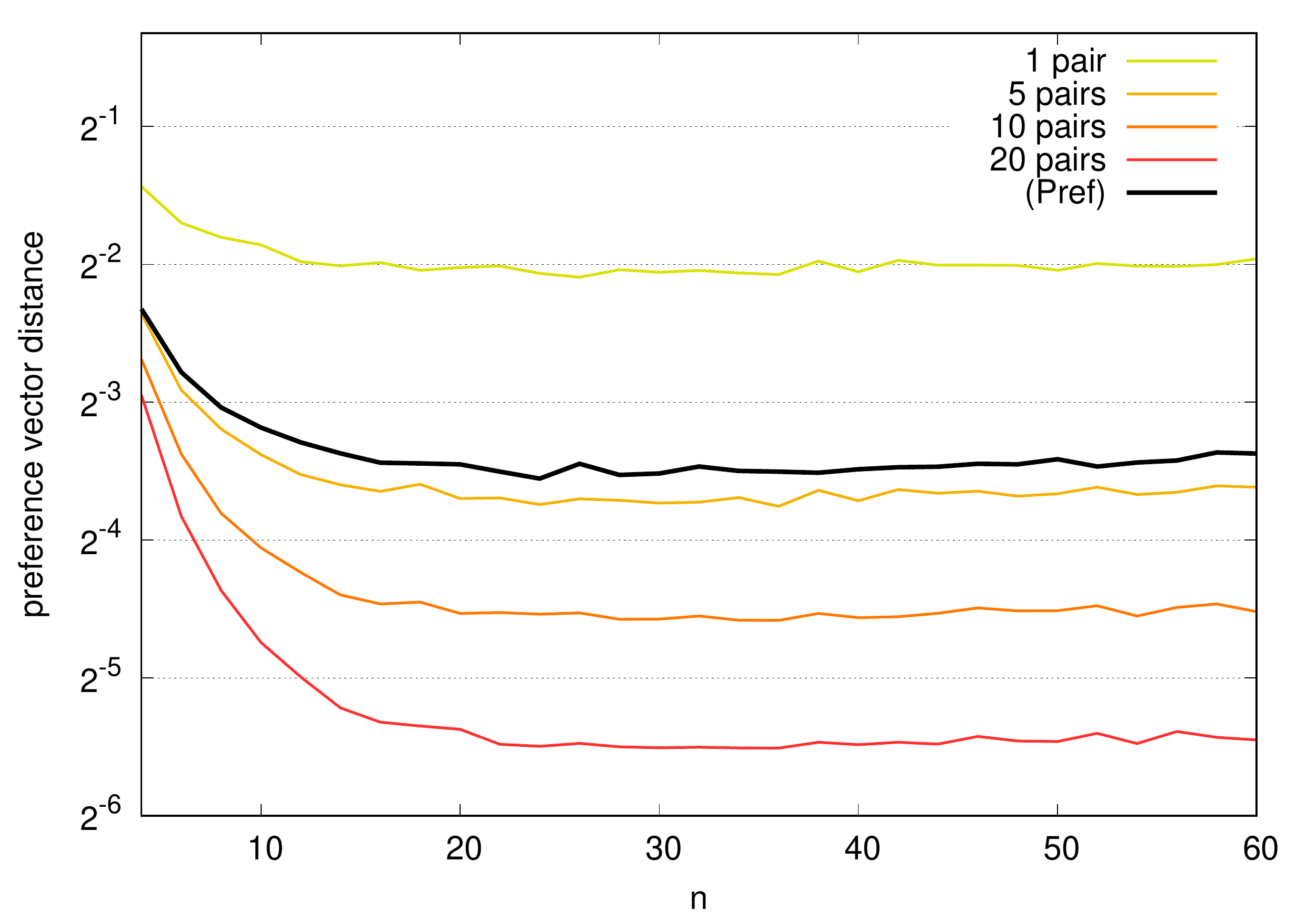}}%
\subfigure[\XInDist\label{fig::ExpNXIn}]{\includegraphics[width=.32\textwidth]{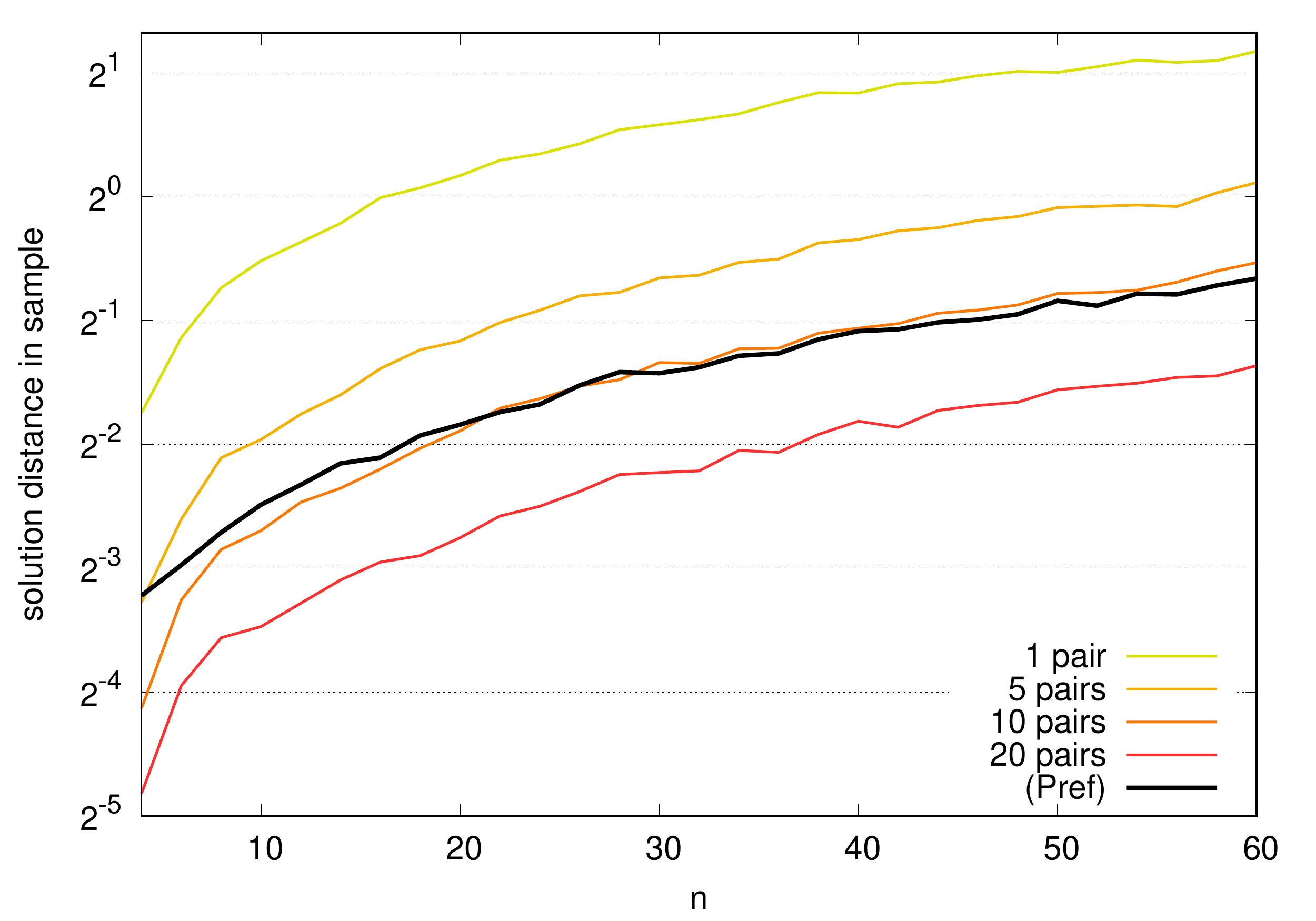}}%
\subfigure[\XOutDist\label{fig::ExpNXOut}]{\includegraphics[width=.32\textwidth]{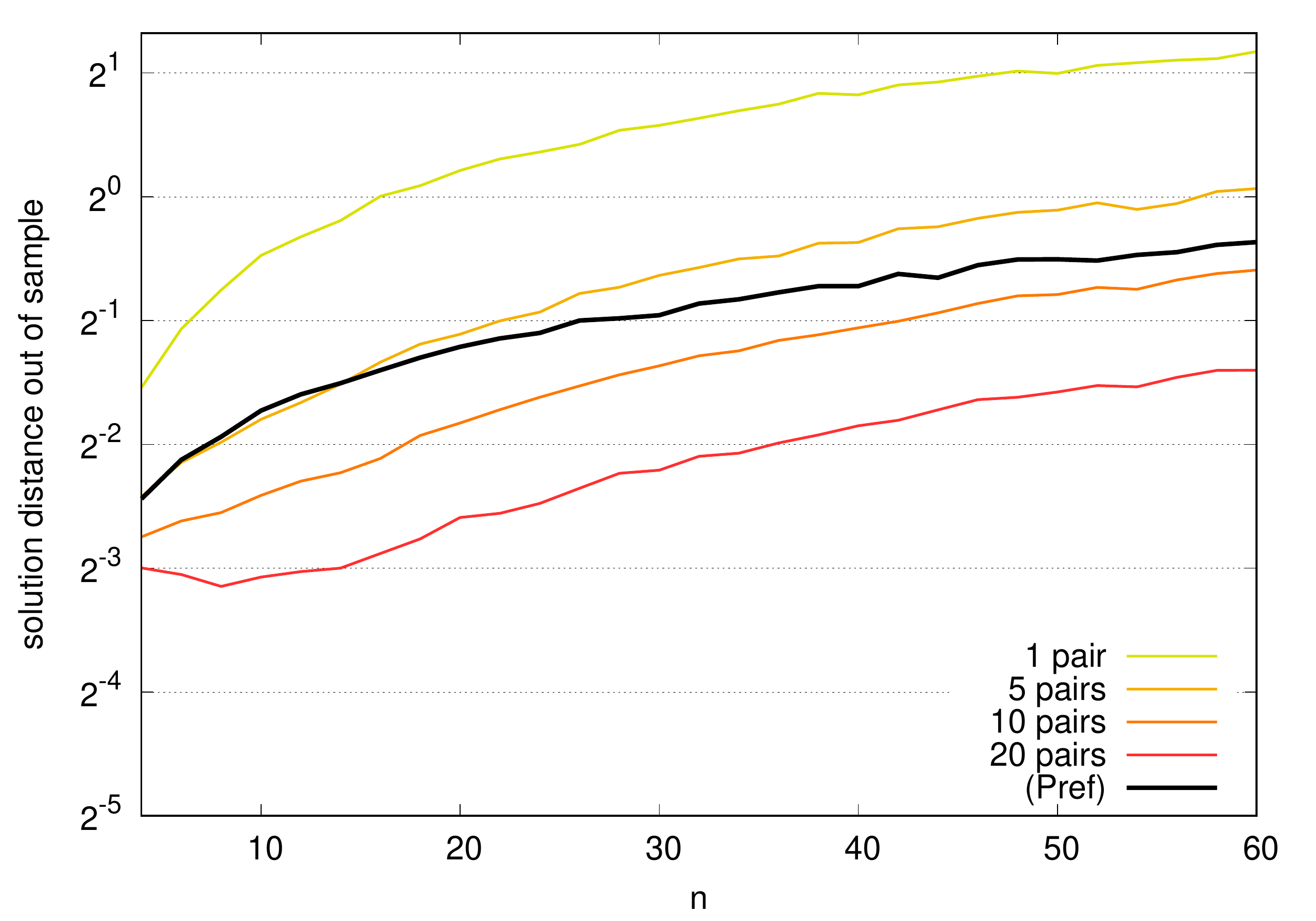}}%
\caption{Comparison of our approach with pairwise comparisons for varying problem size $n$.\label{fig::N}}
\end{figure}

We now examine how the parameters $n$, $S$ and $K$ influence the performance of our approach. To this end, we vary these parameters and again compare them to the pairwise comparison based approach using the Euclidean distance to the true preference vector as well as the in- and out-of-sample distance of solutions as measure. First, for $n\in \{10,11\ldots,60\}$ and $p=\lfloor \frac{n}{2} \rfloor$ we generated each $1000$
random instance ($S=16$, $K=5$) and compared the result of our model to the  pairwise comparison. In Figure \ref{fig::N} the results are shown.
Viewing the distance of the obtained preference vector to the true underlying preference vector our approach performs similarly, but slightly worse, than having five comparisons per observation at hand (Fig.~\ref{fig::ExpNW}). By increasing the number of comparisons the distance to the underlying preference vector can be further reduced. Beyond a certain threshold, the preference vector distance seems to remain unaffected by a further increasing problem size $n$. This is different when considering the Hamming distance of solutions obtained by the true and the obtained preference vector, as for increasing $n$ this distance tends to increase. However, up to $n=60$ the average distance both on in-sample as well as on out-of-sample data is less than $1$. Due to the nature of the selection problem this means that on average every other solution is recreated perfectly. We can also observe, that our approach performs similar to having ten comparisons per observation when it comes to replicating the observations (Fig.~\ref{fig::ExpNXIn}). Furthermore, for increasing $n$ and on new data our approach tends to outperform having five comparisons (Fig.~\ref{fig::ExpNXOut}).   

\begin{figure}[!htb]
\centering
\subfigure[\WDist\label{fig::ExpSW}]{\includegraphics[width=.32\textwidth]{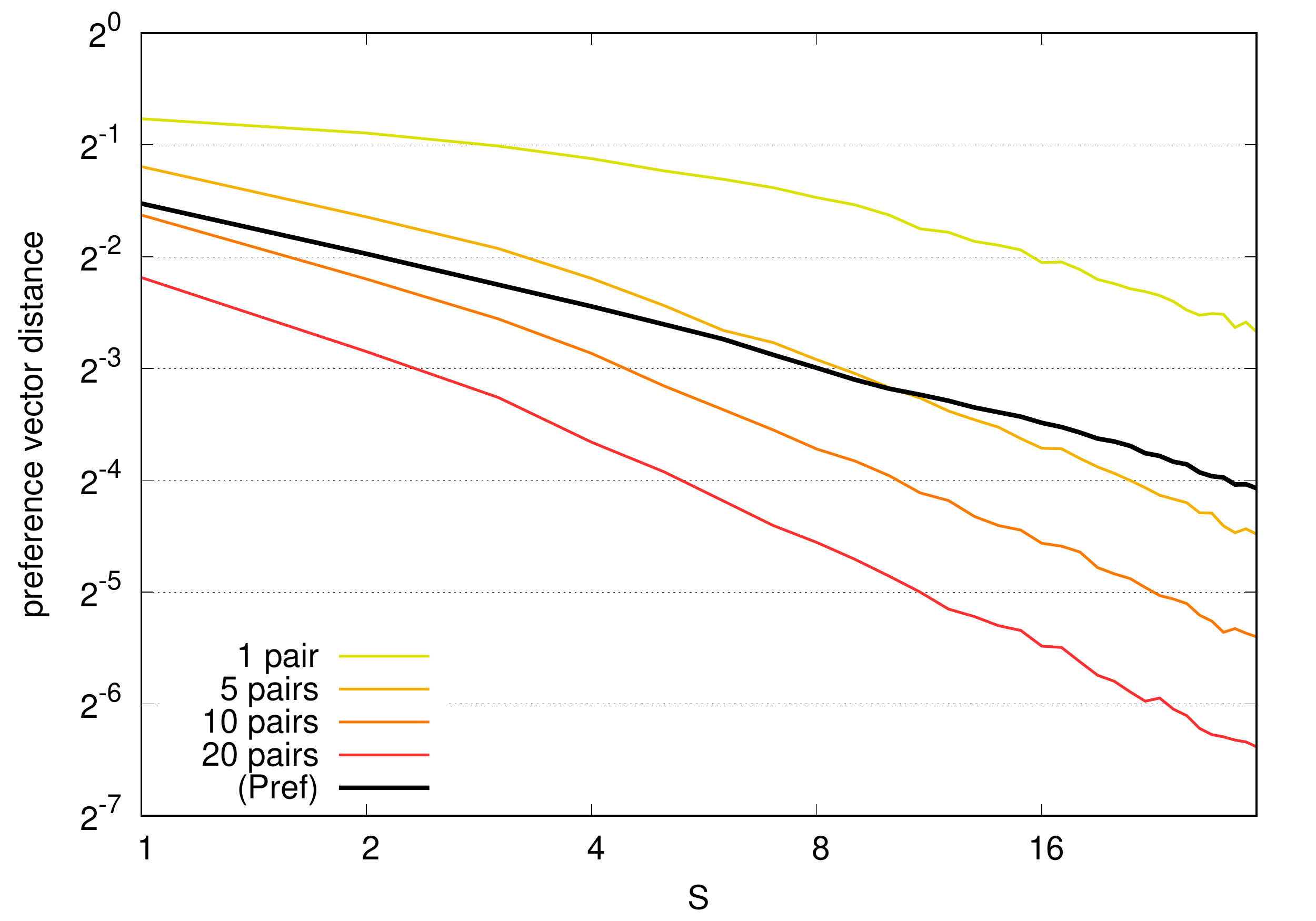}}%
\subfigure[\XInDist\label{fig::ExpSXIn}]{\includegraphics[width=.32\textwidth]{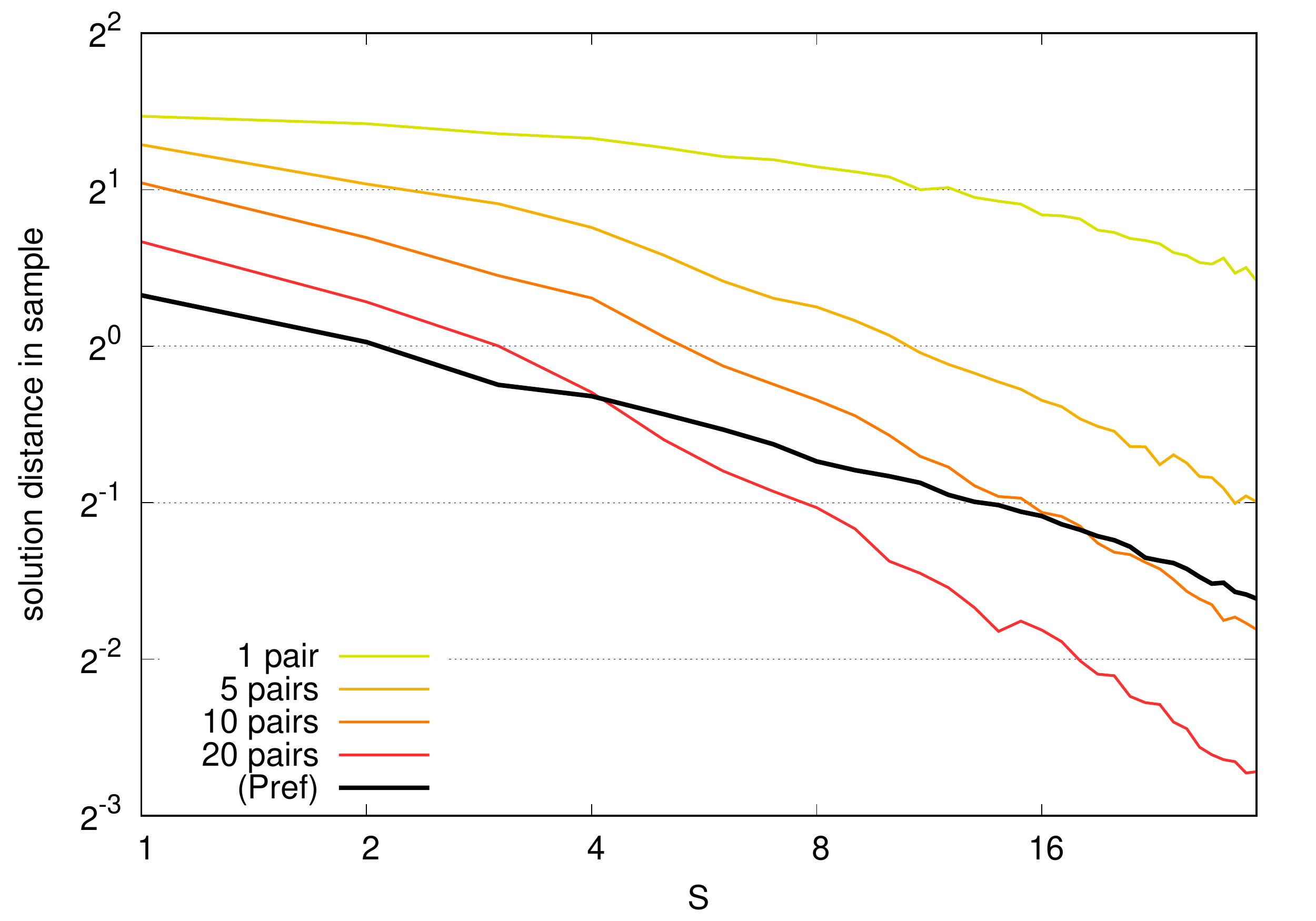}}%
\subfigure[\XOutDist\label{fig::ExpSXOut}]{\includegraphics[width=.32\textwidth]{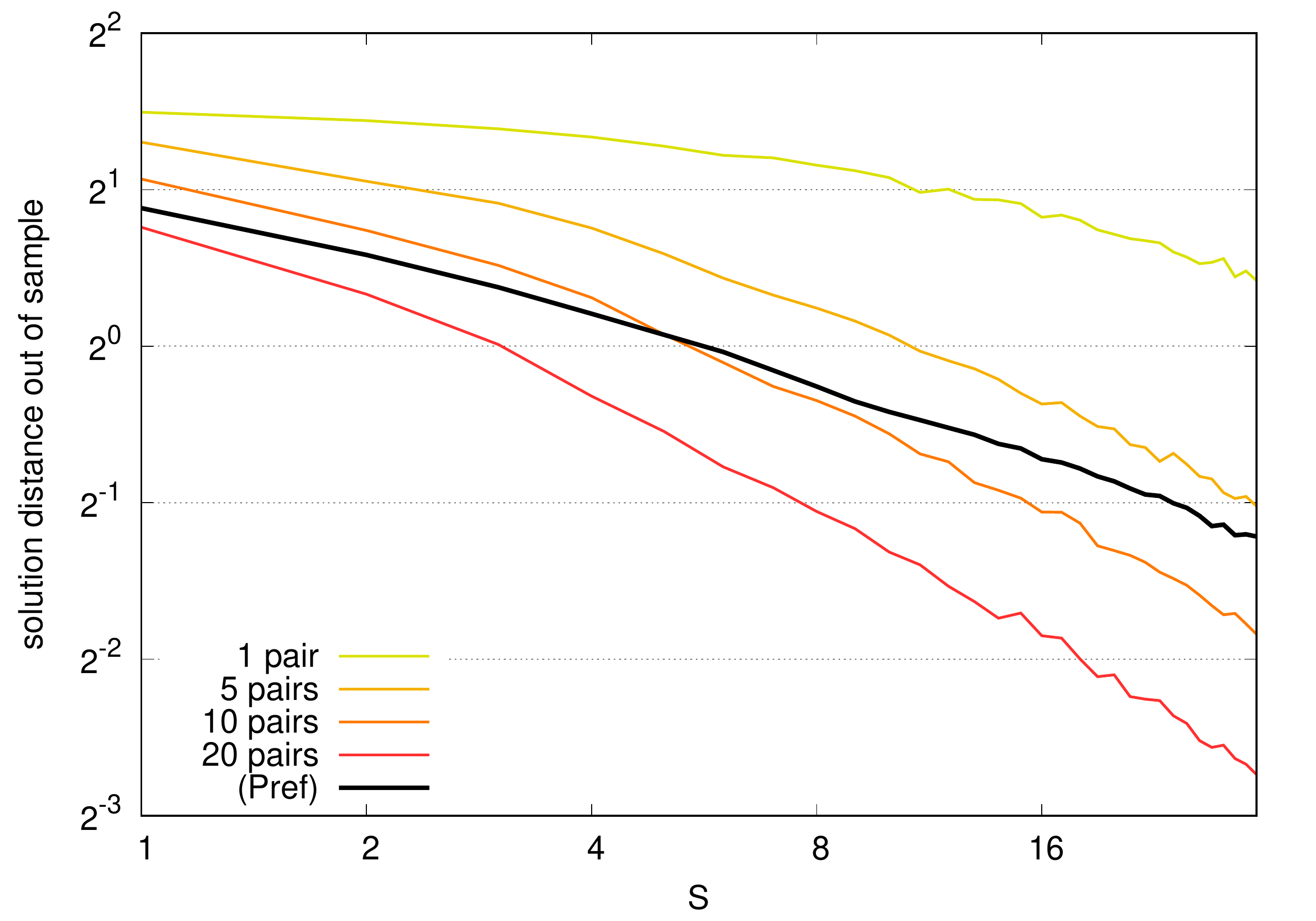}}%
\caption{Comparison of our approach with pairwise comparisons for varying number of  observations $S\in \{1,\ldots,31\}$.\label{fig::S}}
\end{figure}

For $S\in \{1,\ldots,31\}$ observations, $n=40$ and $K=5$ we also generated $1000$
random instances and conducted the same analysis. On the logarithmic horizontal axis of Figure \ref{fig::S} we show the value $S$, i.e., how many of the 31 observations we generated where
available. On the vertical axis we again show the average distance of the predicted preference vector to the actual preference vector as well as the average Hamming distance to the solutions in- and out-of-sample.

As expected, an increasing number of observations available leads to improved performance across all approaches. However, on all three performance measures \eqref{pref} does not improve at the same rate as conducting pairwise comparisons. For example in Figure~\ref{fig::ExpSXIn} it can be seen that our approach is able to recreate the observed solutions even better than using 20 comparisons when there are fewer than four observations. However, \eqref{pref} is not able to maintain this lead and for more than 30 observations the average distance to the observed solutions is similar to using ten comparisons per observation.

\begin{figure}[!htb]
\centering
\subfigure[\WDist]{\includegraphics[width=.32\textwidth]{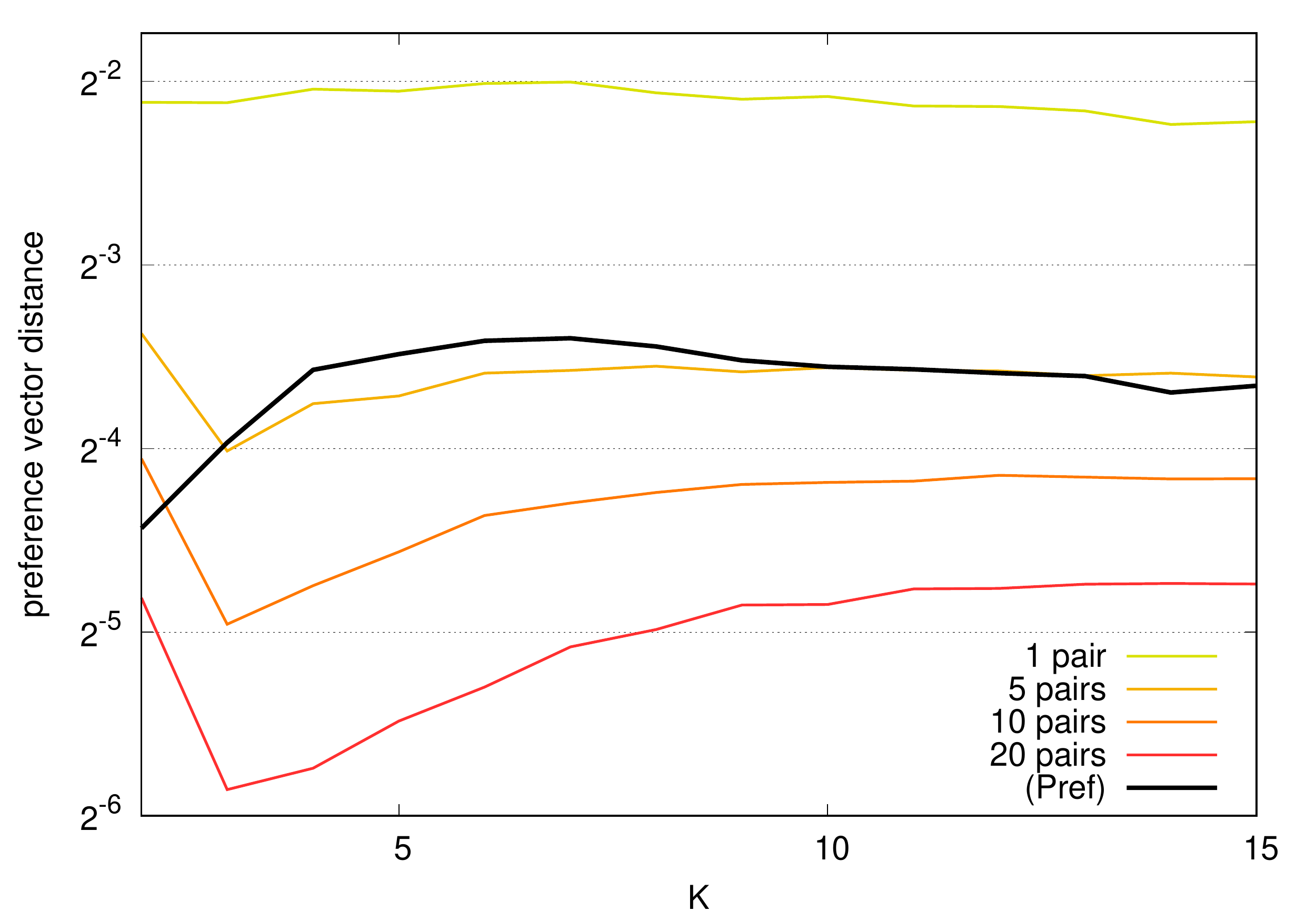}}%
\subfigure[\XInDist]{\includegraphics[width=.32\textwidth]{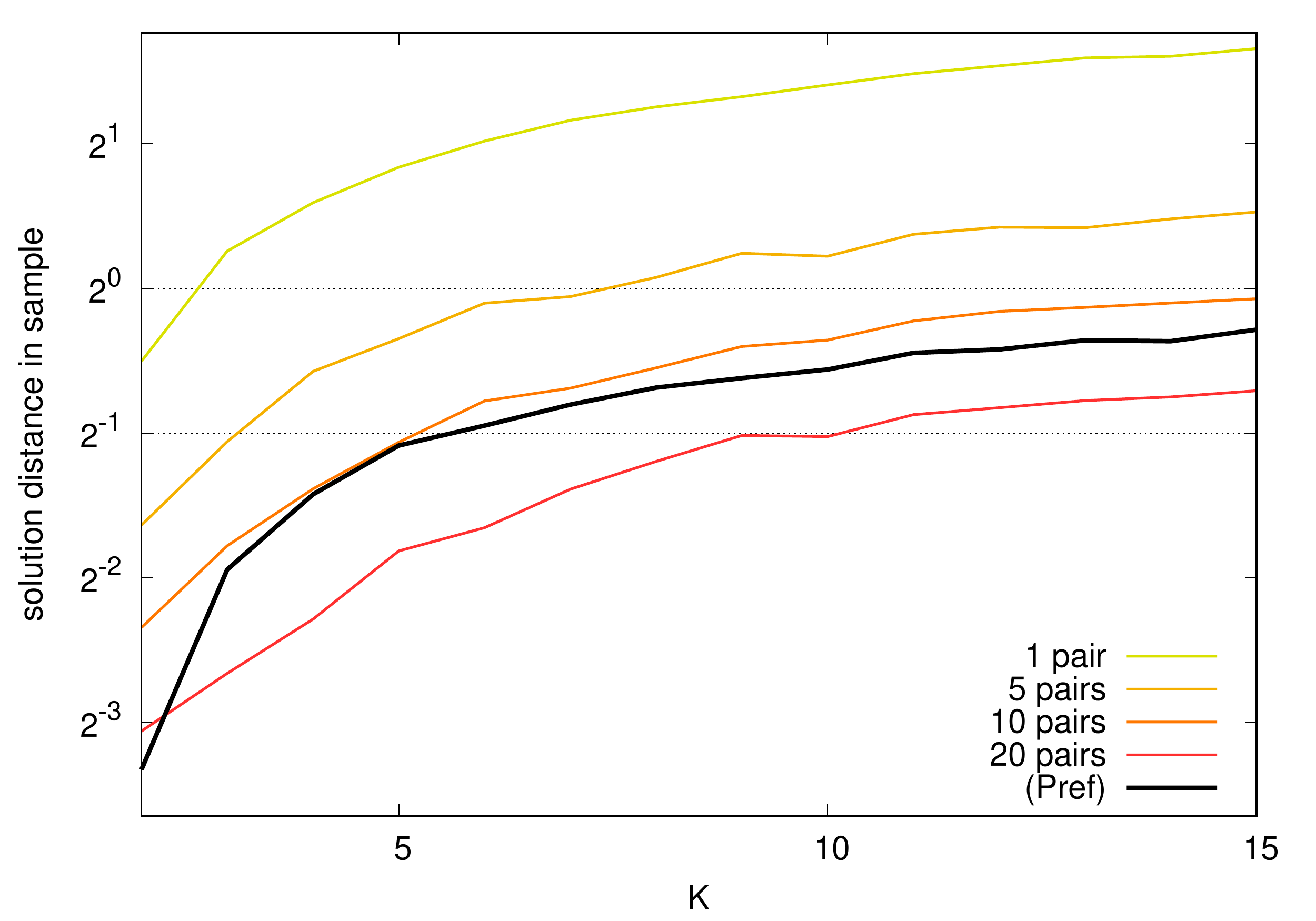}}%
\subfigure[\XOutDist]{\includegraphics[width=.32\textwidth]{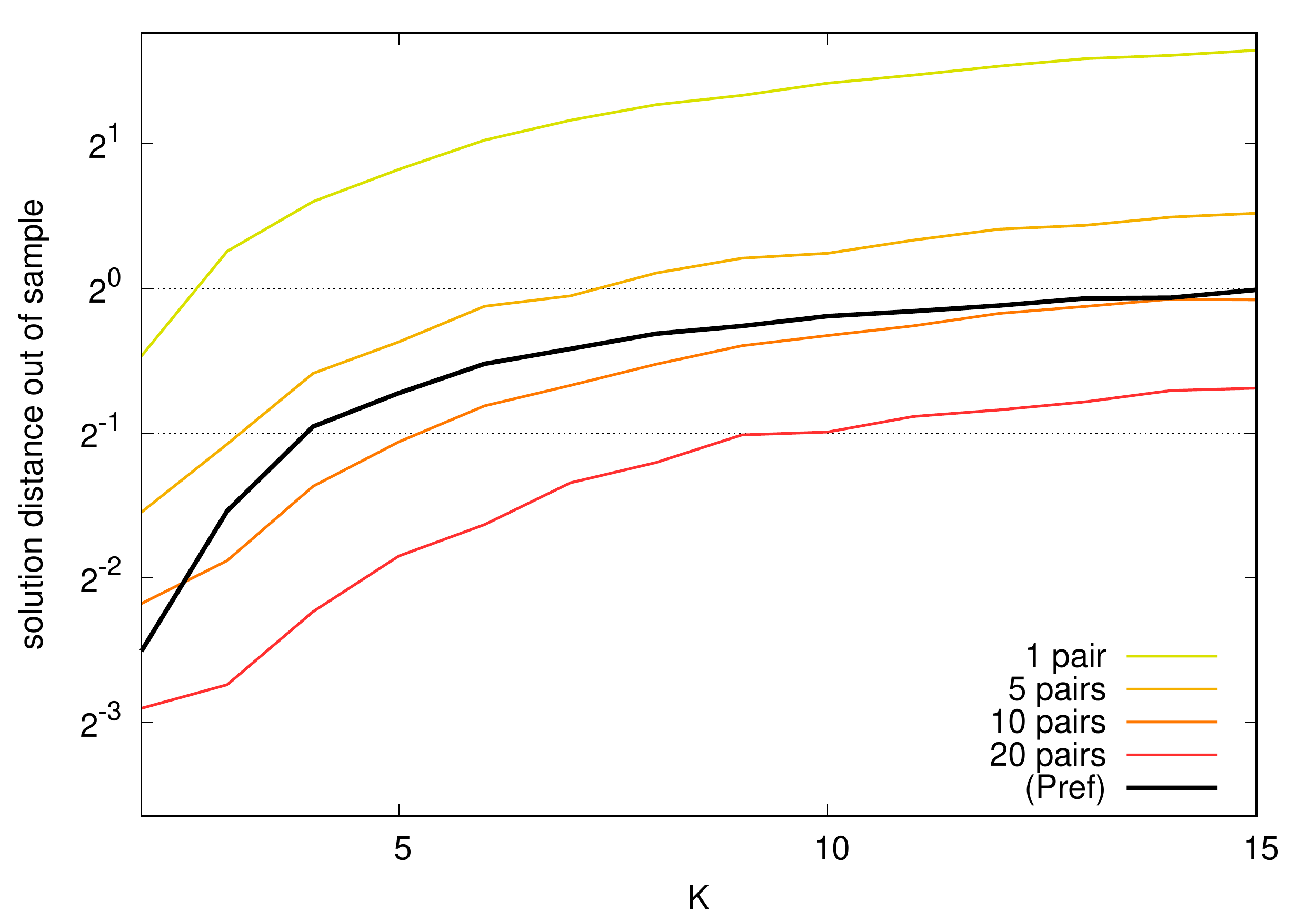}}%
\caption{Comparison of our approach with pairwise comparisons for varying $K\in \{2,\ldots,15\}$.\label{fig::K}}
\end{figure}

Finally, for $K\in\{2,\ldots, 15\}$, $n=40$ and $S=16$ the same experiments are conducted and shown in Figure \ref{fig::K}. 

For an increasing number of criteria, all three distance measures tend to increase. However, for large values of $K$, \eqref{pref} shows a slower decline in performance and regarding the preference vector distance is even able obtain better results for increasing $K$. The peculiar behavior of the approach using pairwise comparisons for very small values $K$ is likely to blame on the small pool of Pareto-optimal solutions to choose from to generate pairs.  It is noteworthy to mention that the result for $K=5$, which is the value we used in all the previous experiments, leads to the least desirable outcome when compared to experiments with varying values of $K$.

\subsubsection{Robustness against Noise}
\begin{figure}[!htb]
\centering
\subfigure[\WDist]{\includegraphics[width=.32\textwidth]{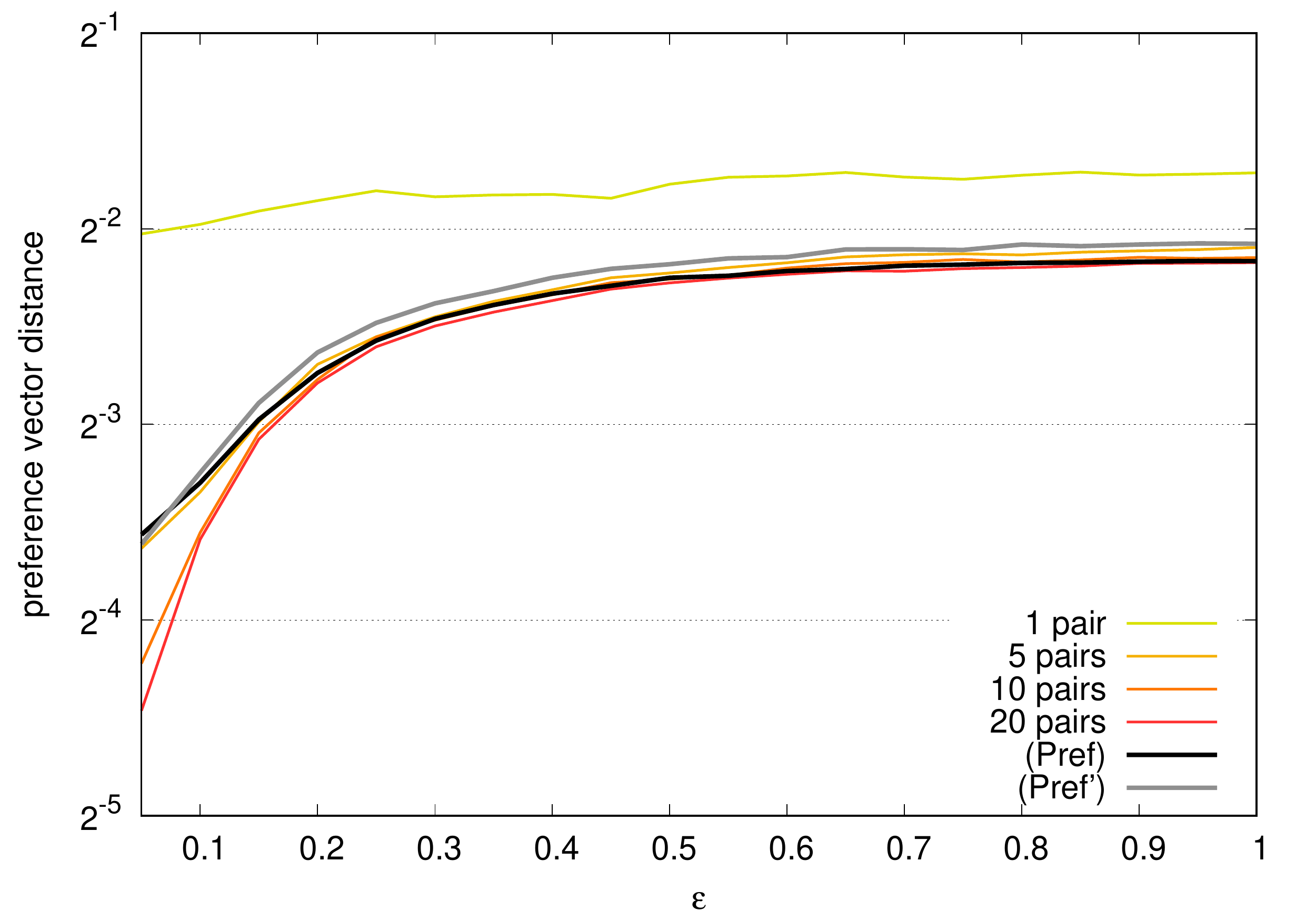}}%
\subfigure[\XInDist]{\includegraphics[width=.32\textwidth]{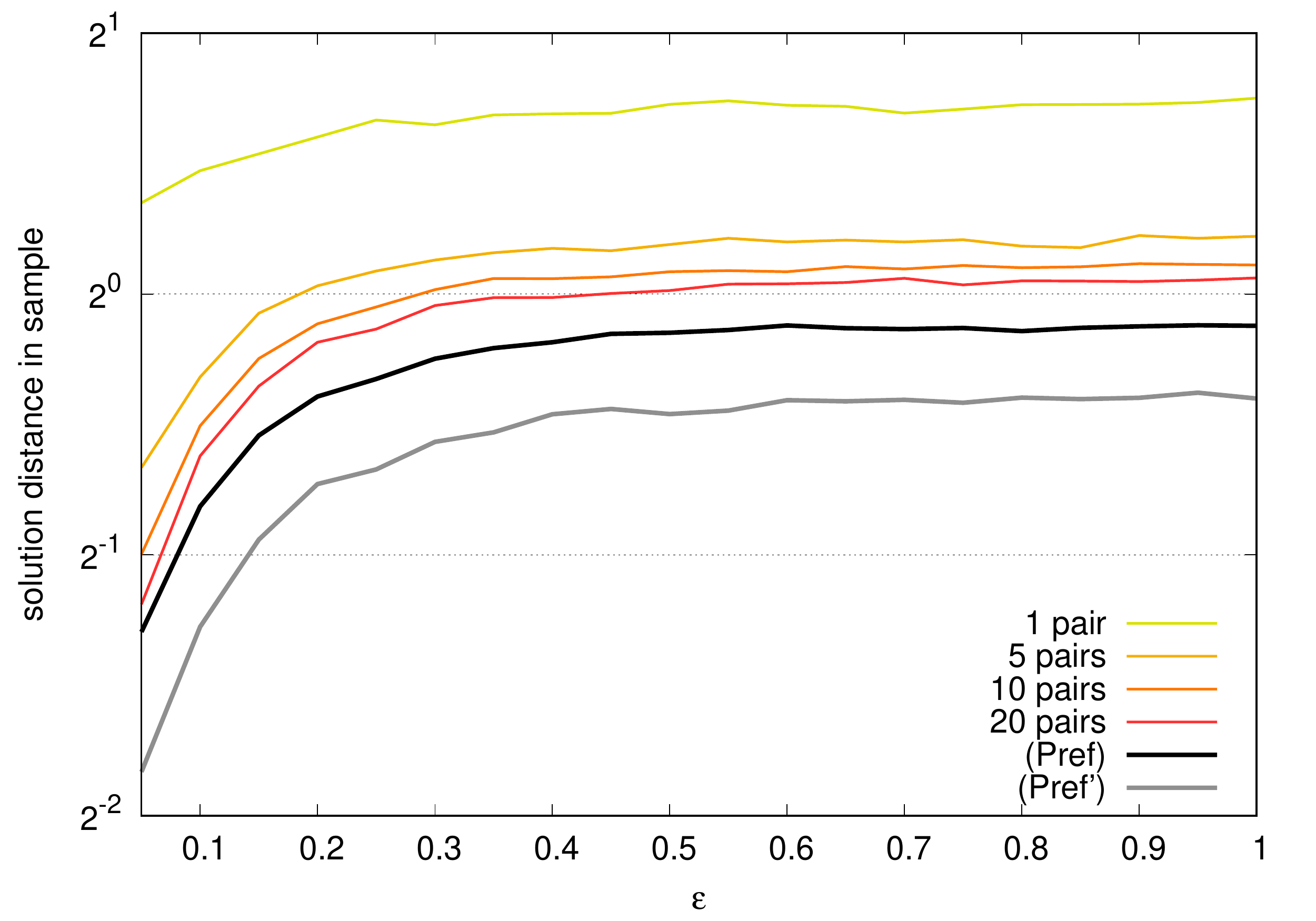}}%
\subfigure[\XOutDist]{\includegraphics[width=.32\textwidth]{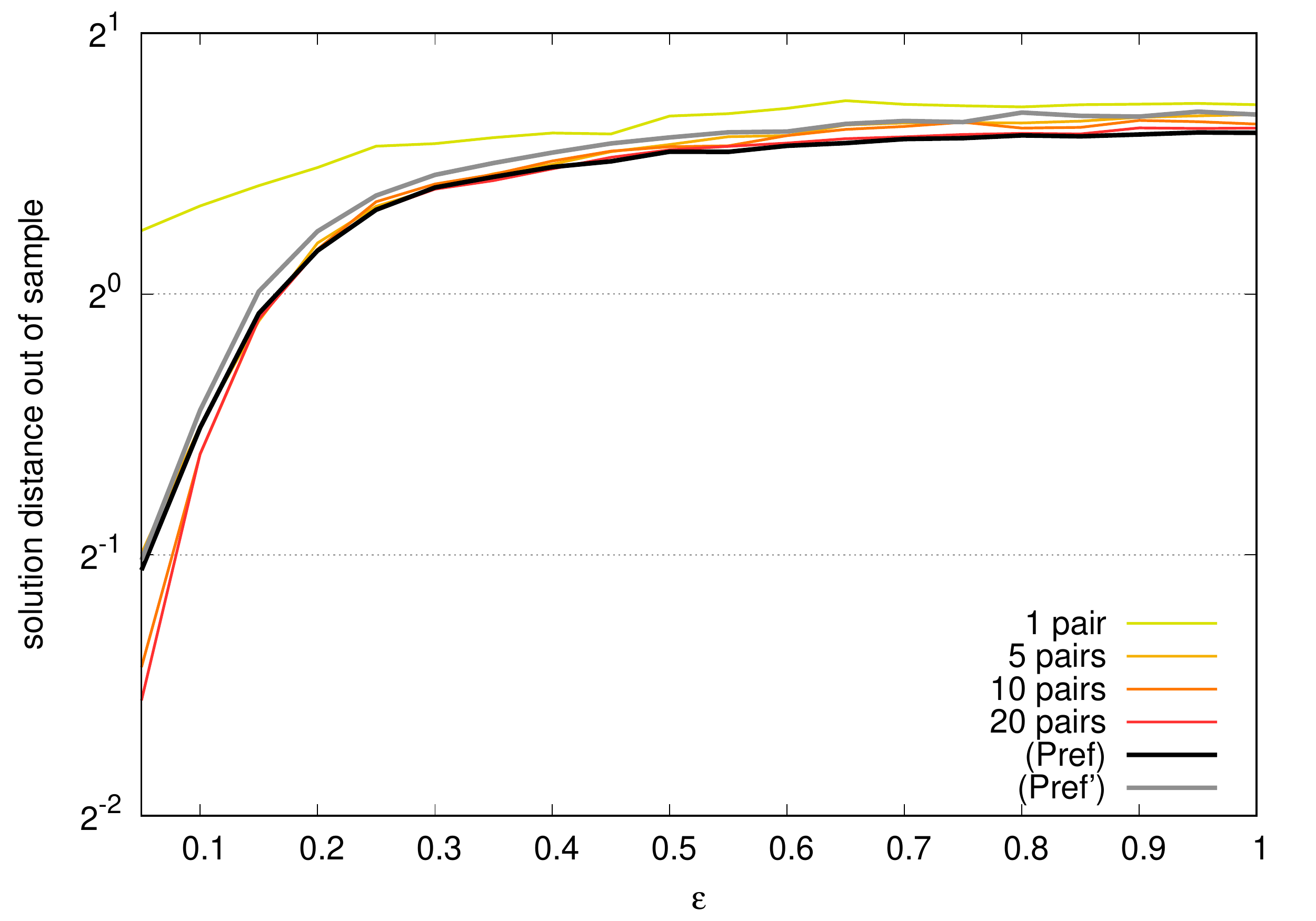}}%
\caption{Comparison of our approach with pairwise comparisons for varying $\epsilon\in \{0,0.05,\ldots,0.95,1\}$.}\label{fig::Eps}
\end{figure}

In this final experimental evaluation we want to test how well our approach behaves when there is noise in the data.  In this setting we now can expect the two proposed models to behave differently, which is why we additionally ran experiments on \eqref{altpref}. As this alternative model also contains integer variables, it exhibits an increased computational complexity making it necessary to briefly discuss handling runtimes. First, we introduced a time limit of 100 seconds for each call to (\ref{altpref}($\mathcal{X}$') (see Figure \ref{fig:alg})) in the iterative approach. In case this time limit was reached by one such call, the iterative approach stops and return the current incumbent solution. Second, in order to reduce the number of instances where this time limit is reached, we significantly decreased the problem size to $n=20$ items and $S=8$ observations, while we kept $K=5$ criteria. We created 1000 preference vectors and for each query that uses this preference vector it is altered by $\epsilon$ as explained before. Hence, each of the $S$ observed solutions as well as each pairwise comparison is obtained using a different preference vector, arising from adding random noise to the true OWA weights. In Figure~\ref{fig::Eps} the average distance to the underlying preference vector and the average solution distances are shown.  
As expected, the average distance to the true OWA weights increases for increasing $\epsilon$. For both the distance of obtained solutions compared to new out-of-sample observations as well as the distance to the true underlying preference vector all approaches behave similarly if $\epsilon$ is too large, with the exception of only using a single comparison. In both performance measures, our models lose accuracy slightly less quickly, which is likely due to being less exposed to noise compared to having every pair disturbed. When trying to rebuild the actually observed solutions both \eqref{pref} and \eqref{altpref} perform significantly better than the approach utilizing pairwise comparisons, as our approaches are based on these observations. Here also the benefit of \eqref{altpref} becomes visible as it is the best approach to recreate the observed solutions, even if they are based on modified preference vectors.

\section{Conclusions and Further Research}\label{Conclusion}

The ordered weighted averaging (OWA) criterion is a popular method to aggregate the performance of a solution over multiple objectives or scenarios. This aggregation is controlled by a vector $\pmb{w}$ which is used to express the decision maker's preference. Some well-known and often-used preference vectors are commonly used to model the worst-case, best-case, average, median, or conditional value at risk criteria. But going beyond these standard vectors, a strength of the OWA criterion is that it also allows for a more nuanced reflection of decision maker preferences. A crucial question then becomes: how can we express this preference through the vector $\pmb{w}$?

One way to approach this question is to prepare a catalog of questions to elicit preferences for a risk-averse decision maker. A drawback of such a method is that it requires interaction with the decision maker.
Other approaches assume observations of historical OWA values, but it is unclear how these should be computed without having a preference vector already available.

In this paper we propose a new and indirect approach to elicit preferences. Instead of asking questions to the decision maker, we observe her preferred choice on a set of example problems which may have come from previous rounds of decision making. Using this set of observations, we then find a preference vector that is capable to explain these choices. Our approach only requires knowledge of the chosen alternative in each observed decision making situation, and assumes no further comparisons between other alternatives, which may not be part of the historical data.

We propose an optimization model to obtain a preference vector that is at minimum distance to the polyhedra of feasible OWA weights capable of explaining the single observations. As the model has a large number of constraints, it can be approached through an iterative solution method, where we alternate between solving a linear program to determine $\pmb{w}$ and solving OWA problems to check its feasibility. Additionally, we propose a second model that finds a preference vector that can recreate the observations most accurately by minimizing the sum of the Hamming distances of the observed and proposed solutions. 

In computational experiments, we compared the performance of our models with an alternative approach from the literature that is based on pairwise preference comparisons between candidate solutions.  Based on the results of the majority of experiments, our approach is comparable or even superior compared to having to conduct around five pairwise comparisons per observation. One particular strength of our approach is its ability to recreate the observed solutions. In particular on observations resulting from varying preference vectors our approach performs particularly well, showcasing its robustness against an inconsistent decision maker. The observed results were consistent even across different problem structures, such as the selection, assignment and min-knapsack problem.

In further research, a study to test our preference elicitation approach with real-world decision makers would be a valuable addition, see, e.g. the recent study by \cite{reimann2017well}. Note that such an experiment is not trivial, as the ''true'' preference of a decision maker cannot be determined. Furthermore, our philosophy may be applied to other decision making criteria. In particular the weighted ordered weighted averaging (WOWA) criterion seems a natural choice as a generalization of the OWA criterion considered in this paper, see \cite{Ogryczak2009915}.

\newpage
\appendix

\section{Connections Between Models \eqref{pref} and \eqref{altpref}}\label{App::ExampleAlternative}

We begin with the following observation: A weight vector $\pmb{w}\in\cW'$ is an optimal solution to Problem~\eqref{pref} with objective value zero, if and only if there are $\pmb{y}^s\in \opt'_s(\pmb{w})$ for all $s\in[S]$ such that $(\pmb{w},(\pmb{y}^s)_{s\in S})$ is an optimal solution to \eqref{altpref} with objective value zero.

To see why this is true, let us first assume that $\pmb{w}\in\cW'$ is an optimal solution to Problem~\eqref{pref} with objective value zero. This means that $\sum_{s\in[S]} D(\pmb{w},\opt_s) = 0$. As $D$ is non-negative, it follows that $D(\pmb{w},\opt_s)=0$ for all $s\in [S]$, and hence, $\pmb{w}\in\opt_s$ for all $s\in [S]$. By definition of $\opt_s$, this means that each $\pmb{x}^s$ is an optimizer of  $\OWA_{\pmb{w}}(\pmb{x},C^s)$. By setting $\pmb{y}^s = \pmb{x}^s$ for all $s\in [S]$, we have $\pmb{y}^s \in \opt'_s(\pmb{w})$ and $d'(\pmb{y}^s,\pmb{x}^s) = 0$. Therefore, $(\pmb{w},(\pmb{y}^s)_{s\in S})$ is an optimal solution to \eqref{altpref} with objective value zero.

For the other direction, let us assume that $(\pmb{w},(\pmb{y}^s)_{s\in S})$ is an optimal solution to \eqref{altpref} with objective value zero, that is, $\sum_{s\in[S]} d'(\pmb{y}^s,\pmb{x}^s)=0$. As $d$ is a metric, this is true if and only if $\pmb{x}^s=\pmb{y}^s$ for all $s\in[S]$. As $\pmb{y}^s\in\opt'_s(\pmb{w})$, it follows that $\pmb{x}^s\in\opt'_s(\pmb{w})$ and therefore $\pmb{w}\in\opt_s$. We conclude that $\sum_{s\in[S]} D(\pmb{w},\opt_s)=0$.

However, if there is no single $\pmb{w}$ that can explain all observations, optimal solutions may differ, as the following example illustrates.

We consider the selection problem where one has to select $2$ out of $4$ items. There are $K=3$ objectives and $S=2$ observations. In this case the two observed solutions stem from two different preference vectors, i.e.~an inconsistent decision maker. The objective coefficients of the objectives as well as the two observed solutions are given in the upper part of Table \ref{tab::ExampleAlternative}.

\begin{table}[htb]
\caption{Data for the example. \label{tab::ExampleAlternative}}
\centering
\begin{tabular}{rrr}
\toprule
&Observation 1 & Observation 2\\\midrule
Objective 1 &$\pmb{c}^{1,1}=\left(0.6,1,0.5,0\right)$&$\pmb{c}^{2,1}=\left(0.8,1,0.3,0\right)$\\
Objective 2&$\pmb{c}^{1,2}=\left(1,0.7,0,0.3\right)$&$\pmb{c}^{2,2}=\left(0.8,0,1,0.6\right)$\\
Objective 3&$\pmb{c}^{1,3}=\left(0,0,0.8,1\right)$&$\pmb{c}^{2,3}=\left(0.8,1,0.1,0\right)$\\
Preference vector&$\left(0.35,0.33,0.32\right)$&$\left(0.54,0.33,0.13\right)$\\
Optimal solution &$\pmb{x}^1=\left(0,0,1,1\right)$&$\pmb{x}^2=\left(0,1,0,1\right)$\\
\midrule
\eqref{pref}\\
Preference vector &$\left(0.46,0.29,0.25\right)$\\
Solution &$(1,0,1,0)$&$(0,0,1,1)$\\
\midrule
\eqref{altpref}\\
Preference vector &$\left(1,0,0\right)$\\
Solution &$(0,1,0,1)$&$(0,1,0,1)$\\
\bottomrule
\end{tabular}
\end{table}

In this example, all optimal $\pmb{x}$-solutions are unique. While Model \eqref{pref} does not explicitly know the preference vectors used to create the solutions in the two observations, this example shows that the returned preference vector of \eqref{pref} indeed manages to find a preference vector that is rather close to the two underlying vectors. \eqref{altpref}, on the other hand aims at rebuilding the solutions of the observations. The preference vector returned by \eqref{altpref} has no similarities to the underlying preference vectors, but is able to rebuild one of the solutions and only has a cumulated Hamming distance of 2, while the solutions returned by \eqref{pref} have have a cumulated Hamming distance of 4 to the observed solutions. 

\section{Model to Generate OWA Weights with Given Orness\label{App::OWAOrness}}
This model was used by \cite{wang2005minimax} and creates a OWA preference vector $\pmb{w} \in [0,1]^K$ with given orness $\alpha$, where the auxiliary variable $\delta$ linearizes the maximum disparity $\max_{k \in \{1,\ldots,K-1\}} |w_k-w_{k+1}|$ between two adjacent weights.
\begin{align*}
\min \ & \delta \\
\textnormal{s.t.}\ & \alpha = \frac{1}{K-1}\sum_{k\in [K]}(K-k)w_k\\
& -\delta \leq w_k - w_{ k+1} \leq \delta \quad \forall k\in{1,\ldots,K-1}\\
& \pmb{w} \in \cW' \\
&\delta\geq 0
\end{align*}

\section{Model for Obtaining an OWA Weights that Minimize Violations\label{App::OwaAhn}}
This model was proposed by \cite{ahn2008preference} and aims at finding OWA weights $\pmb{w} \in \cW'$ that are consistent with the
decision-maker's judgments between alternatives $(\pmb{x},\pmb{y}) \in \X \times \X$, where $\pmb{x}$ is the preferred alternative. Let $\Theta \subseteq \X \times \X$ be the set of such ordered pairs and let $a_1(\pmb{x}),\ldots,a_K(\pmb{x})$ be the objective values sorted from largest to smallest of solution $\pmb{x} \in \X$.
\begin{align*}
\min \ & \sum_{(\pmb{x},\pmb{y}) \in \Theta} \delta_{\pmb{x},\pmb{y}} \\
\textnormal{s.t.}\ & \sum_{k\in [K]}\left( a_k(\pmb{y})- a_k(\pmb{x})\right)w_k + \delta_{\pmb{x},\pmb{y}} \geq \epsilon \quad \forall (\pmb{x},\pmb{y}) \in \Theta \\
& \pmb{w} \in \cW' \\
&\delta_{\pmb{x},\pmb{y}}\geq 0 \quad \forall (\pmb{x},\pmb{y}) \in \Theta \\
\end{align*}
Here, $\epsilon > 0$ is a small constant to express strict preference between two alternatives.

\section{Experiments on Assignment and Knapsack Instances}\label{Sec::AssignmentKnapsack}

\begin{figure}[!htb]
\centering
\rotatebox{90}{\hspace{.7cm} Assignment}
\subfigure[\WDist]{\includegraphics[width=.32\textwidth]{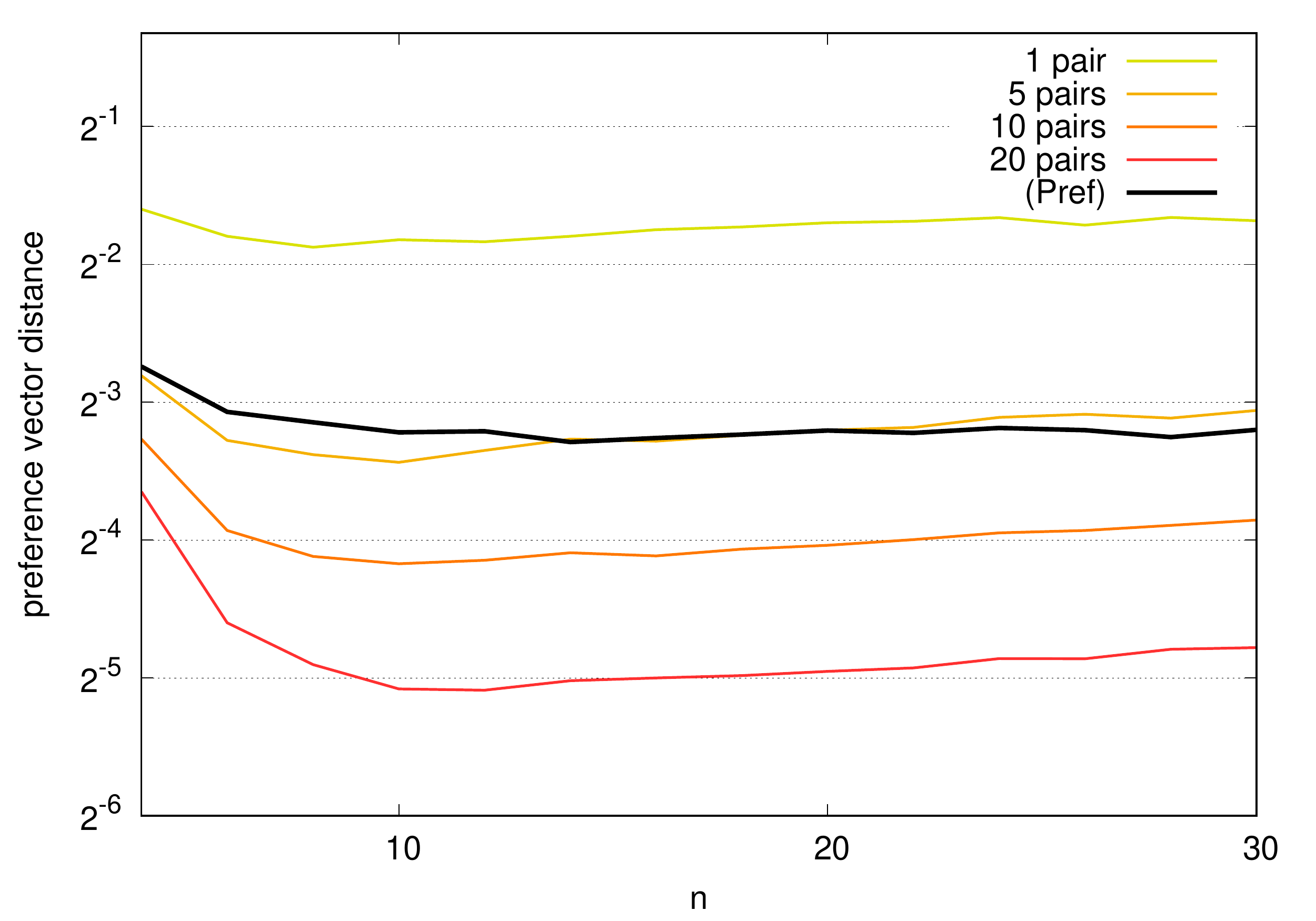}}%
\subfigure[\XInDist]{\includegraphics[width=.32\textwidth]{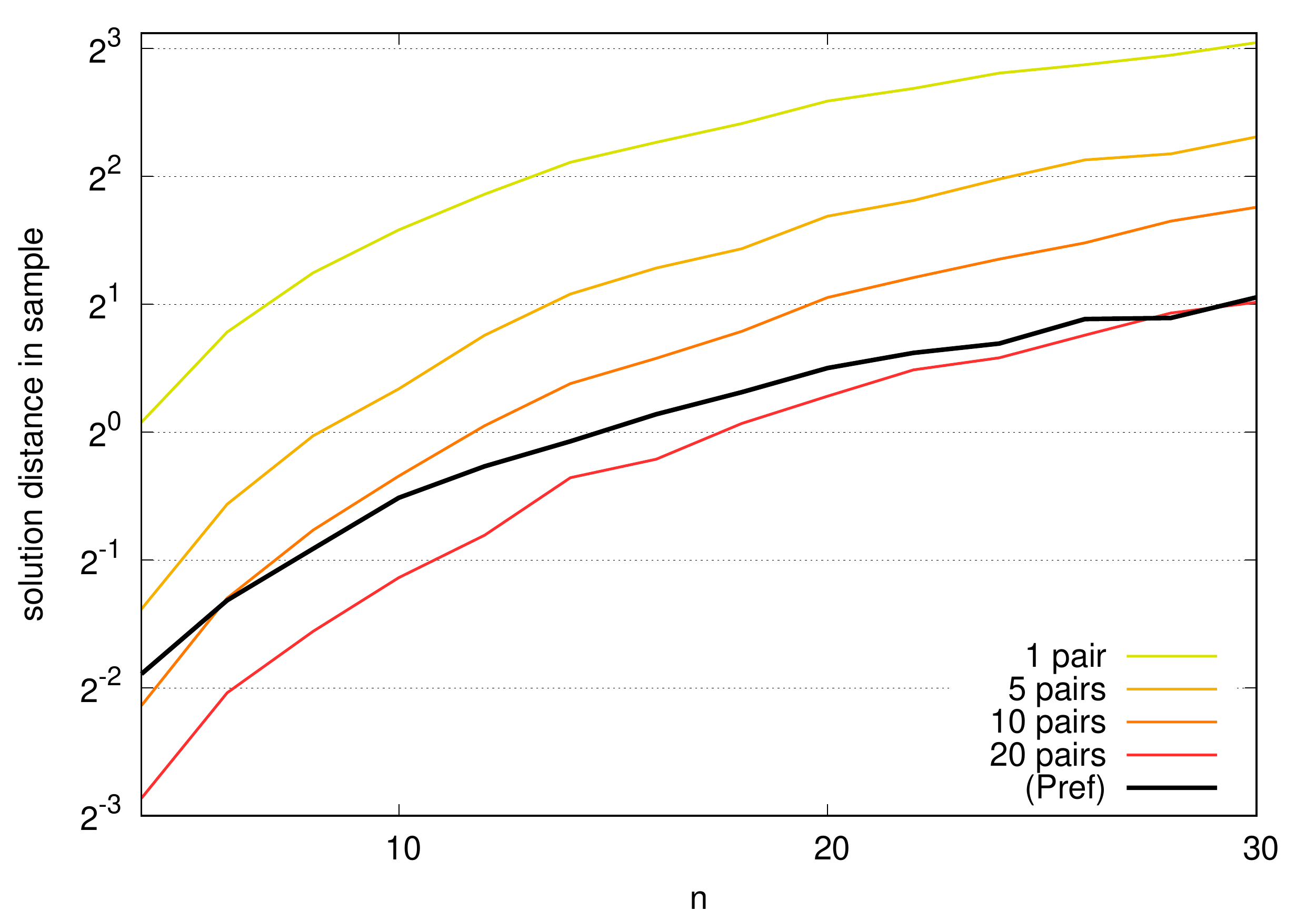}}%
\subfigure[\XOutDist]{\includegraphics[width=.32\textwidth]{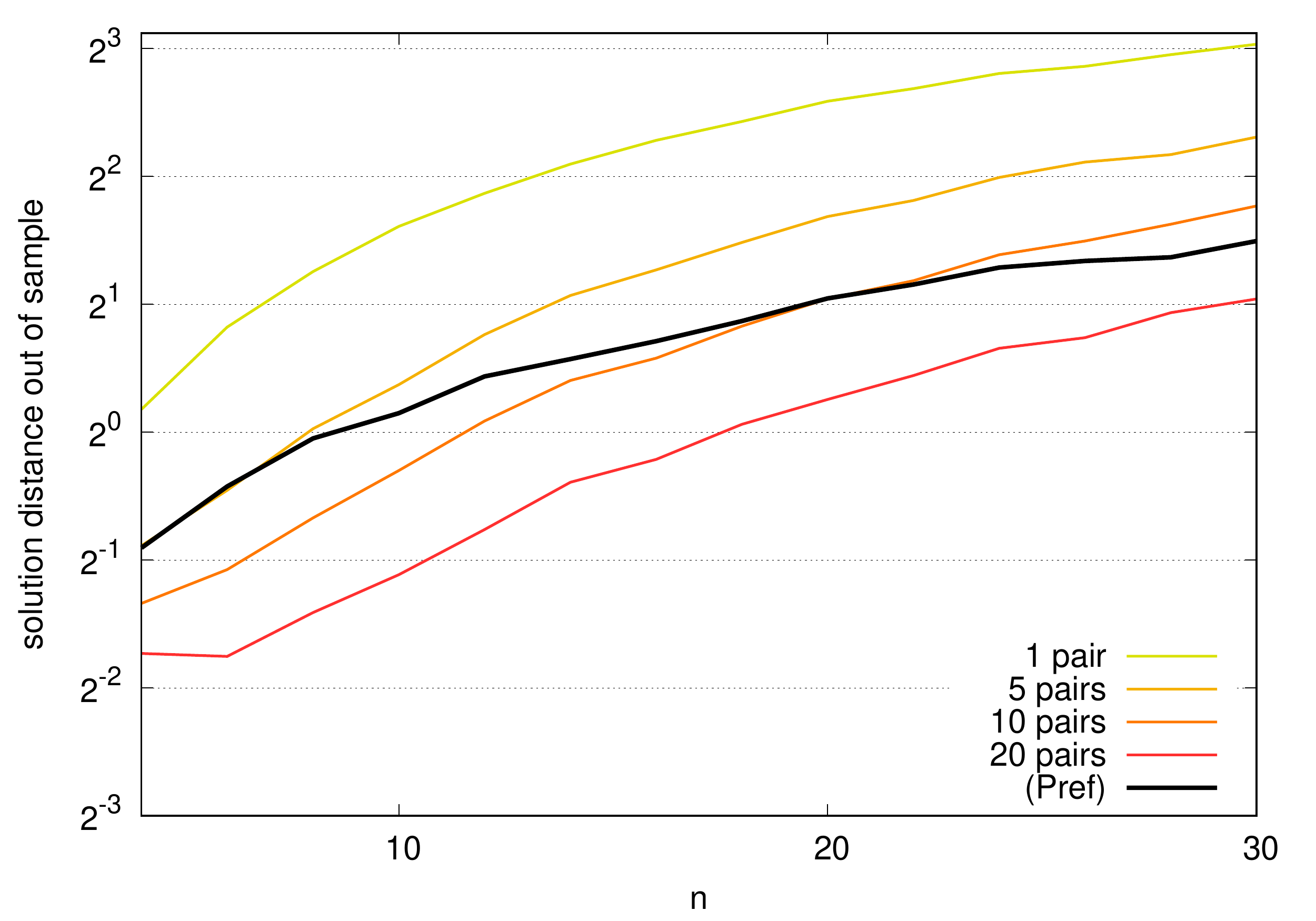}}%

\rotatebox{90}{\hspace{.9cm} Knapsack}
\subfigure[\WDist]{\includegraphics[width=.32\textwidth]{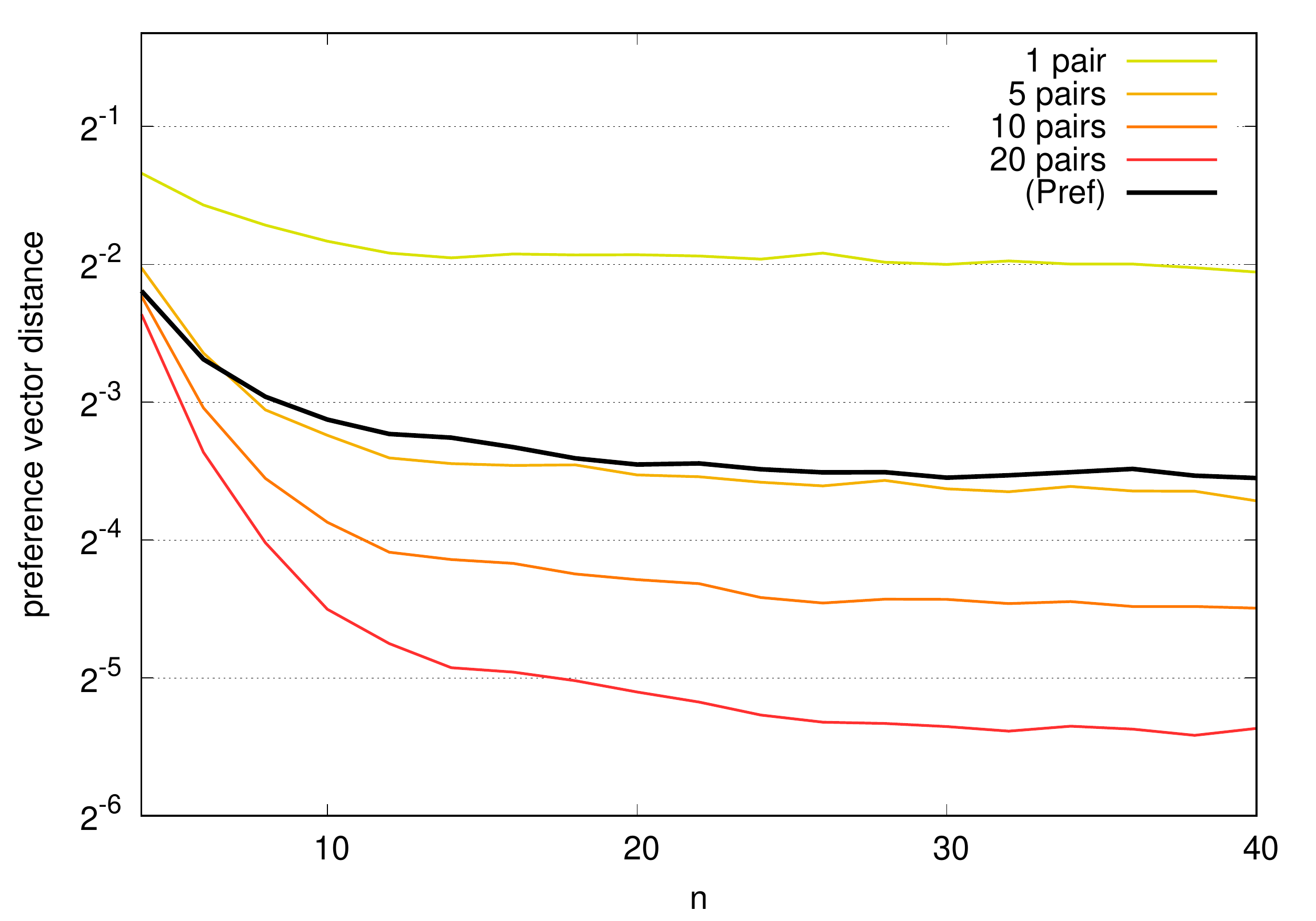}}%
\subfigure[\XInDist]{\includegraphics[width=.32\textwidth]{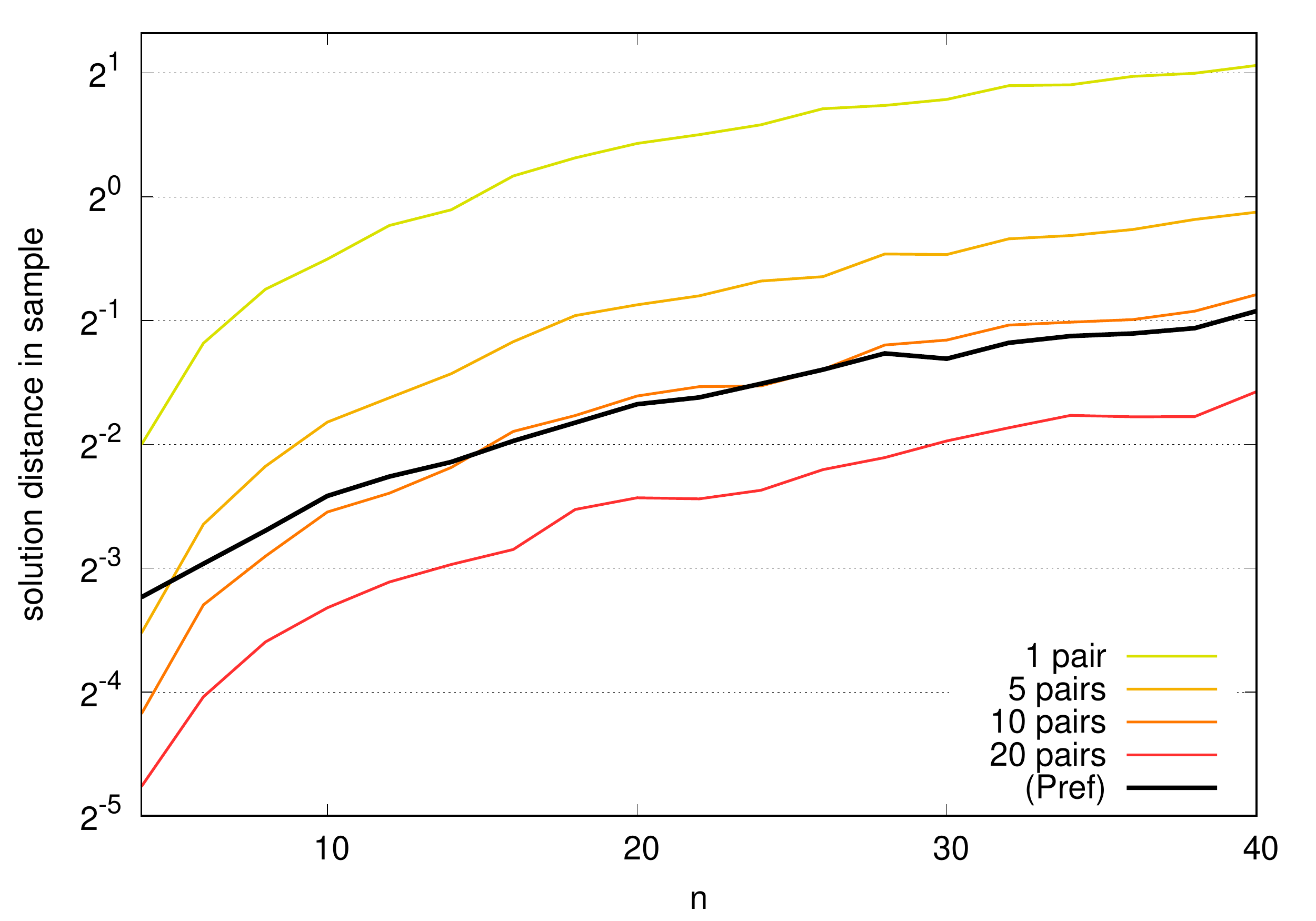}}%
\subfigure[\XOutDist]{\includegraphics[width=.32\textwidth]{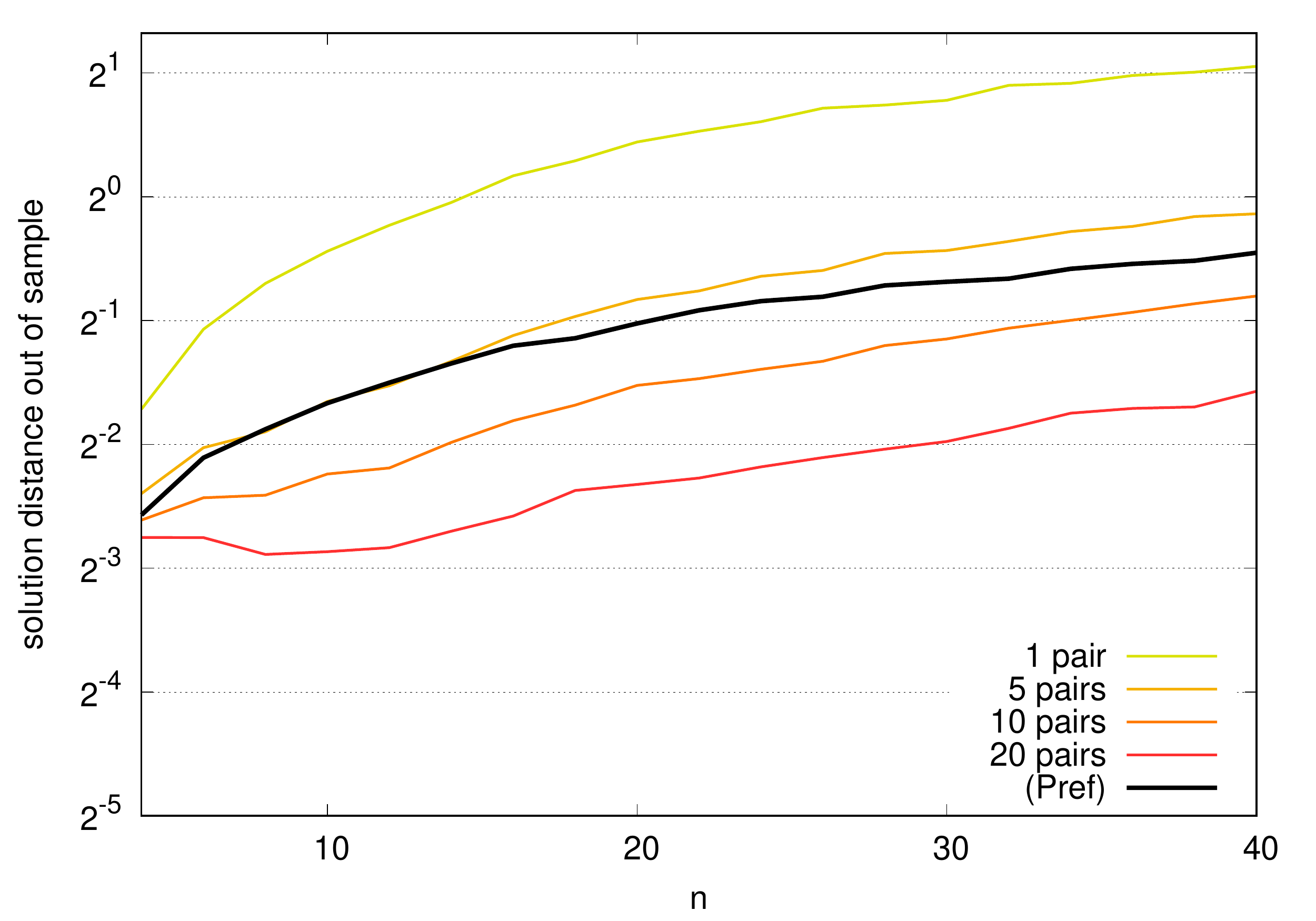}}%
\caption{Comparison of our approach with pairwise comparisons for varying problem size number $n$  for assignment (top) and knapsack (bottom) instances, with $S=16$ and $K=5$. \label{fig::Full_N}}
\end{figure}

\begin{figure}[!htb]
\centering
\rotatebox{90}{\hspace{.7cm} Assignment}
\subfigure[\WDist]{\includegraphics[width=.32\textwidth]{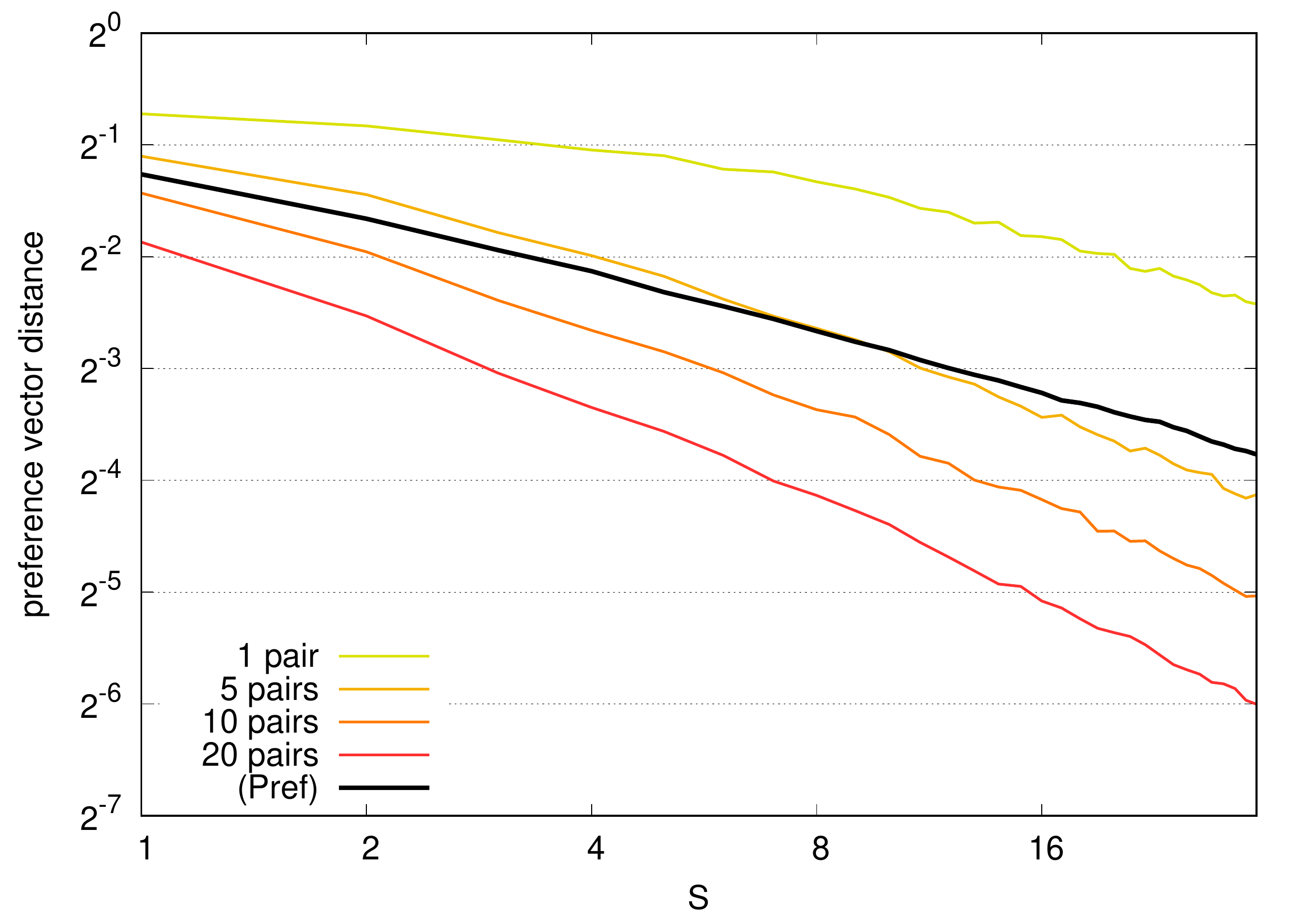}}%
\subfigure[\XInDist]{\includegraphics[width=.32\textwidth]{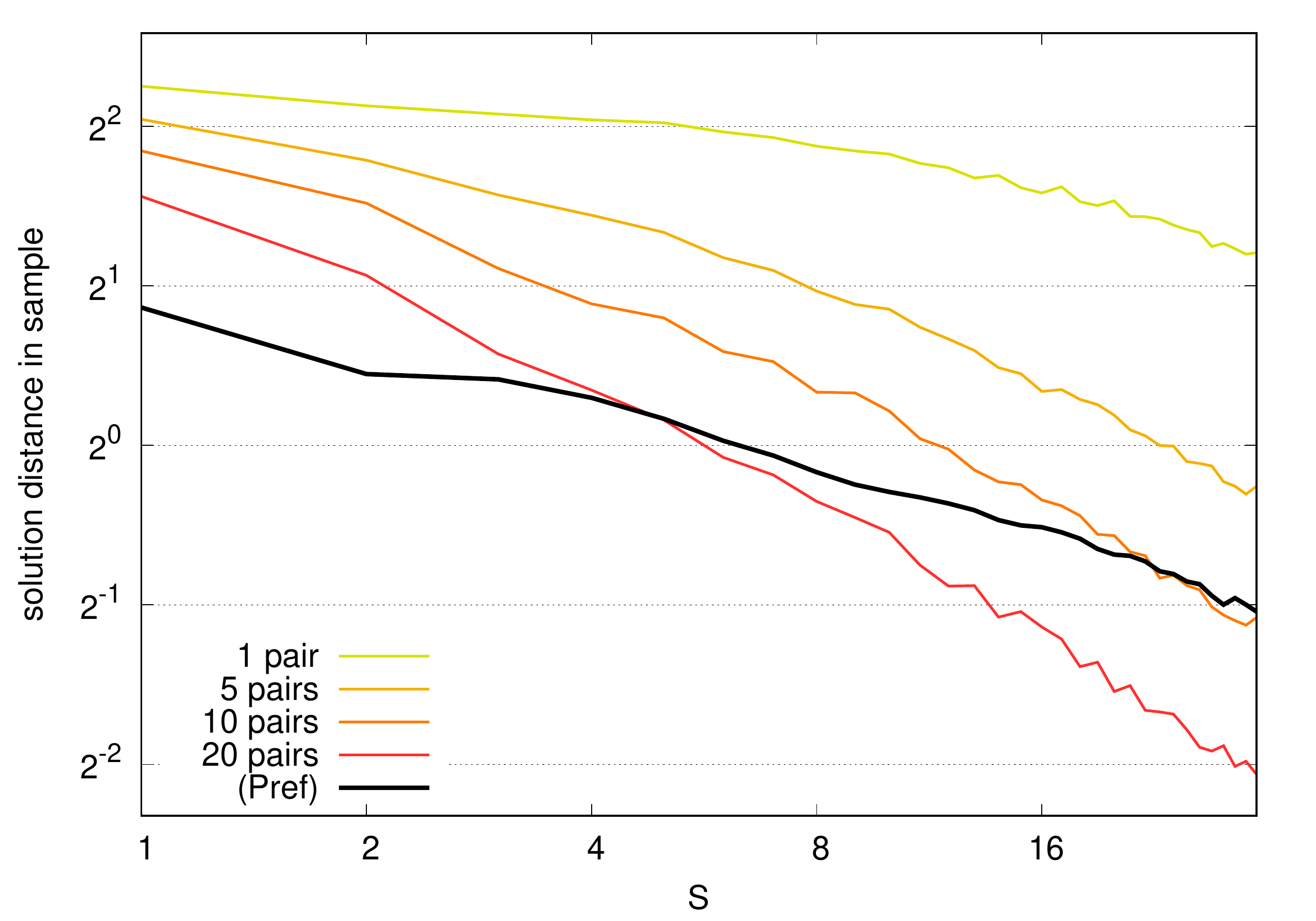}}%
\subfigure[\XOutDist]{\includegraphics[width=.32\textwidth]{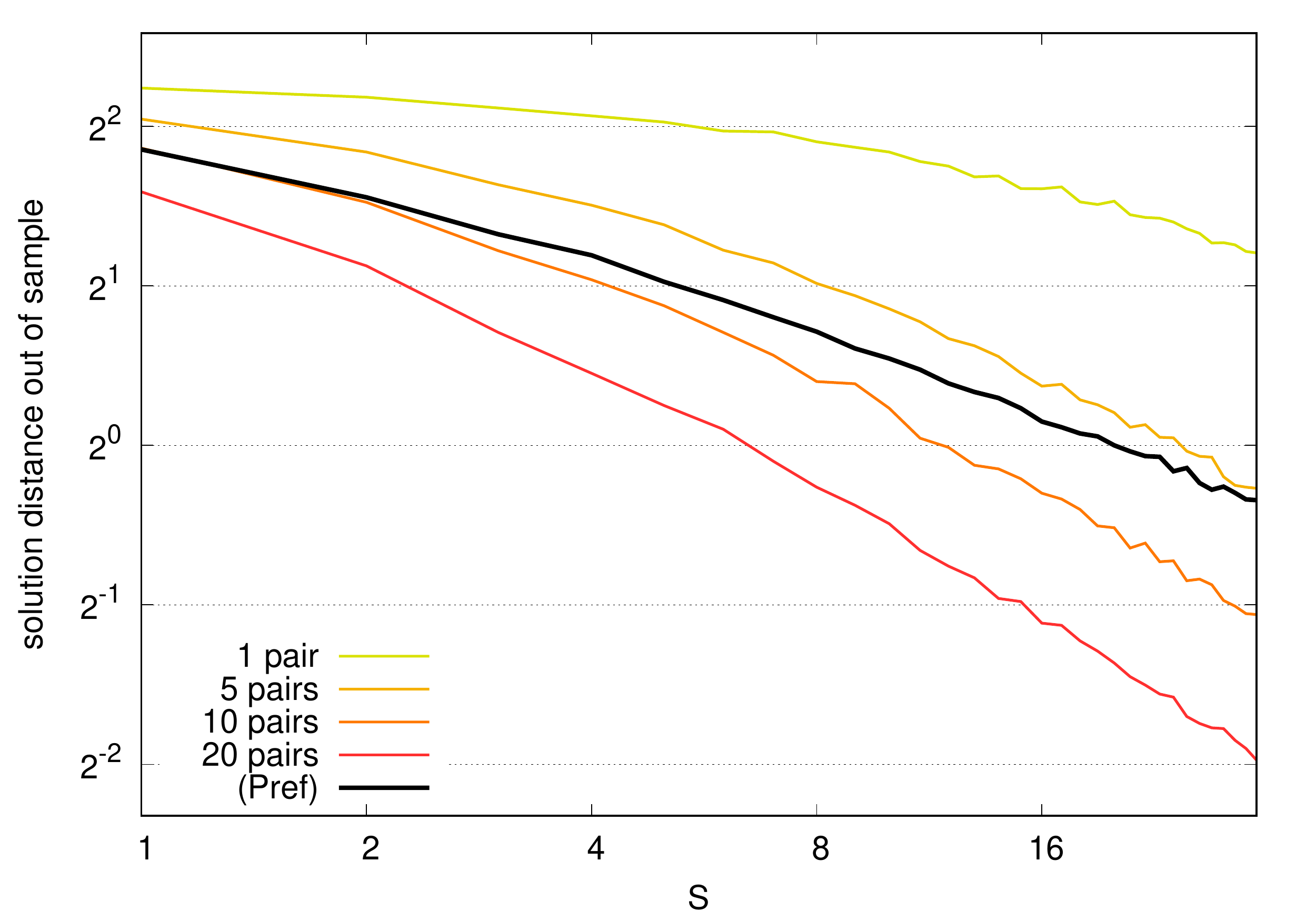}}%

\rotatebox{90}{\hspace{.9cm} Knapsack}
\subfigure[\WDist]{\includegraphics[width=.32\textwidth]{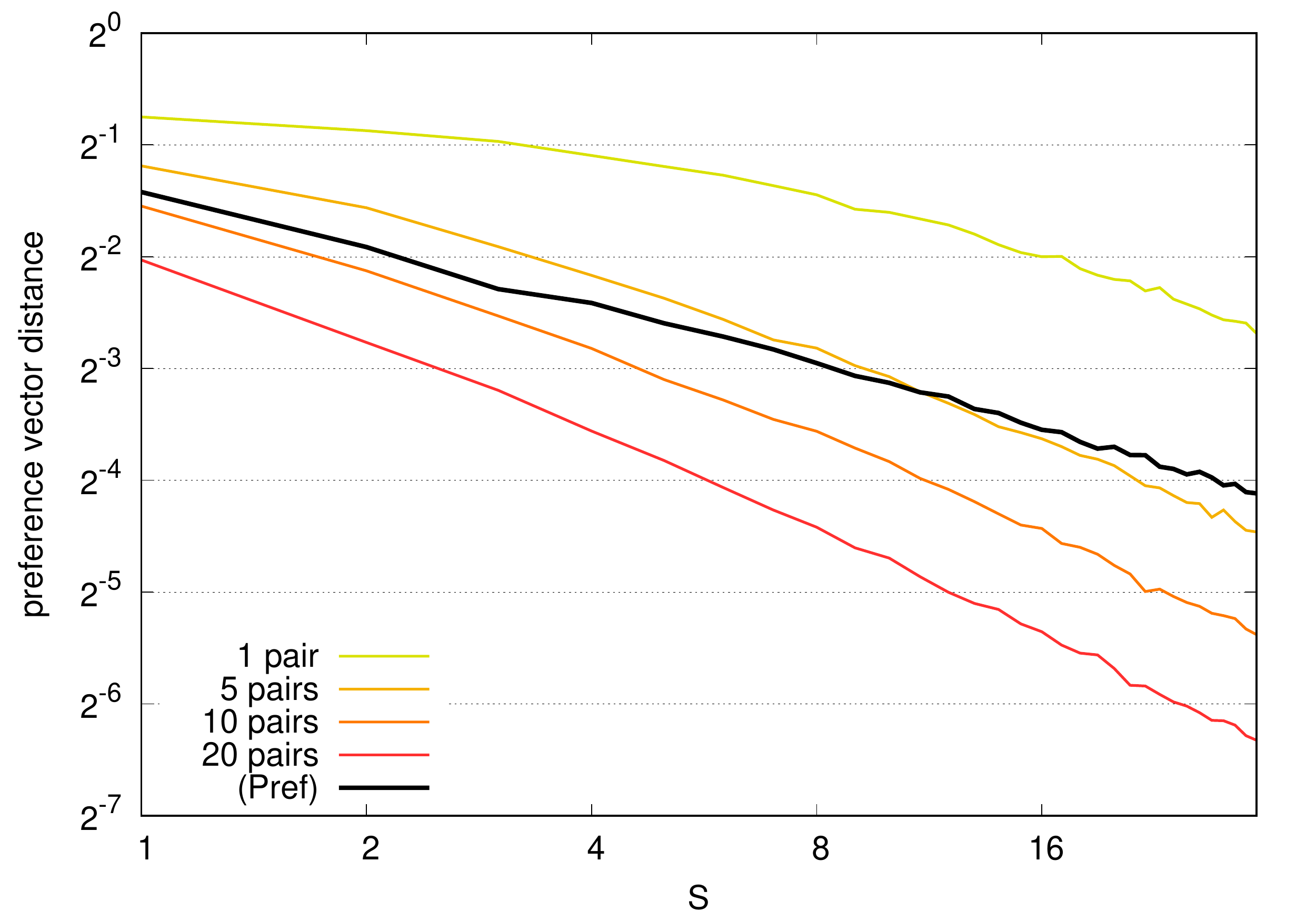}}%
\subfigure[\XInDist]{\includegraphics[width=.32\textwidth]{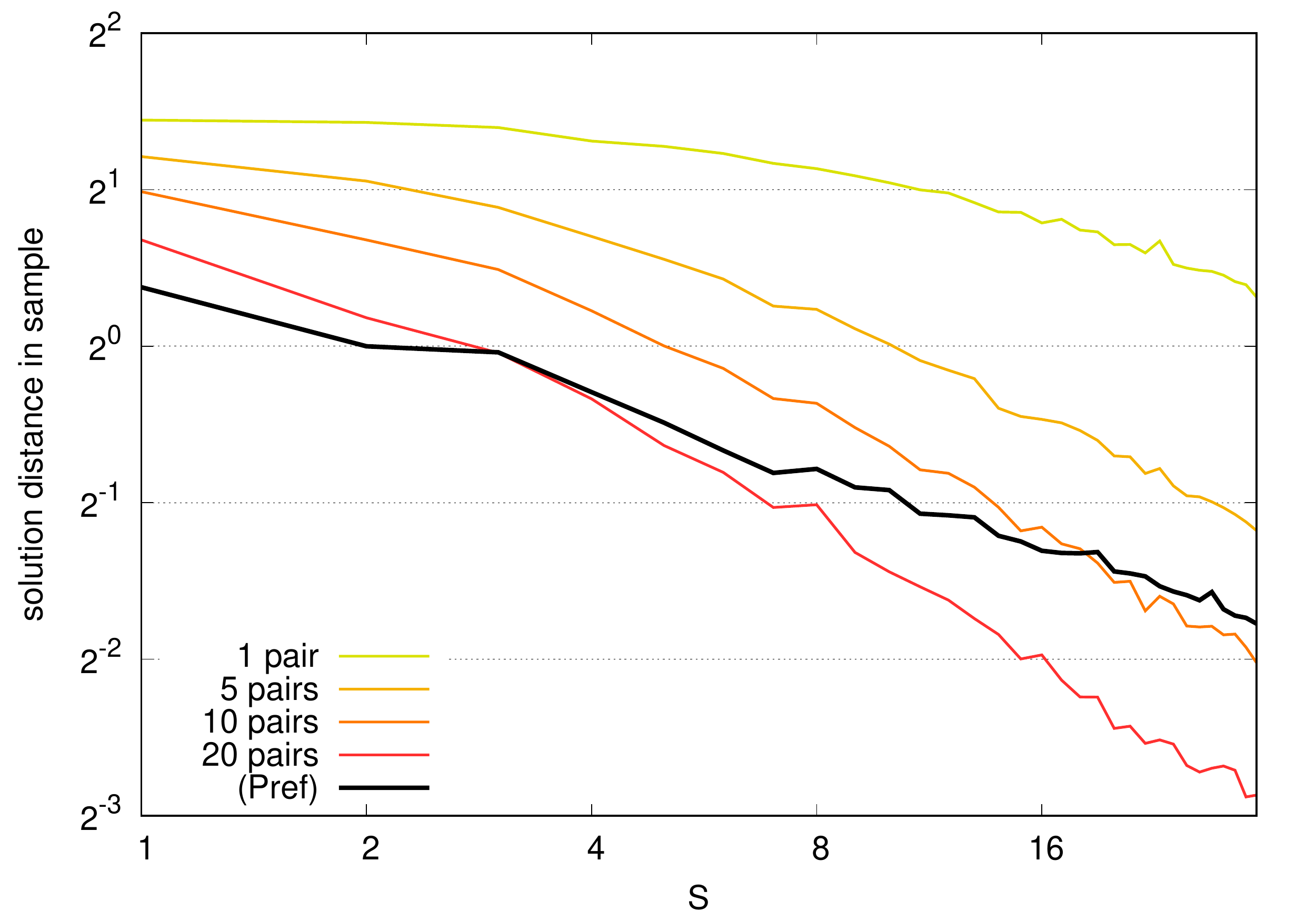}}%
\subfigure[\XOutDist]{\includegraphics[width=.32\textwidth]{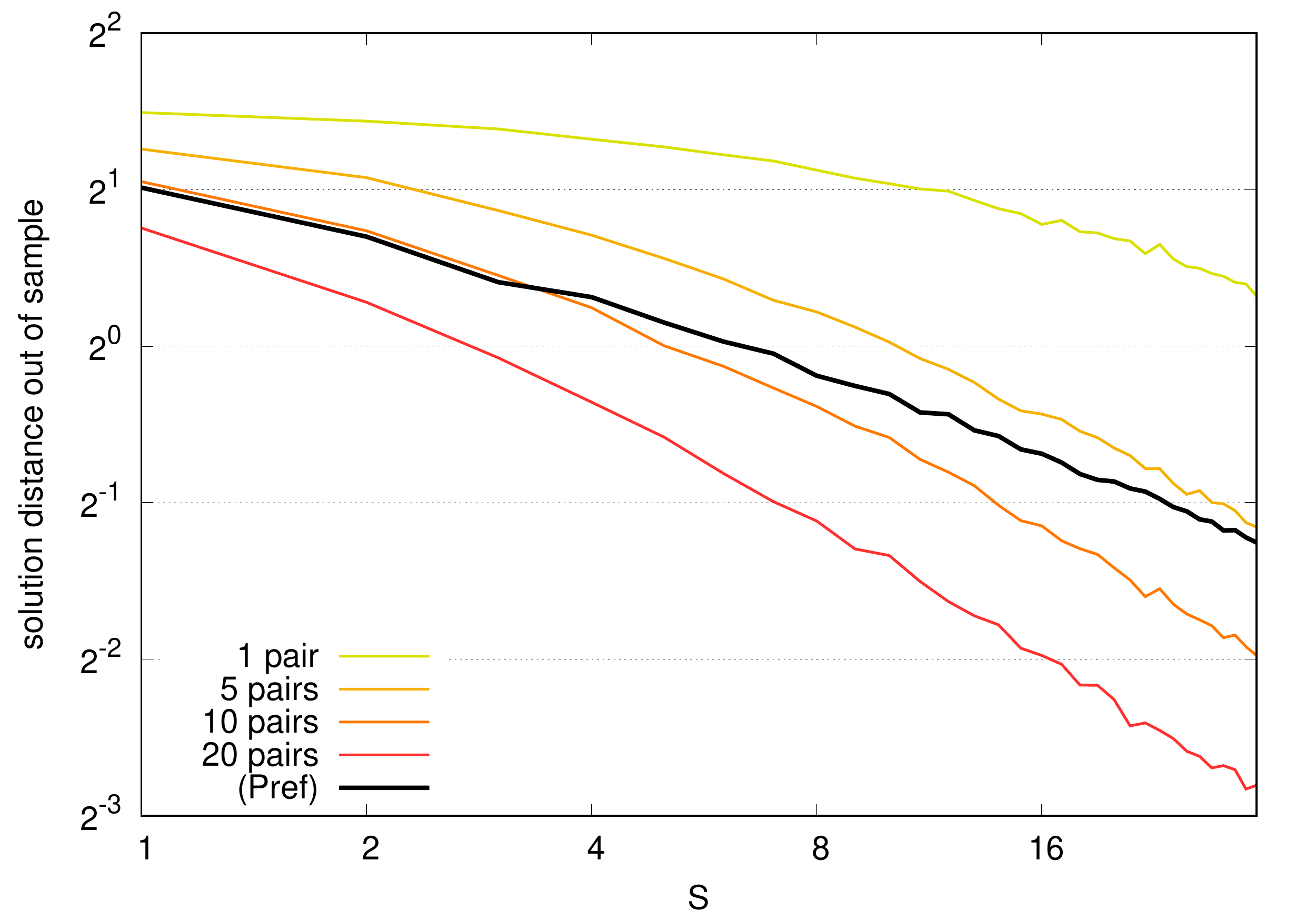}}%
\caption{Comparison of our approach with pairwise comparisons for varying number of  observations $S\in \{1,\ldots,31\}$ for assignment (top) and knapsack (bottom) instances, with $K=5$ and $n=10$ and $n=30$, respectively. \label{fig::Full_S}}
\end{figure}

\begin{figure}[!htb]
\centering
\rotatebox{90}{\hspace{.7cm} Assignment}
\subfigure[\WDist]{\includegraphics[width=.32\textwidth]{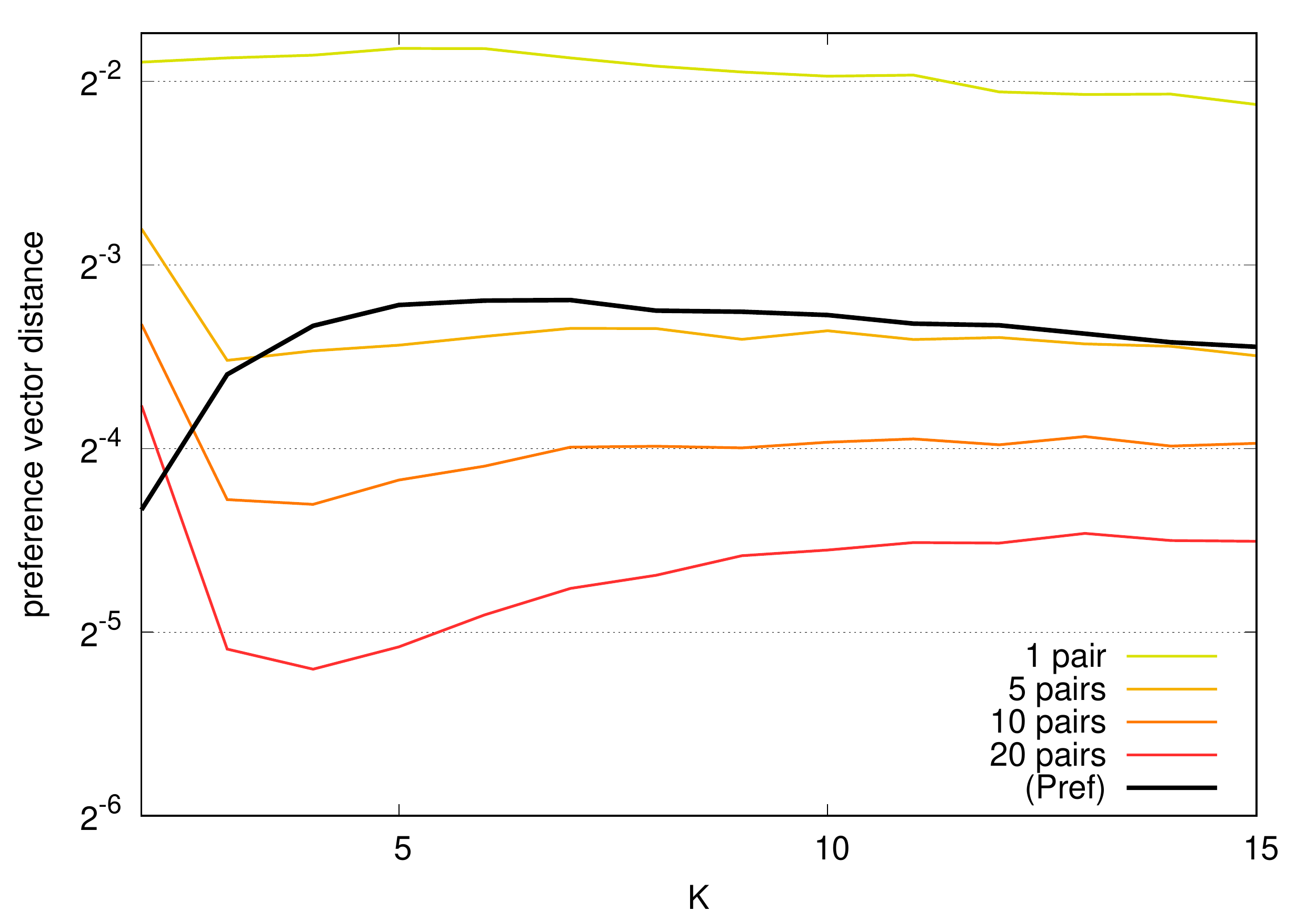}}%
\subfigure[\XInDist]{\includegraphics[width=.32\textwidth]{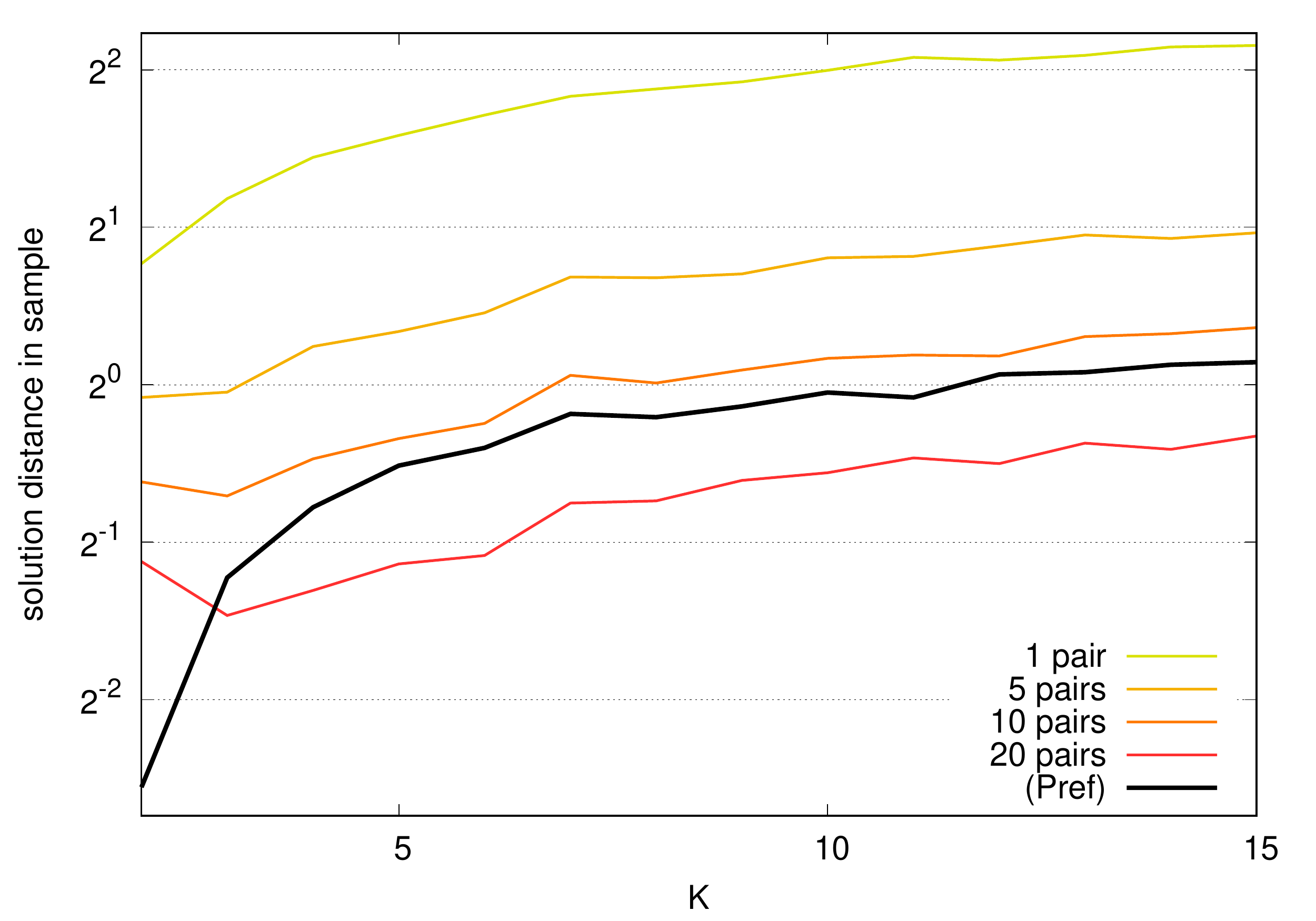}}%
\subfigure[\XOutDist]{\includegraphics[width=.32\textwidth]{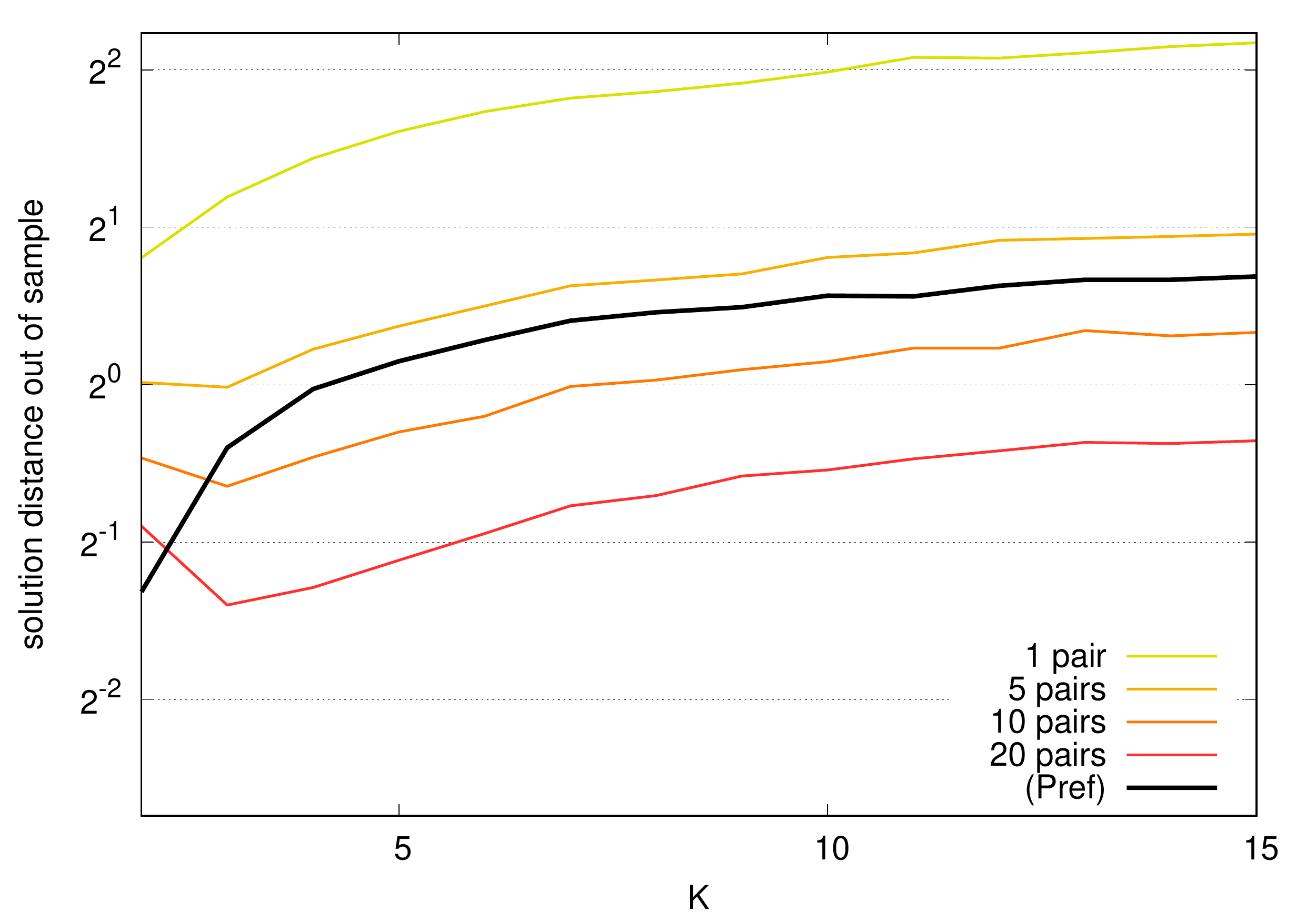}}%

\rotatebox{90}{\hspace{.9cm} Knapsack}
\subfigure[\WDist]{\includegraphics[width=.32\textwidth]{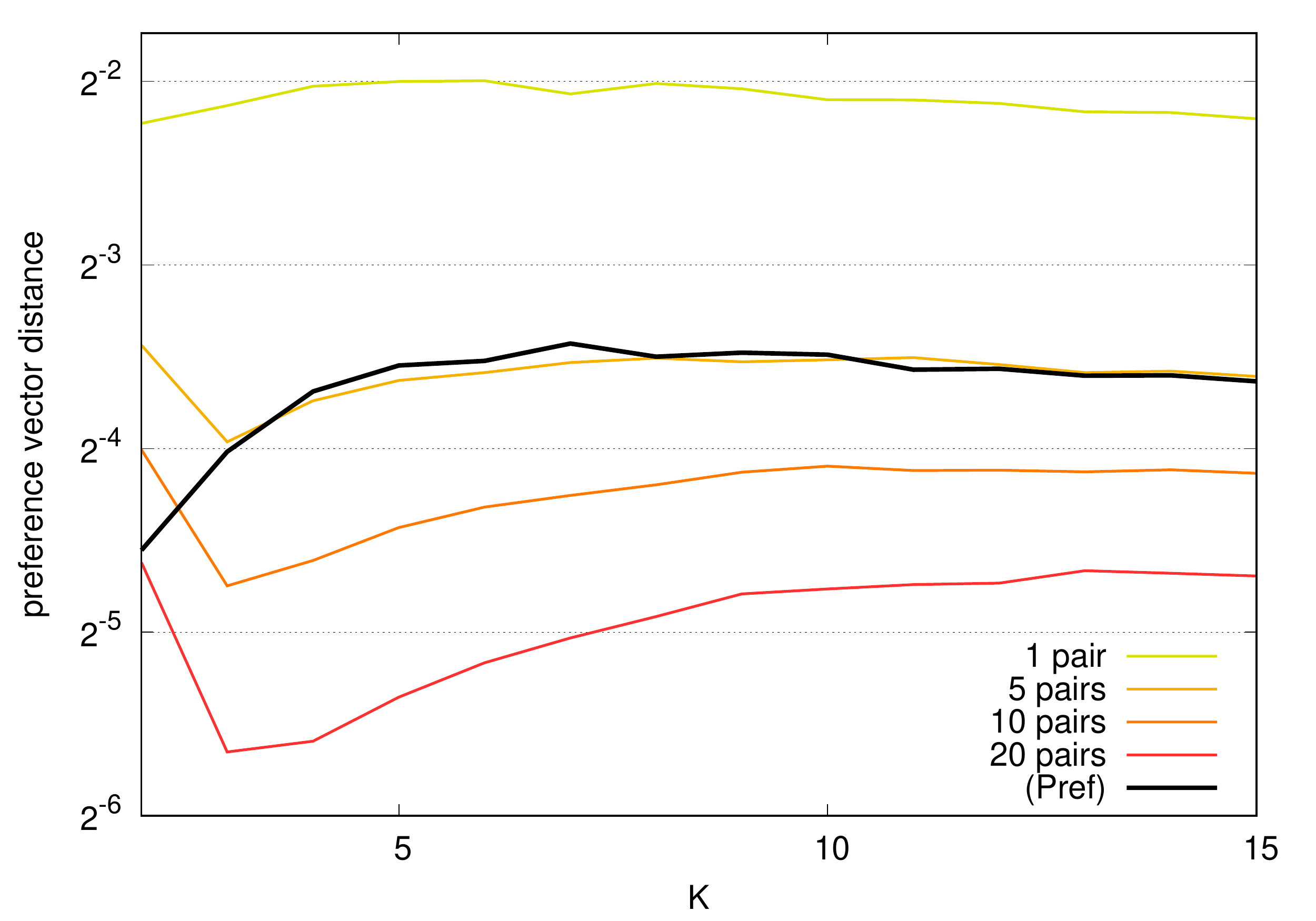}}%
\subfigure[\XInDist]{\includegraphics[width=.32\textwidth]{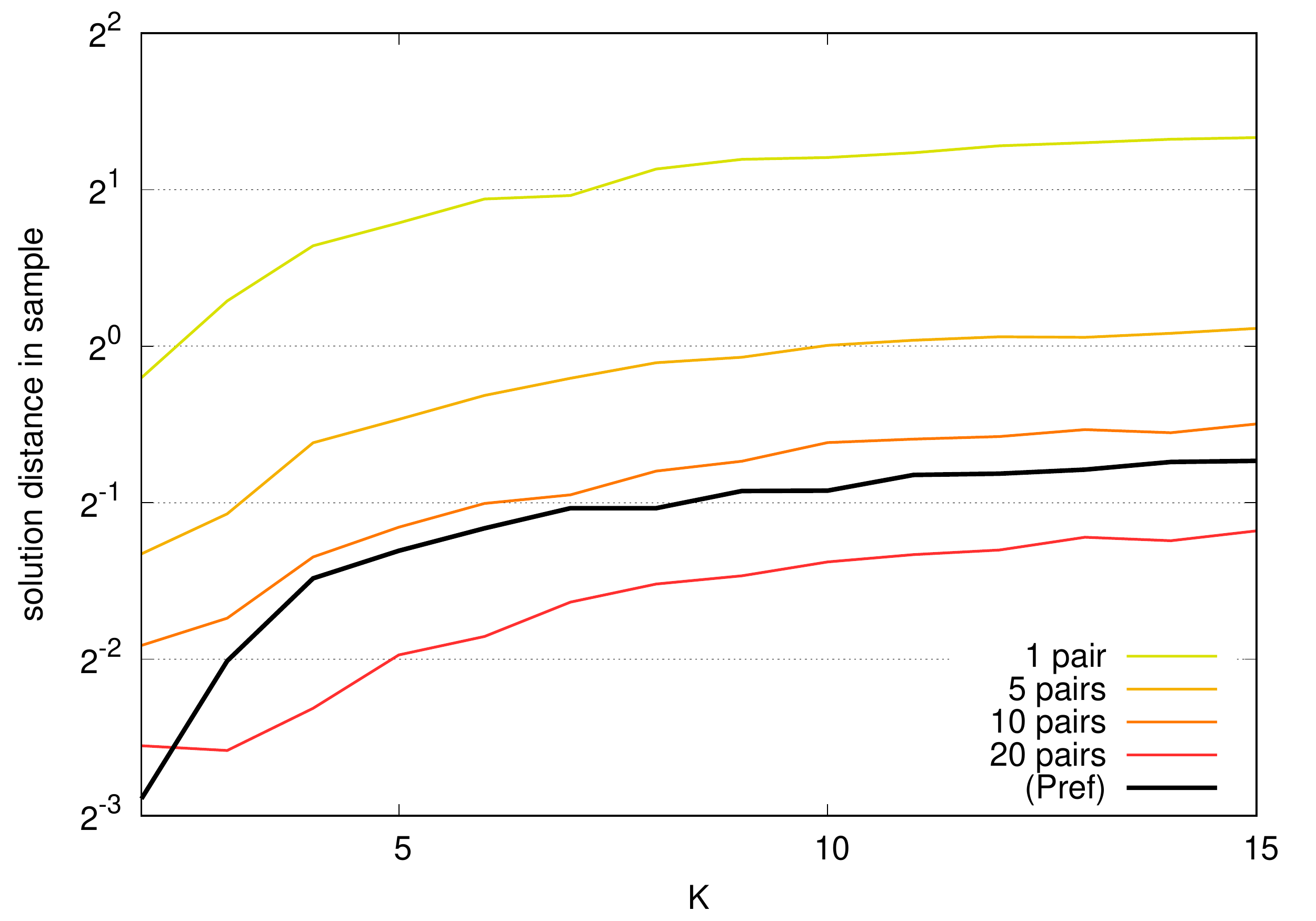}}%
\subfigure[\XOutDist]{\includegraphics[width=.32\textwidth]{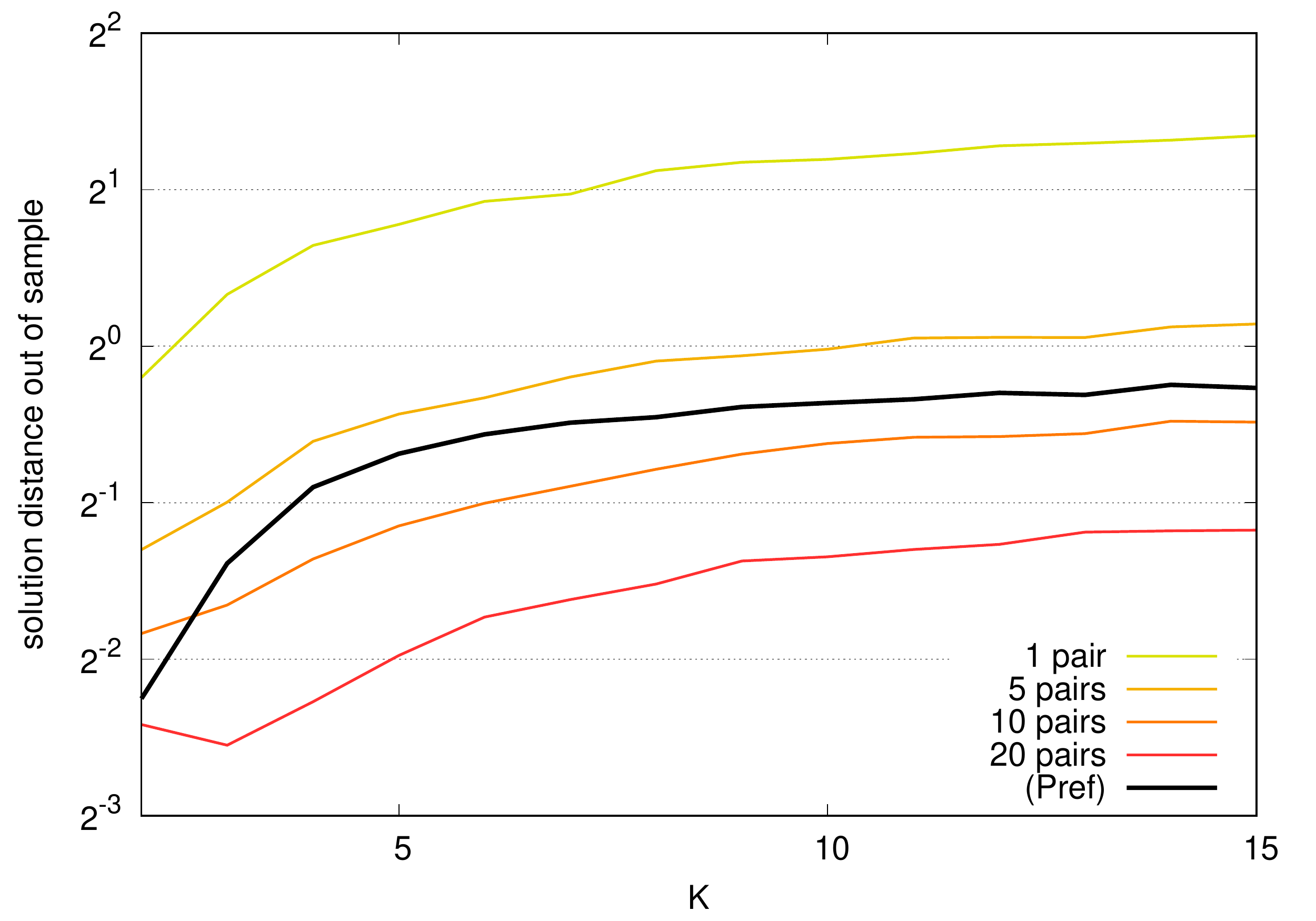}}%
\caption{Comparison of our approach with pairwise comparisons for varying numbers of criteria $K \in\{2,\ldots, 15\}$ for assignment (top) and knapsack (bottom) instances, with $S=16$ and $n=10$ and $n=30$, respectively.\label{fig::Full_K}}
\end{figure}

\begin{figure}[!htb]
\centering
\rotatebox{90}{\hspace{.7cm} Assignment}
\subfigure[\WDist]{\includegraphics[width=.32\textwidth]{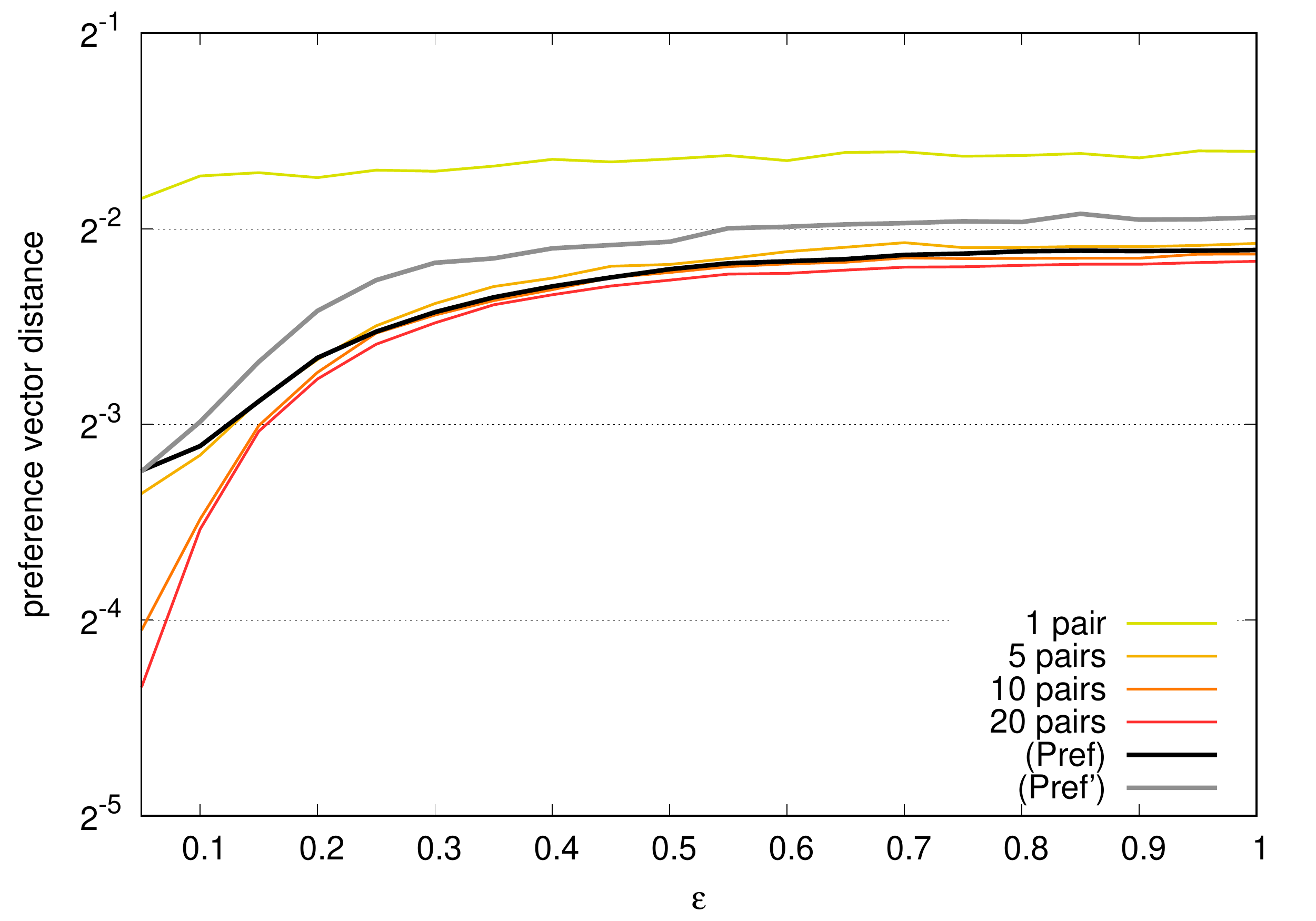}}%
\subfigure[\XInDist]{\includegraphics[width=.32\textwidth]{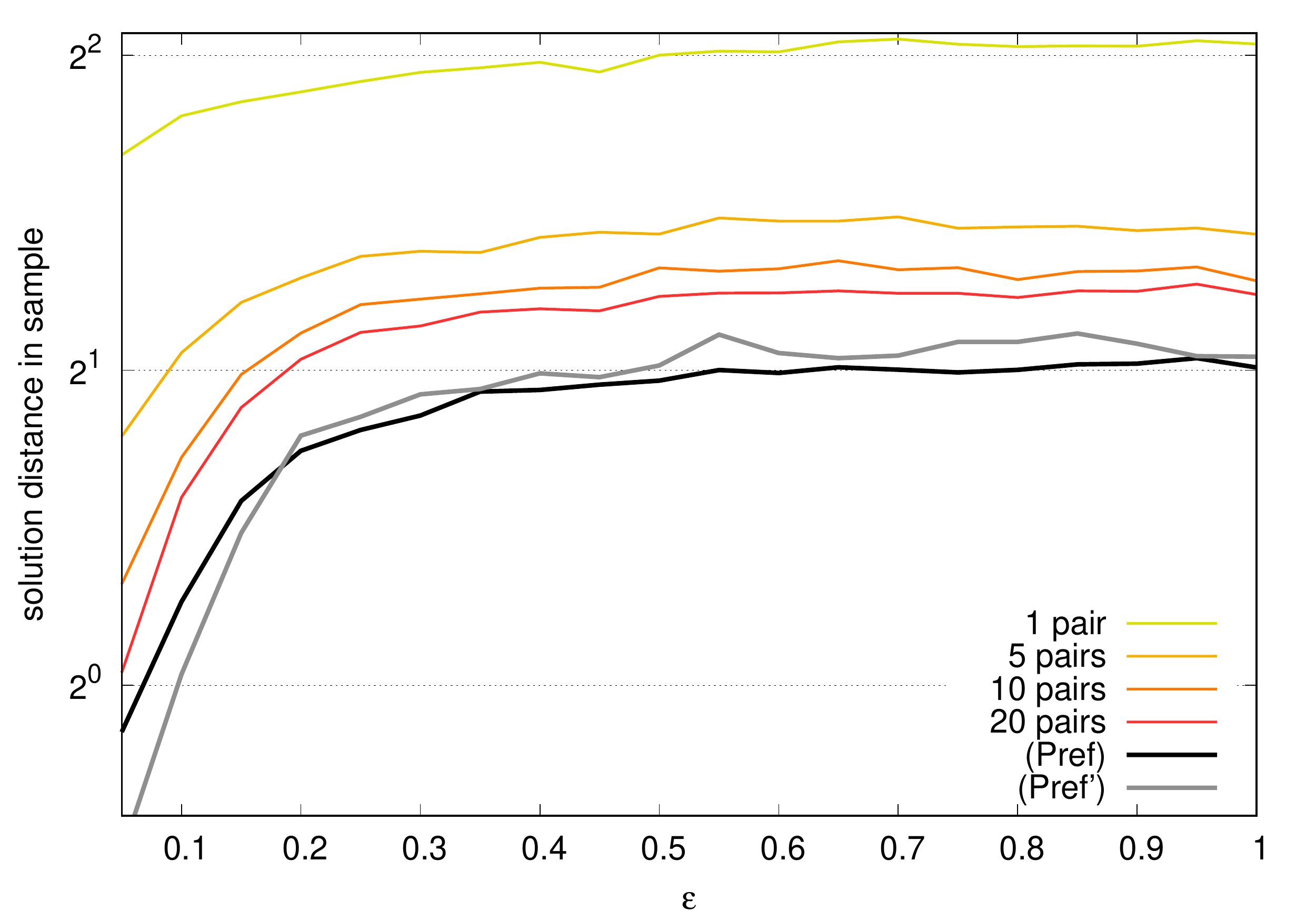}}%
\subfigure[\XOutDist]{\includegraphics[width=.32\textwidth]{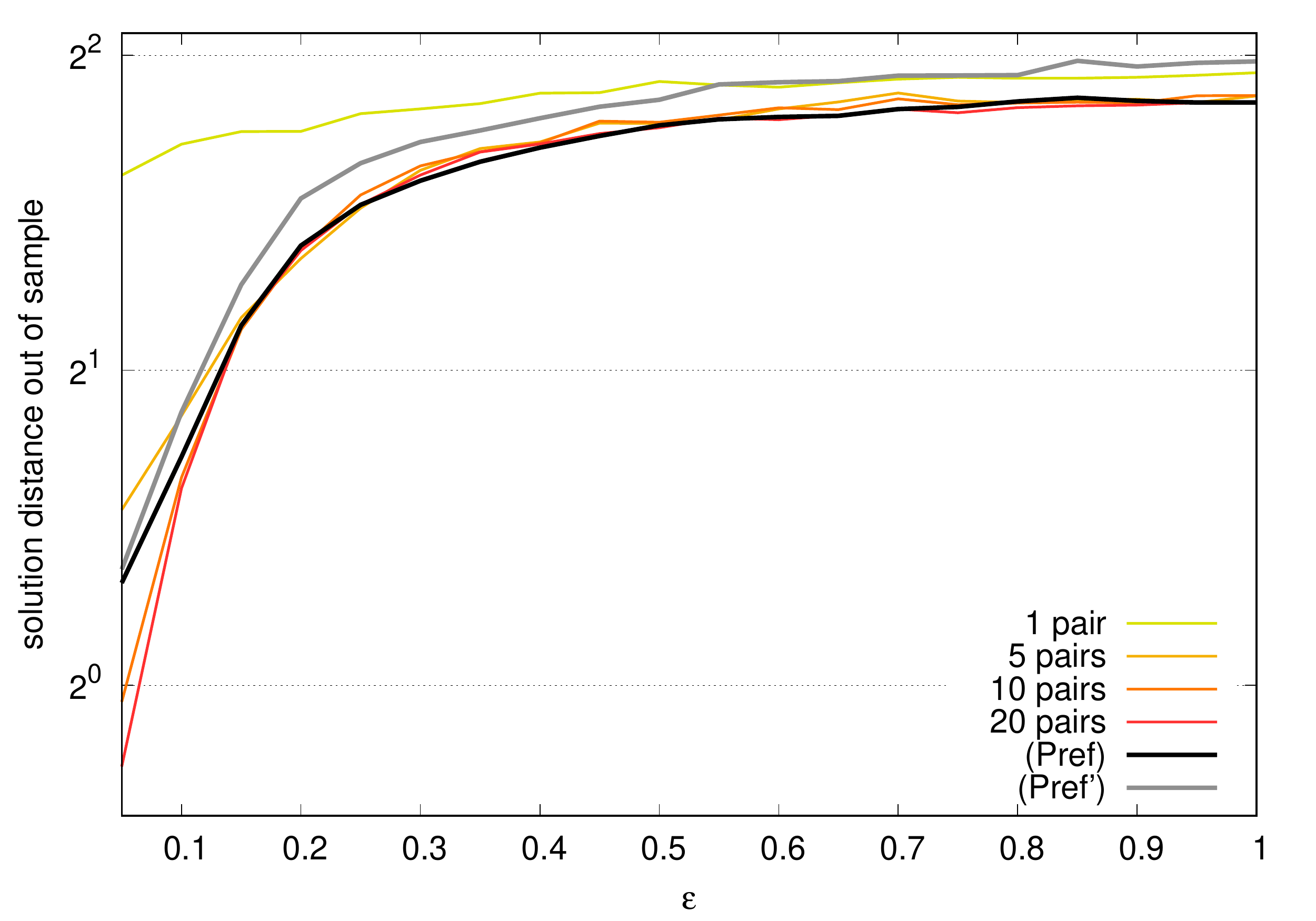}}%

\rotatebox{90}{\hspace{.9cm} Knapsack}
\subfigure[\WDist]{\includegraphics[width=.32\textwidth]{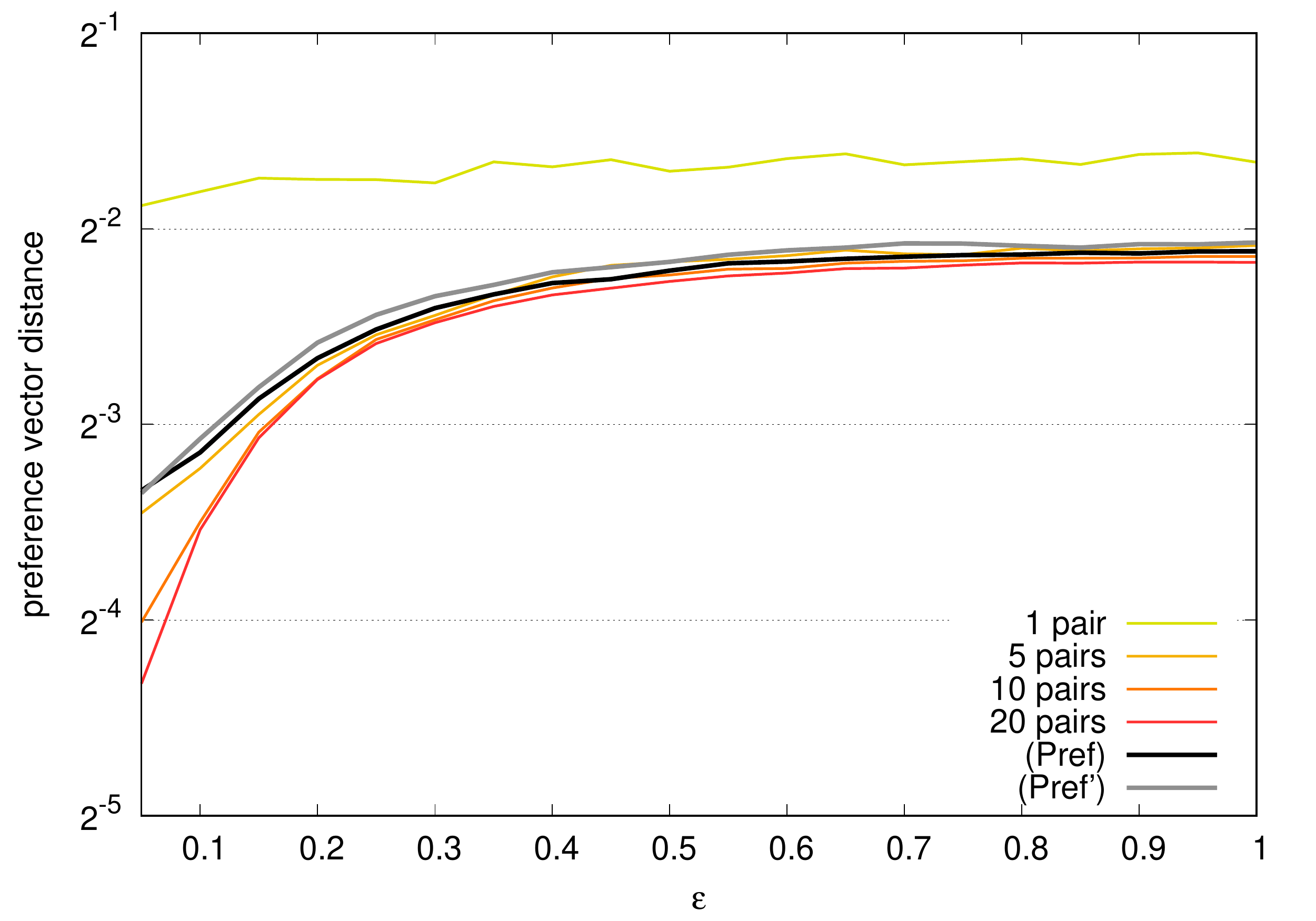}}%
\subfigure[\XInDist]{\includegraphics[width=.32\textwidth]{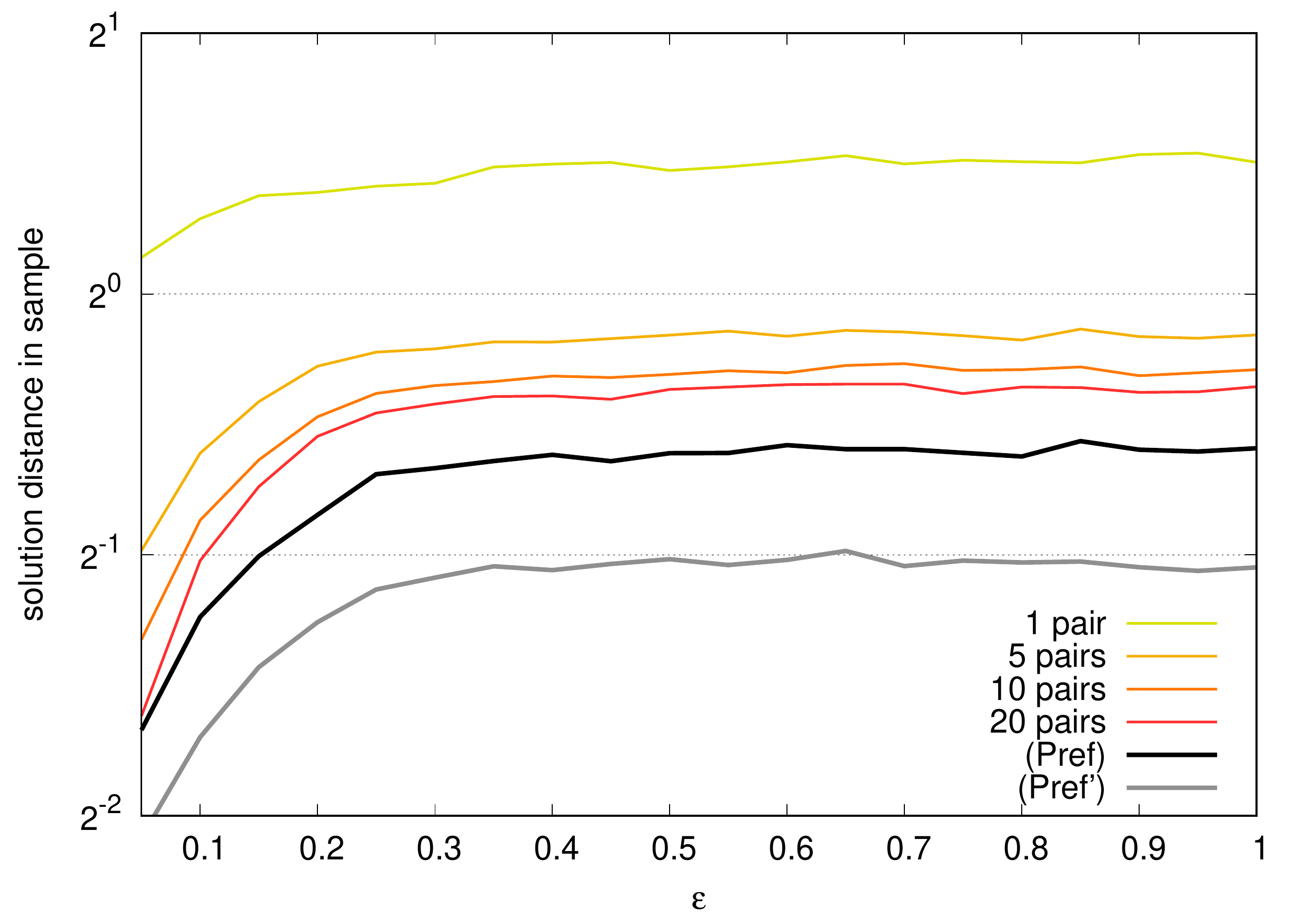}}%
\subfigure[\XOutDist]{\includegraphics[width=.32\textwidth]{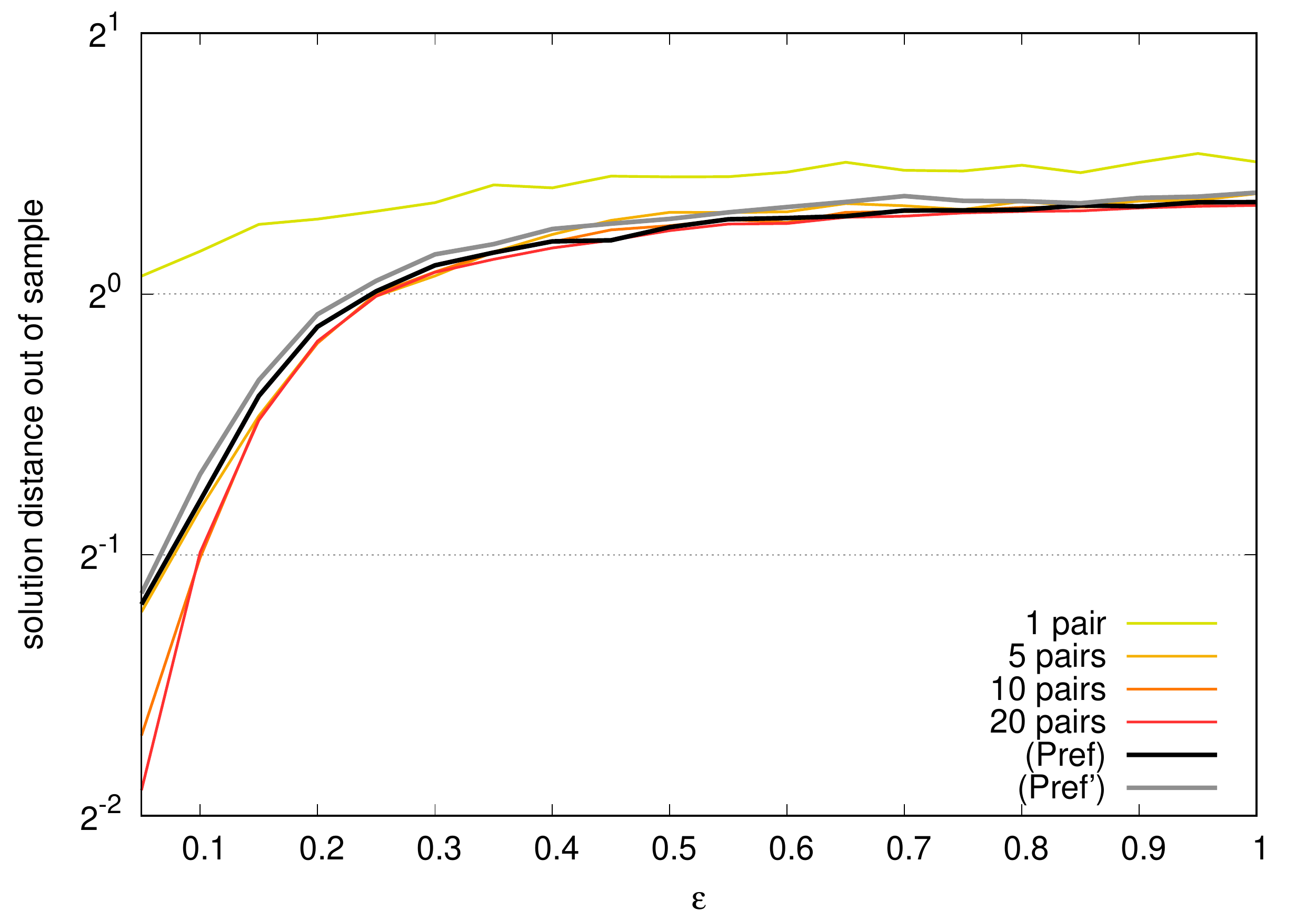}}%
\caption{Comparison of our approach with pairwise comparisons for varying noise $\epsilon$, for assignment (top) and knapsack (bottom) instances, with $S=8$, $K=5$ and $n=5$ and $n=15$, respectively.\label{fig::Full_Eps}}
\end{figure}

We also performed the experiments described in Section \ref{sec:exp} for the assignment problem and the min-knapsack problem, i.e.~$\mathcal{X}=\{\pmb{x}\in\{0,1\}^{n \times n} : \forall j \in [n]:\ \sum_{i\in[n]} x_{i,j}=1 ,\ \forall i \in [n]:\  \sum_{j\in[n]} x_{i,j}=1\} $ and  $\mathcal{X}=\{\pmb{x}\in\{0,1\}^{n} : \sum_{i\in[n]} \omega_i x_{i}\geq \Omega\} $, respectively. Objective function values $C^{s}_{k,i,j}$ and $C^{s}_{k,i}$, respectively, are created in the same fashion as for the selection problem. For the min-knapsack problem the weights $\omega^s_i$ are drawn from a uniform distribution over $[0.7,1.3]$
and the capacities $\Omega^s$ are set to $\frac{1}{2}\sum_{i\in [n]}\omega_i$. In Figures~\ref{fig::Full_N}, \ref{fig::Full_S}, \ref{fig::Full_K}, and \ref{fig::Full_Eps} we show the results for varying values of $n$, $S$, $K$ and $\epsilon$, respectively. The remaining parameters are given in the figure captions.

\end{document}